\documentclass[10pt,a4paper]{article}

\usepackage{comment}
\usepackage{epsfig}
\usepackage{amssymb}
\usepackage{amsmath}
\usepackage{amsmath}
\usepackage{algorithmic}
\usepackage{sectsty}
\usepackage{float}
\usepackage{array}
\usepackage{multirow}
\usepackage{dsfont}
\usepackage{mathrsfs}
\usepackage[numbers]{natbib}
\newtheorem{thm}{Theorem}
\newtheorem{lem}{Lemma}
\newtheorem{algo}[lem]{Algorithm}

\newtheorem{rem}[thm]{Remark}
\usepackage[breaklinks,bookmarks=false]{hyperref}
\hypersetup{
  breaklinks,bookmarks=false,colorlinks = true,
             linkcolor = blue,
             urlcolor  = blue,
             citecolor = blue,
             anchorcolor = blue
}
\usepackage[textsize=tiny,color=yellow!60!green]{todonotes}
\newcommand{\bse}{\begin{subequations}}
\newcommand{\ese}{\end{subequations}}
\usepackage{mathtools}
\usepackage{comment} 
\usepackage[linesnumbered, boxed]{algorithm2e}
\usepackage{enumitem}
\usepackage{capt-of}
\usepackage{graphicx}
\usepackage{textcomp}
\newcommand{\jump}[1]{[\![#1]\!]}

\usepackage{amssymb}

\usepackage[figuresright]{rotating}
\usepackage[labelfont=bf]{caption}

\DeclareSymbolFontAlphabet{\amsmathbb}{AMSb}%
\DeclareMathAlphabet{\mathbfit}{OML}{cmm}{b}{it}

\usepackage{color}
\def\R{{\mathbb R}}
\newcommand{\dsp}{\displaystyle}
\renewcommand{\div}{\nabla\cdot}

 \newcommand{\Hdiv}{\ensuremath{\mathbf{H}(\textnormal{div};\mathrm{\Omega})\;}}
  \newcommand{\RTN}{{\mathbf{RTN}}}
\newcommand{\p}{\partial}

\newcommand{\vf}{{\bf f}}

\newcommand{\vn}{{\bf n}}

\newcommand{\vr}{{\bf r}}

\newcommand{\vu}{{\bf u}}
\newcommand{\vv}{{\bf v}}

\newcommand{\vecv}{\mathbf {v}}
\newcommand{\vecx}{\mathbf {x}}
\graphicspath{{Figures/}}
\usepackage{tabularx}

\newcommand{\vecn}{\mathbf {n}}

\newcommand{\gammaTih}{\mathcal{T}_{h,i}}
\newcommand{\RR}{{\amsmathbb R}}
\newcommand{\textn}[1]{\textnormal{#1}}
\newtheorem{df}[thm]{Definition}

\newcommand{\vK}{\mathbf K}

\newcommand{\divvv}{{\textnormal{div}}}
\DeclareMathOperator*{\argmin}{argmin}
\DeclareMathOperator*{\argmax}{argmax}

\usepackage[breaklinks,bookmarks=false]{hyperref}
\usepackage[textsize=tiny,color=yellow!60!green]{todonotes}
\title{  Splitting-based domain decomposition methods  for two-phase
flow with different rock types\thanks{This work was partially funded by the Hydrinv Inria Euro Med 3+3: HYDRINV project. 
It has also received funding from the EPIC project (within LIRIMA: http://lirima.inria.fr) and the Tunisian Ministry of Higher Education and Scientific Research}}

\author{Elyes Ahmed\footnotemark[2]\ \footnotemark[3] \footnotemark[1]
}
\date{\today}

\begin{document}

\maketitle

\renewcommand{\thefootnote}{\fnsymbol{footnote}}

\footnotetext[2]{Universit\'{e} Tunis El Manar, LR99ES20, LAMSIN-ENIT, B.P. 37, 1002 Tunis-Belv\'{e}d\`{e}re, Tunisia. }
\footnotetext[3]{Inria Paris, 2 rue Simone Iff, 75589 Paris, France.}
\footnotetext[1]{Department of Mathematics, University of Bergen, P. O. Box 7800, N-5020 Bergen, Norway (Current adress).
\href{mailto:elyes.ahmed@uib.no}{elyes.ahmed@uib.no}}

\renewcommand{\thefootnote}{\arabic{footnote}}

\numberwithin{equation}{section}

\begin{abstract}
In this paper, we are concerned with  the global pressure formulation of  immiscible incompressible two-phase flow   
between different rock types. We develop for  this problem  two robust schemes based on  domain decomposition (DD) methods and operator-splitting techniques.  The  first scheme follows a sequential procedure in which the (global) pressure, the saturation-advection and the saturation-diffusion  problems are fully decoupled. In this scheme,  each problem is treated  individually using various DD approaches and specialized numerical methods. The coupling between the different problems is explicit  and the time-marching is with no iterations. To adapt to  different time scales of  problem components and different rock types, the novel scheme uses a multirate time stepping strategy, by taking multiple finer time steps for saturation-advection within one coarse  time step for saturation-diffusion and pressure, and permits independent time steps for the advection step in the different rocks.  In the second scheme,  we  review the  classical Implicit Pressure--Explicit Saturation (IMPES) method (by decoupling only pressure and saturation)  in the context of multirate coupling schemes and nonconforming-in-time DD approaches.  For the discretization,  the   saturation-advection problem is approximated  with the explicit Euler method in time,   and  in space with the cell-centered finite volume method of first order of Godunov type. The saturation-diffusion problem is approximated in time with the implicit Euler method
and in space with the mixed finite element method, as in the pressure problem. Finally, in a series of numerical experiments, we investigate the practicality of the proposed schemes,  the accuracy-in-time of the multirate and nonconforming time strategies, and compare the convergence of various  DD methods within each approach.
\end{abstract}

\vspace{3mm}

\noindent{\bf Key words:}
Two-phase flow in porous
media;  continuous capillary pressure; non-overlapping domain decomposition; operator-splitting;  mixed finite element method; 
 finite volume method; multirate and non-conforming time grids.


\pagestyle{myheadings} \thispagestyle{plain} 

%
\section{Introduction}
Numerical simulation of two-phase flow in  porous formations  has been the subject of investigation of many researchers 
owing to important applications in e.g., the management of petroleum reservoirs or  environmental  remediation. Two-phase flow modeling has received  increased attention in connection with  porous medium  with 
different rock types, so that  the  permeability and the capillary  pressure field are changing across the interfaces
between the rocks. 

Two-phase   flow  in  porous  media  can  be  modeled  by mass  balance  laws  for  each  of  the  fluids \cite{chavent1976new,chen2006computational}.
In particular, an equivalent  formulation can be obtained for the system of equations governing two-phase immiscible, incompressible flow in porous
media,  by introducing an artificial variable called the \textit{global pressure}. The dependent variables in such formulation  are the global  pressure and
 the wetting  phase saturation \cite{chavent1976new}. Considering these variables, the governing equations consist 
 of a \textit{nonlinear} elliptic Darcy's law for the global pressure and a \textit{nonlinear  degenerate} parabolic equation of \textit{advection-diffusion} type for  the wetting  phase saturation.
 Flow simulation of such systems in porous media with \textit{different  rock types} is very difficult because of  variations in physical parameters in  the different rock types 
 that require coupling~\cite{alboin2000domain,chavent1986mathematical}. Particularly, due to changes  in permeability and   capillary force,
the global pressure and saturation  exhibit \textit{strong discontinuities} across the interfaces 
between the rocks (cf. \cite{MR2989844,MR2997427}).  
Several numerical schemes have been developed for two-phase flow
in heterogeneous media. Finite volume schemes have been proposed in \cite{MR2465972,MR2206441} for a simplified two-phase model involving only diffusion effects. 
One can see~\cite{MR2559741} for the extension to the full two-phase flow between different rock types. A discontinuous Galerkin (DG) method or    a combination of  DG method and mixed finite elements (MFE) 
method was employed in \cite{MR3142429,MR2989844,hoteit2008numerical}. For domain decomposition methods applied to two-phase flow problems, one can cite in particular  
the work in~\cite{ganis2014global}, where  a  fully implicit-multiscale mortar mixed finite element method 
for two-phase flow in a heterogeneous porous medium is presented. 
In~\cite{ahmed:hal-02275690}, a space--time domain decomposition method formulated using Robin and Ventcell 
type coupling conditions was applied to a simple two-phase flow model between different rock types (cf.~\cite{ahmed:hal-01540956,ahmed2018global}). 
We  also refer the reader to linear domain decomposition methods 
for two-phase flow problems in homogeneous porous media presented in~\cite{MR3771899,seus2017lineartwophase}. Therein, the problem of two-phase flow   is  decomposed into a set of subproblems  in different subdomains, and solved   at each time step with  the 
$L$-scheme linearization method;  a  robust  quasi-Newton  method with a parameter $L >0$ mimicking the Jacobian from the Newton method~\cite{MR3489128,MR3926838,MR2079503,MR3800042}. Other related works can be found in \cite{ahmed2018multiscale,ahmed2019robust,MR3144798,HOANG2017366,peszynska2002mortar,reichenberger2006mixed}.

Our  contribution in this paper is to formulate   
robust  \textit{domain decomposition methods}  for the two-phase flow 
between  \textit{different rock types}. The methods are motivated  by the fact that the (global) pressure and saturation as well as 
advection and diffusion effects  act on \textit{different  time scales} within the same rock type and also between different rock types, 
and it is easy to implement this in the domain decomposition context. However,  classical domain  
decomposition  methods  are known to  perform  poorly  if  diffusion  is  dominated  by  advection as well as in purely advective problems~(cf.~\cite{MR1370250} for more details). From these observations, we  investigate    
new  procedures based on a divide-and-conquer (operator-splitting) strategy for the numerical solution of  the two-phase flow  between  different rock types~\cite{MR3564686,brun2019monolithic,hoteit2008efficient,MR2837398}. In the first approach,  we  split   the global pressure,   the saturation-advection and the saturation-diffusion problems  and  treat them individually using specialized numerical methods (see Scheme~\ref{flowchart} for the workflow strategy). This workflow will then apply various domain decomposition approaches adapted to each problem, transforming  the  original problem  into a sequential solves of reduced problems of~\textit{pressure,  saturation-diffusion} and \textit{saturation-diffusion} posed only on the physical interfaces between rock types,   through the use of trace operators. The DD approaches are based either on the Steklov-Poincar\'{e} type methods~\cite{MR1297465,MR3144798,MR1208381} or  the optimized Schwarz methods~\cite{MR2218966,MR2344706,HOANG2017366}. 

The resulting interface problems of the first approach can  be solved by various iterative methods. For the  discretization in the subdomains,  the saturation-diffusion is approximated  with the implicit Euler method in time, and in space with the \textit{hybrid mixed finite element method}, as in the pressure problem.   To handle efficiently  the advection step, a numerical scheme is presented that uses  the \textit{hybrid  finite volume method}
 based on the \textit{Godunov scheme} (cf. \cite{MR2051062,MR3208750,Mishra2017479}), with the explicit Euler method in time.  The coupling strategy
 between the three processes is explicit and  the time marching is without any iterations. It is known that the time step must be sufficiently small for the stability of such  explicit-coupling schemes~\cite{Douglas1983}. To overcome this difficulty and reveal the multi-physical process  of the  problem, the algorithm  allows  a multirate time stepping strategy, in which  adequate time steps are used for each of the  problems~(see Scheme~\ref{flowchart}). Here, we consider   applications,  where  advection effects are often much greater than capillary diffusion effects. To reduce the computational cost (i) we solve the (global) pressure and  saturation-diffusion problems in the coarse time step and solve the advection problem in the finer time step (ii) we  employ  \textit{nonconforming time grids} between different rock types for the advection step. For an overview and further insight on multirate schemes, we refer the reader to~\cite{MR3564686,MR3771343,MR3351781} where various coupled problems are considered (see also the adaptive strategy in~\cite{MR3983155}). For  local/nonconforming time discretizations in the context of DD approaches, we refer to~\cite{MR3144798,HOANG2017366} and the references therein.

 In the  second approach of this paper, we  review  the classical Implicit Pressure--Explicit Saturation (IMPES) method (cf. \cite{2017arXiv170901644L}) in the context of DD methods. This approach consists on decoupling only  (global) pressure  and   saturation  problems    in  order  to  solve   \textit{sequentially}  each  problem   for  its  main  variable.  The pressure problem is discretized as in the first approach with MFE method, and the coupled advection-diffusion problem for the saturation,  is discretized with a scheme \textit{explicit} in time and   with a cell-centered finite volume method  in space. In this approach,  we employ  \textit{larger time steps} for the pressure problem with respect to the  saturation problem, and to further save computational effort, we use  \textit{different time-steps}  for the saturation  in  the different rock types. 
 
 The  two  approaches of this paper are iteration-free and allow  easy \textit{reuse of existing codes}  specialized to each component of the problem,  which minimizes the computational effort. Furthermore,   the approaches  can easily \textit{integrate}  more \textit{complex problems};    
 the    first approach is \textit{extended}  to address \textit{reduced fracture model}  between different rock types presented in~\cite{Ahmed2016}.

 The remainder of this paper is organized as follows. In Section~\ref{flow_rock_types}, we recall briefly  the global pressure formulation for incompressible 
two-phase flow problem. We then rewrite this  problem for the case of flow between different rock types and describe  the physical conditions occurring at the interface between the rocks. In~Section~\ref{DD_FVMF}, we present the  numerical approaches.  We show numerical results in Section~\ref{section:numerics}.  Section~\ref{section:conclusion} contains concluding remarks. In Appendix~\ref{section:appendix}, we present  
 a brief   application of the first approach on a  reduced fracture model  for two-phase flow between different rock types.
\label{intro}
\section{Two-phase Darcy
flow model between different rock types}\label{flow_rock_types}
Let $\Omega$ be an open bounded domain of $\RR^{d}$ ($d= 2$ or $3$) which is
assumed to be polygonal if $d=2$ and polyhedral if $d=3$. We denote by
$\partial \Omega$ its boundary (supposed to be Lipschitz-continuous) and by
$\vn$ the unit normal to $\partial \Omega$, outward to $\Omega$.
\subsection{Global pressure formulation for incompressible two-phase flow}
We consider the immiscible incompressible two-phase flow in a porous medium. Assuming that
there are only two phases occupying the porous medium $\mathrm{\Omega}$, and that
each phase is composed of a single component, the mathematical form of this problem is as follows~\cite{chavent1986mathematical}:
\bse\label{Physical_problem}\begin{alignat}{3}
\label{conservation-ellph1}
&\mathrm{\Phi} \displaystyle{\frac{\p s_{\ell}}{\p t}} + \div\left(- \vK k_{\ell}(s_{\ell})(\nabla p_{\ell}-\rho_{\ell} \vu_{g})\right) = 0,\quad&&\mbox{in }\mathrm{\Omega}\times(0,T), \\
\label{Volume_pore}
&s_{w}+s_{n}=1,\quad&&\mbox{in }\mathrm{\Omega}\times(0,T),\\
\label{Cap_press_eq}
&p_{n}-p_w=\pi(s_w),\quad&&\mbox{in }\mathrm{\Omega}\times(0,T),
\end{alignat}\ese
where the unknowns are $s_{\ell}$, the \textit{phase saturation}, and $p_{\ell}$, the \textit{phase pressure},  $\ell\in\{w,n\}$. The subscripts $w$,
$n$ stand for wetting and non-wetting, respectively. Typically, the non-wetting  phase is oil and the wetting one is water.    The  function $\mathrm{\Phi}=\mathrm{\Phi}(\vecx)$ denotes the porosity of the rock ($\mathrm{\Phi}\in(0,1)$ in the domain $\mathrm{\Omega}$) and  $\vu_{g}$ 
denotes the \textit{gravity field}; the permeability of the porous medium $\vK$ is assumed to be a positive scalar function,  
the \textit{mobility} $ k_{\ell}$   of the phase $\ell$ is an increasing function of the saturation $s_{\ell}$, 
satisfying $ k_{\ell}(0)=0$ and $ k_{\ell}(1)=1$, and $\rho_{\ell}$ denotes the density of phase $\ell$, $\ell\in\{w,n\}$.  
 The function $\pi$ is the \textit{capillary pressure},   defined  as  the  difference  between  the phase pressures  and assumed to be a function (positive and decreasing) of the saturation $s_{w}$. 
More details on capillary pressure and relative permeability can be found in \cite{chavent1986mathematical,van1995effect}.
We assume  that  the  initial  phase  saturation  is
known, say 
\begin{equation}\label{Initial_cond}
 s_{w}(\cdot,0)=s_{0},\,\, \textnormal{ in }\mathrm{\Omega}.
\end{equation}
 We also assume homogeneous Neumann boundary conditions on the phase fluxes, i.e.,
\begin{eqnarray}
\label{B_cd}&&\vK\, k_{\ell}(s_{\ell})(\nabla p_{\ell}-\rho_{\ell} \vu_{g})\cdot\vecn=0,\,\, \textnormal{ on }\partial\mathrm{\Omega}\times(0,T).
\end{eqnarray}

In order to avoid some of the known difficulties related to the degeneracy of the problem~\eqref{Physical_problem}--\eqref{B_cd},  we follow the classical idea of
introducing the so-called global pressure formulation \cite{chavent1976new}. For the simplicity of notation let $s:=s_{w}$, so that $s_{n}:=1-s$. 
Define the total mobility function $\mathbf{M}(s):=\vK\left(k_w(s)+k_{n}(s)\right)$ and introduce  the fractional flow function 
$f(s):={\displaystyle \frac{k_w(s)}{k_w(s)+k_{n}(s)}}$ and let $\beta(s):=\displaystyle{\int_1^s} (\frac{1}{2}-f(a))\pi'(a) 
\textnormal{d}a,\, \forall s\in(0,1)$. We then introduce
the global pressure function $p$ by
\begin{equation}
\label{globalpressureeq}
  p := \frac{1}{2}(p_w + p_{n}) + \beta(s),
\end{equation}
and the Kirchoff transform
\begin{eqnarray*}
\alpha(s):= \dsp{\int_0^s} -\displaystyle\frac{k_w(a)k_{n}(a)}{k_w(a)+k_{n}(a)}\pi'(a) \textnormal{d}a,\,\quad \forall s\in(0,1).\\ 
\end{eqnarray*}
Following \cite{chavent1976new},  the problem~\eqref{Physical_problem}--\eqref{B_cd} can be  rewritten  under the form
\bse\label{mathematical_problem}\begin{alignat}{3}
\label{System_Conserv_wett}&\mathrm{\Phi} \dfrac{\p {s}}{\p t} + \div \vu_{w} = 0,\quad&&\mbox{in }\mathrm{\Omega}\times(0,T),\\
\label{System_Darcy_wett}&\vu_{w}=f(s)\vu +  f_{g}(s) \vu_{g}- \vK\nabla \alpha(s),\quad&&\mbox{in }\mathrm{\Omega}\times(0,T),\\
\label{System_Conserv_tota}&\div\vu =0,\quad&&\mbox{in }\mathrm{\Omega}\times(0,T),\\
\label{System_Darcy_total}&\vu=-\mathbf{M}(s)\left(\nabla p-\rho(s)\vu_{g}\right),\quad&&\mbox{in }\mathrm{\Omega}\times(0,T),\\
&\vu_{w}\cdot\vecn=0,\;\;\textnormal{and}\;\;\vu_{i}\cdot\vecn=0,\quad&&\textnormal{on } \partial\Omega\times(0,T),\\
&s(\cdot,0)=s^{0},\quad&&\textn{in }\mathrm{\Omega},
\end{alignat}\ese
where
\begin{equation}
f_{g}(s):=\vK\displaystyle{\frac{k_w(s)k_{n}(s)}{k_w(s)+k_{n}(s)}
(\rho_w-\rho_{n})},\quad\textnormal{and }\quad
\rho(s):=\dfrac{ k_w(s)\rho_w+k_{n}(s)\rho_{n}}{ k_w(s)+k_{n}(s)}.
\end{equation}
In view of the above formulation, in equation~\eqref{System_Darcy_total} the saturation appears  only through the coefficients $\mathbf{M}(s)$ and $\rho(s)$, 
and not through its gradient  $\nabla s $. Hence, the equation~\eqref{System_Darcy_total} looks very much like  \textit{Darcy's law} 
for a fictitious fluid representing the
global (water + oil) flow pattern. The global pressure is not a physical pressure but a good numerical and methematical unknown~\cite{chavent1976new,chavent1986mathematical}. Precisely, it is a smooth function defined in the whole domain, 
whether a phase vanishes or not, while this is not the case for the phase pressures. Still following~\cite{chavent1976new}, the global pressure satisfies  $p_{w} \leq p \leq p_{n}$, and such that $p=p_{w}$ if $s_{w}=1$ and $p=p_{n}$ if $s_{w}=0$. Also, by  plugging equation~\eqref{System_Darcy_wett} into~\eqref{System_Conserv_wett}, the resulting \textit{saturation equation} gives a \textit{nonlinear
parabolic equation} of \textit{advection-diffusion} type given by  the sum of 
a diffusion contribution  due to  capillary effects and an advection contribution due to gravity  and  total flow effects.

\subsection{Flow between different rock types}
In this part,  we particularize the model problem~\eqref{mathematical_problem} to a porous medium with different rock types, following~\cite{alboinDD}. Precisely, we suppose that $\mathrm{\Omega}$ is divided  into two non-overlapping polyhedral subdomains
$\mathrm{\Omega}_{i}$, $ i\in\{1,2\}$. The exterior boundary of $\mathrm{\Omega}_{i}$ is denoted by $\mathrm{\Gamma}_{i}^{\textnormal{N}}:=\partial \mathrm{\Omega}_{i}\cap\partial\mathrm{\Omega}$, 
$ i\in\{1,2\}$. Let $\mathrm{\Gamma}$ be the interface between the subdomains, i.e., 
$\mathrm{\Gamma}:=\partial \mathrm{\Omega}_{1}\cap\partial\mathrm{\Omega}_{2}$ (see Figure~\ref{cap_rock_types} (left)).    Let $\vn_{i}$    denotes the unit outward normal on $\partial \mathrm{\Omega}_i$.  In the following, we consider that each subdomain $\mathrm{\Omega}_{i} $ corresponds to a (homogeneous) rock type.   This means that across $\Gamma$, not only  the porosity and the absolute permeability differ, but also the \textit{relative permeability} 
and \textit{capillary pressure} curves. For every  function~in~\eqref{mathematical_problem}  involving  a  space-dependent quantity   (e.g., $f(s)$)  that  is  valid  separately  in  each  of  the  two  subdomains (see Figure~\ref{cap_rock_types} (left)), we will use subscripts (e.g., $f_{i}(s_{i})$).
\begin{figure}[hbtp]\centering
\includegraphics[scale=0.75]{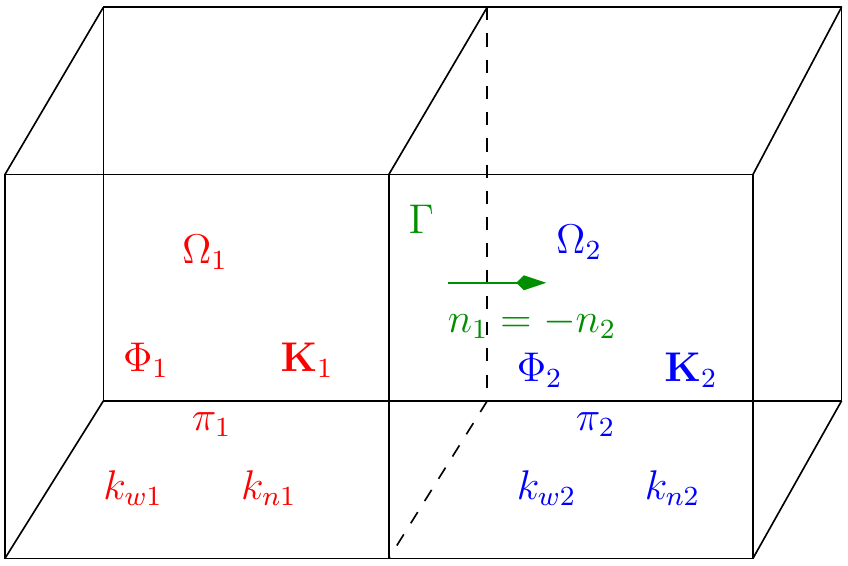}
\includegraphics[scale=0.3]{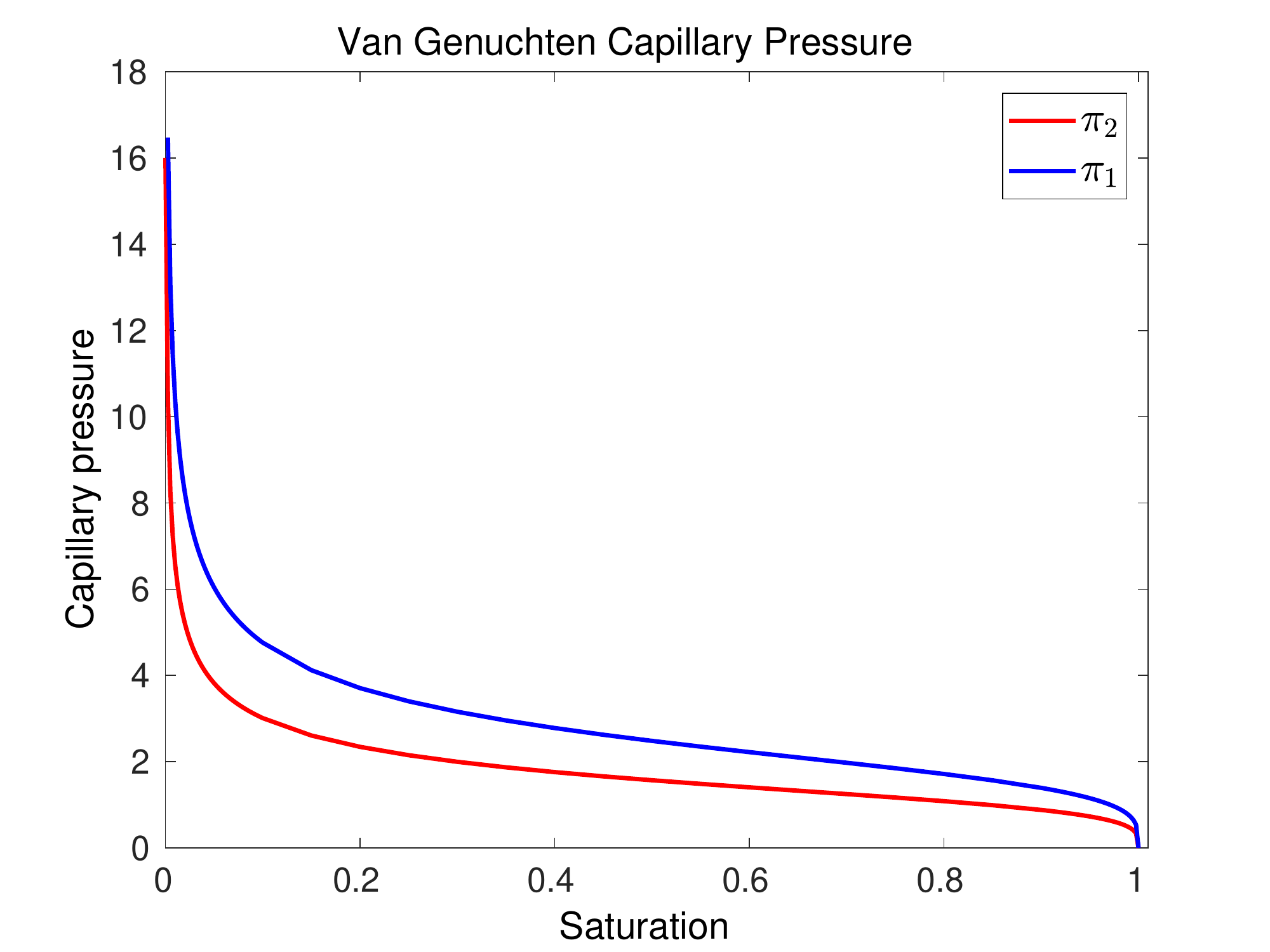}
\caption{Illustration of the computational domain with different   rock types (left) and  associated capillary pressures curves (right).}
\label{cap_rock_types}
\end{figure}

Now that the decomposition of $\Omega$ is made, the two-phase flow equations~\eqref{mathematical_problem} remain valid in each  subdomain $\mathrm{\Omega}_{i}$, $ i\in\{1,2\}$, i.e.,
\bse\label{mathematical_problem_DD}\begin{alignat}{3} 
\label{Chap2_System_Conserv_wett}&\mathrm{\Phi}_{i} \dfrac{\p {s_{i}}}{\p t} + \div \vu_{wi} =  0,\quad&&\textnormal{in }\mathrm{\Omega}_{i}\times(0,T),\\
\label{Chap2_System_Darcy_wett}&\vu_{wi}=f_{i}(s_{i})\vu_{i} +  f_{gi}(s_{i})\vu_{g}- \vK_{i}\nabla \alpha_{i}(s_{i}),\quad&&\textnormal{in }\mathrm{\Omega}_{i}\times(0,T),\\
\label{Chap2_System_Conserv_tota}&\div\vu_{i} =0,\quad&&\textnormal{in }\mathrm{\Omega}_{i}\times(0,T), \\
\label{Chap2_System_Darcy_total}&\vu_{i}=-\mathbf{M}_{i}(s_i)\left(\nabla p_{i}-\rho_{i}(s_{i})\vu_{g}\right),\quad&&\textnormal{in }\mathrm{\Omega}_{i}\times(0,T),\\
&\vu_{wi}\cdot\vecn_{i}=0,\;\;\textnormal{and}\;\;\vu_{i}\cdot\vecn_{i}=0,\quad&&\textnormal{on } \mathrm{\Gamma}_{i}^{\textnormal{N}}\times(0,T),\\
&s_{i}(\cdot,0)=s^{0},\quad&&\textn{in }\mathrm{\Omega}_{i},
\end{alignat}\ese
where  the unknown functions for the saturation equations~\eqref{Chap2_System_Conserv_wett}-\eqref{Chap2_System_Darcy_wett} are $s_{i}$, the subdomain saturation,  and $\vu_{wi}$, 
the subdomain velocity for the wetting phase, $ i\in\{1,2\}$. That of the global pressure equations~\eqref{Chap2_System_Conserv_tota}-\eqref{Chap2_System_Darcy_total} are  $p_{i}$,
 the subdomain  (global) pressure,  and $\vu_{i}$,  the subdomain total velocity, $ i\in\{1,2\}$.

Now we come to the \textit{transmission conditions} across the  interface $\mathrm{\Gamma}$.
  The subdomain saturation problems~\eqref{Chap2_System_Conserv_wett}-\eqref{Chap2_System_Darcy_wett}, for $i\in\{1,2\},$ are coupled through the following interface conditions~\cite{alboinDD}: 
\bse\label{Chap2_Matching_saturation}
\begin{alignat}{2}
 \label{Chap1_Matching_capilarity}&\jump{\pi(s)}=0,&\quad\textnormal{on }\mathrm{\Gamma}\times(0,T),\\
  \label{Chap2_Matching_wettingflow}&\jump{\vu_{w}\cdot \vn}=0,&\quad\textnormal{on }\mathrm{\Gamma}\times(0,T),
\end{alignat}\ese
where the jumps on  $\mathrm{\Gamma}$  are defined as
\begin{alignat*}{1}
\jump{v}:=v_{1}-v_{2},\textnormal{ and } \jump{\vv \cdot \vn}:=\vv_{1} \cdot \vn_{1}+\vv_{2} \cdot \vn_{2}.
\end{alignat*}

The  condition on the capillary pressure~\eqref{Chap1_Matching_capilarity} stems from the continuity  of  phase pressures since $\pi(s)=p_{n}-p_{w}$, resulting in a \textit{jump} on the  saturation as depicted in Figure~\ref{cap_rock_types} (right). The  condition~\eqref{Chap2_Matching_wettingflow}   imposes the continuity of  the normal component of the wetting flow $\vu_{w}$,  to preserve phase conservation.   The subdomain pressure problems~\eqref{Chap2_System_Conserv_tota}-\eqref{Chap2_System_Darcy_total}, for $i\in\{1,2\},$ are coupled 
through the following interface conditions~\cite{Ahmed2016}:
\bse\label{Chap1_Matching_pressure}\begin{alignat}{2} 
\label{Chap1_Matching_pressure1}&\jump{p-\beta(s)}=0,&\quad\textnormal{ on }\mathrm{\Gamma}\times(0,T),\\
\label{Chap1_Matching_pressure2}&\jump{\vu \cdot \vn}=0,&\quad \textnormal{ on }\mathrm{\Gamma}\times(0,T).
\end{alignat}\ese
Equation (\ref{Chap1_Matching_pressure1})  justifies  the \textit{continuity} of the \textit{phase pressures} since $p - \beta(s)=\frac{1}{2}(p_w + p_{n})$. 
Equation (\ref{Chap1_Matching_pressure2}) represents \textit{conservation} of the \textit{total fluid}.
\begin{rem}[Compatibility condition]
We present the solution procedures for the case of continuous capillary pressure fields. That case  can be realized under 
 the compatibility condition that the capillary pressure curves have the same endpoints as  in Figure~\ref{cap_rock_types} (right).  The approaches can be applied to discontinuous capillary pressure  by extending the curves of these functions into  monotone graphs(cf.~\cite{MR3208750,MR2559741}). 
\end{rem}
\section{Non-overlapping DD methods with  explicit coupling}\label{DD_FVMF}
We present in this section two numerical approaches for the two-phase flow problem~\eqref{mathematical_problem_DD}--\eqref{Chap1_Matching_pressure} based on domain decomposition methods and operator splitting techniques. 
\subsection{A brief description  of  operator-splitting methods}
 The natural  IMPES-based approach for solving~\eqref{mathematical_problem_DD}--\eqref{Chap1_Matching_pressure} is to use a sequential splitting in which one fixes the saturation in the pressure equations~\eqref{Chap2_System_Conserv_tota}-\eqref{Chap2_System_Darcy_total} with the  corresponding coupling conditions~\eqref{Chap1_Matching_pressure} and solves for the unknown pressure and velocity, which are subsequently held fixed when evolving the saturation  a time step  according to~(\ref{Chap2_System_Conserv_wett})-(\ref{Chap2_System_Darcy_wett}) with the 
coupling conditions~\eqref{Chap2_Matching_saturation}. Here, we involve a further operator splitting of this  problem and \textit{split} the advection-diffusion operator~(see Remark~\ref{rem:advection_dominated}), and apply a multirate time stepping strategy   exploiting the different time scales of the problem components~(details on the flowchart of this method and its modules are given in~Scheme~\ref{flowchart}). In \textit{grosso modo},   once the pressure system is solved through the module \texttt{PRESSURE}\_\texttt{SOLVER}   at the \textit{coarse time step}, say  $t^{n}$, with
time-lagging the saturation,  we advance this later through the module~\texttt{ADVECTION}\_\texttt{SOLVER}  for multiple finer time steps, say  $t^{n,l}$, $1\leq l \leq \mathcal{N}^{\textn{a}}$ (with $t^{n,\mathcal{N}^{\textn{a}}}=t^{n+1}$), and finally advance saturation through the module~\texttt{DIFFUSION}\_\texttt{SOLVER}   at the \textit{coarse time step} $t^{n+1}$. Such  approach    allows for the use of various schemes and DD techniques taking the advantage of the specific properties of each problem,   permits the use of different time steps for advection effects to adapt the different time scales in the different rocks,  and generally reduce memory requirements.
\begin{figure}[!ht]
\renewcommand{\figurename}{Scheme}
\setcounter{figure}{0}    
   \includegraphics[scale=0.69]{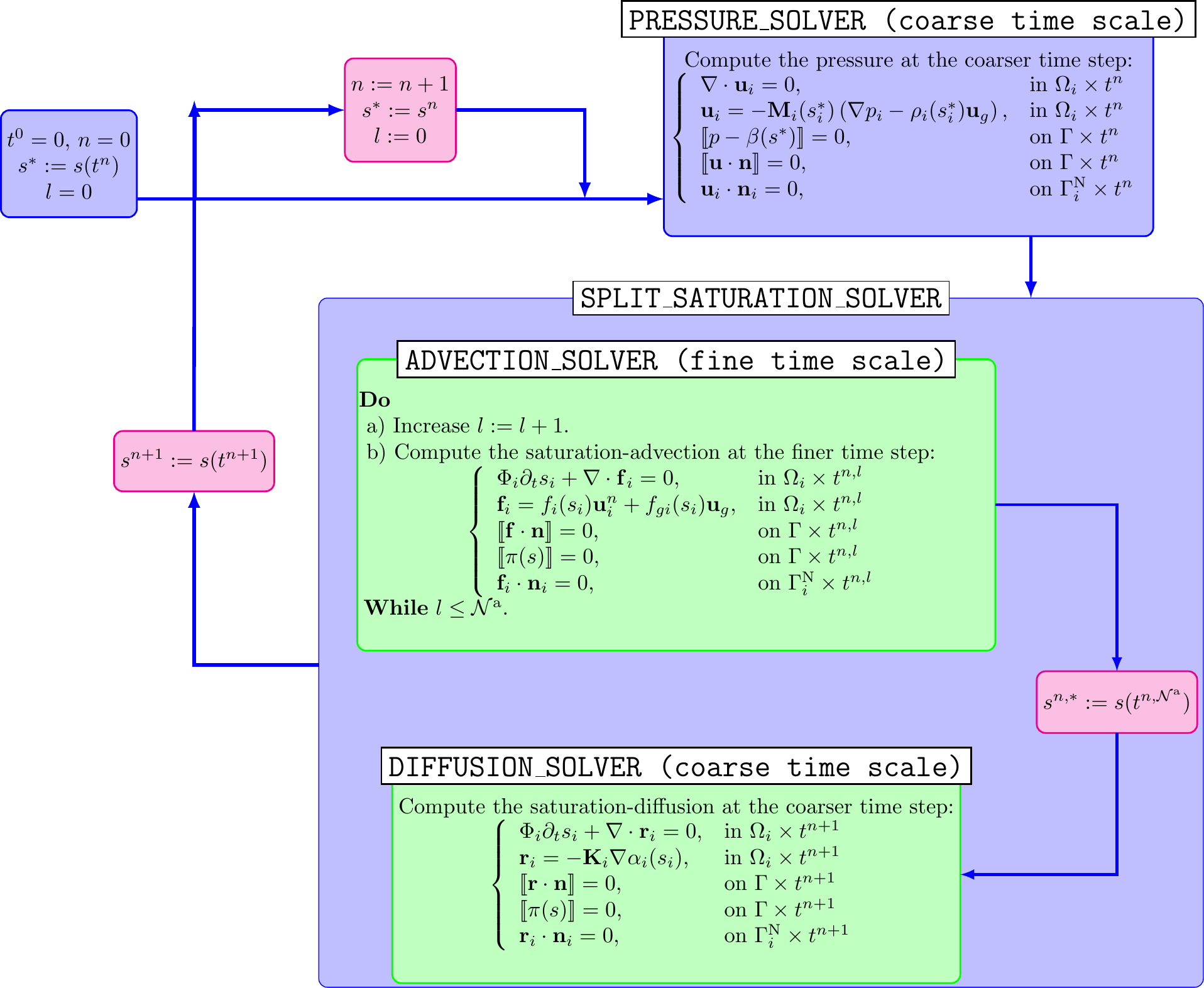}
  \caption{Flowchart of the two-stage splitting scheme with multirate time stepping for two-phase flow problem. The  module~\texttt{ADVECTION}\_\texttt{SOLVER} allows  nonconforming time steps in the different rocks.}\label{flowchart}
\end{figure}
\begin{rem}[Advection-dominated effects]\label{rem:advection_dominated}
Recalling that plugging~\eqref{Chap2_System_Darcy_wett} into~\eqref{Chap2_System_Conserv_wett} results on  a nonlinear  advection-diffusion problem for the saturation. For  many  applications,    advection effects are much greater than  capillary-diffusion effects  through  porous formations. Therefore, it is efficient to compute  fast (advection-effects) solutions on a fine mesh in time and  slow (total velocity + capillarity-effects) solutions on a coarse mesh in time. 
\end{rem}
\subsection{Meshes and function spaces}
We introduce here the partitions of $\mathrm{\Omega}$, time discretization, notation, and function spaces; 
see~\cite{Ahmed2016,ahmed:hal-02275690,MR2206441} for the standard part of the notation. 
\subsubsection{Temporal meshes}
\label{subsubsec:TimeDiscretization}
For integer values $N\geq 0$, let $\left(\tau^{n}\right)_{1\leq n\leq N}$ denote a sequence of positive real numbers corresponding to the discrete time steps 
such that $T:=\sum_{n=1}^{N}\tau^{n}$. Let $t^{0}:=0$, and $t^{n}:=\sum_{j=1}^{n}\tau^{j}, \ 1 \leq n\leq N$ be the discrete times. Let
$I^n:=(t^{n-1},t^n], \ 1 \le n \le N$. 
 For every time step $0\leq n \leq N$, we let  $v_{h}^{n} := v_{h}(\cdot,t^{n})$ for any sufficiently smooth function $v_{h}$.
\subsubsection{Spatial meshes}
\begin{figure}[!ht]
  \centering
   \includegraphics[scale=0.25]{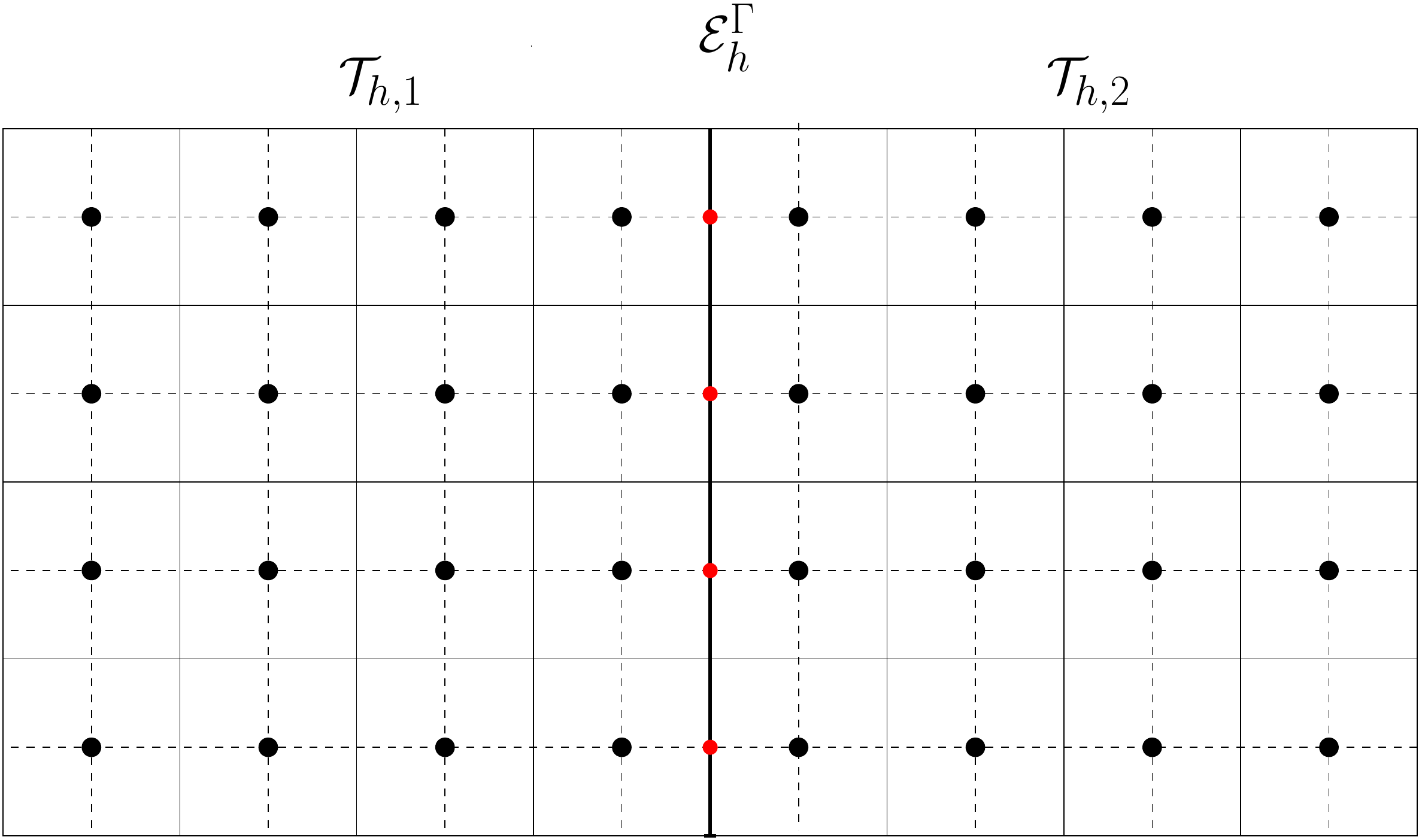}
  \caption{Example of conforming meshes in two subdomains in 2D.}\label{compatible_meshes}
\end{figure}

We consider a spatial discretization $\mathcal{T}_{h,i}$ of the  domain $\mathrm{\Omega}_{i}$ 
consisting  of  either simplicial or rectangular elements $K\in\mathcal{T}_{h,i}$.  We assume that
these meshes are such that $\mathcal{T}_{h}:=\displaystyle{\cup_{i=1}^{2}\mathcal{T}_{h,i}}$ forms a conforming finite element mesh on all of $\Omega$~(see Figure~\ref{compatible_meshes});
we assume that this partition is 
matching in the sense that for two distinct elements of $\mathcal{T}_h$ their
intersection is either an empty set or their common vertex or edge.  We denote by $h$   the maximal
mesh size of $\mathcal{T}_{h,i}$. The interior mesh faces in $\mathcal{T}_{h,i}$ are collected into
the set $\mathcal{E}_{h,i}^{\textnormal{int}}$, and we denote by $\mathcal{E}_{h,i}$  all the faces of
$\mathcal{T}_{h,i}$  and we set $\mathcal{E}_{h}:=\displaystyle\cup_{i=1}^{2}\mathcal{E}_{h,i}$. We denote by $\mathcal{E}_{h,i}^{\textnormal{N}}$ the sides of
$\mathcal{E}_{h}$ on $\mathrm{\Gamma}^{\textnormal{N}}_{i}$.   Finally,   let $\mathcal{E}^{\mathrm{\Gamma}}_{h}$ be a partition of $\mathrm{\Gamma}$ given by the sides of
$\mathcal{T}_{h}$ on $\mathrm{\Gamma}$ and  we denote by $\mathcal{E}_{K}$ the faces of the element $K\in\mathcal{T}_h$. The volume of an element $K$ is
denoted by $|K|$ and that of a face $\sigma$ by $|\sigma|$. Finally, we use the
notation $\vecx_K$ to denote the ``center'' of the cell $K\in\mathcal{T}_{h}$. If
$\sigma=K|L \in \mathcal{E}_{h}$ separates the cells $K$ and $L$, $d_{K,L}$ denotes
the Euclidean distance between $\vecx_K$ and $\vecx_L$, and $d_{K,\sigma}$
for $\sigma \in \mathcal{E}_{K}$ denotes the distance from $\vecx_K$ to $\sigma$. We assume that the composite mesh $\mathcal{T}_{h}$ satisfies the following
orthogonality condition: for a face $\sigma=K|L\in \mathcal{E}_{h}$, the line segment
$\vecx_{K}\vecx_{L}$ is orthogonal to $\sigma$ (see~\cite{reichenberger2006mixed}). 
 \subsubsection{Function spaces}
\begin{figure}[hbtp]
\centering
\includegraphics[scale=0.2]{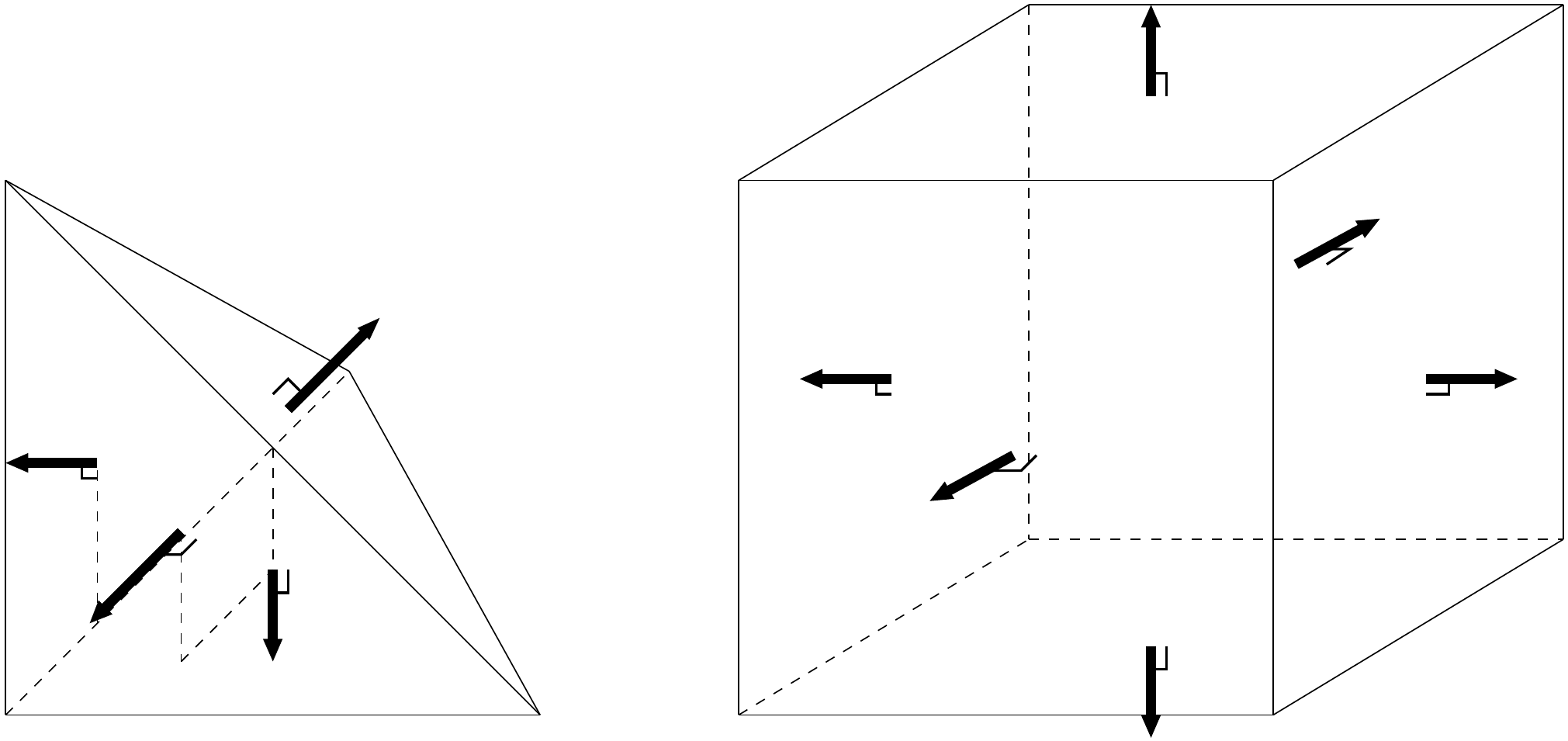}
\caption{Local degrees of freedom for the space $\RTN_0$ on a tetrahedron or a parallelepiped.}
\label{RTN_0}
\end{figure}

We denote by $\mathbb{P}_{l}(S)$ the space of polynomials on a subdomain $S \subset \mathrm{\Omega}$ 
of total degree less than or equal to $l$, and by
$ \mathbb{P}_{l}(\mathcal{T}_{h}) := \{p_{h}\in L^{2}(\Omega);\, p_{h}|_{K}\in \mathbb{P}_{l}(K), \ \forall K \in\mathcal{T}_{h}\}$. Let $\Hdiv$ be the space of
vector-valued functions  from $\left[L^{2}(\mathrm{\Omega})\right]^{d}$ that admit a
weak divergence in $L^{2}(\mathrm{\Omega})$.  We additionally consider the piecewise lowest-order Raviart--Thomas--N\'{e}d\'{e}lec space 
$\mathbf{RTN}_{0}(\mathcal{T}_{h})\subset L^{2}(\Omega)$ defined by $\mathbf{RTN}_{0}(\mathcal{T}_{h}):=\{\vecv_{h}\in L^{2}(\Omega);\, \vecv_{h}|_{K}\in \mathbf{RTN}_{0}(K),\;\forall K\in\mathcal{T}_{h}\}$, where $\mathbf{RTN}_{0}(K):= \left[\mathbb{P}_{0}(\mathcal{T}_{h})\right]^{d}+\mathbb{P}_{0}(\mathcal{T}_{h})\vecx$ (see Figure~\ref{RTN_0}). 
For the discretization of the  subdomain variables, we introduce the  discrete spaces:
\begin{alignat*}{2}
\nonumber& \mathbf{W}_{h,i}:=\dsp\{ \vv_{h}\in \mathbf{RTN}_{0}(\mathcal{T}_{h,i})\cap\mathbf{H}(\divvv,\mathrm{\Omega}_{i});\, \vv_{h}\cdot\vn=0\,\textnormal{ on }\mathrm{\Gamma}_{i}^{\textnormal{N}}\},\\
\label{Chap1_discre_space_scalar}&M_{h,i}:=\mathbb{P}_{0}(\mathcal{T}_{h}).
\end{alignat*}
 The DD  algorithms utilize  \textit{Lagrange multipliers} on the interfaces between the rocks 
 to impose weakly interface conditions. Thus,  to impose the  matching conditions~\eqref{Chap2_Matching_saturation}-\eqref{Chap1_Matching_pressure}, 
 we introduce the  mortar finite element space:
 \begin{equation*}
\label{Chap1_discre_space_Lagrange}\mathrm{\Lambda}_{h}:=\{\mu_{h}\in L^{2}(\mathrm{\Gamma});\,\mu_{h}|_{\sigma}\in\mathbb{P}_{0}(\sigma),\;\forall \sigma \in \mathcal{E}_{h}^{\mathrm{\Gamma}} \}. 
 \end{equation*}
 \subsubsection{Multirate and nonconforming time grids}\label{seb:projection}
\begin{figure}[h]
\begin{center}
{\includegraphics[scale=0.405]{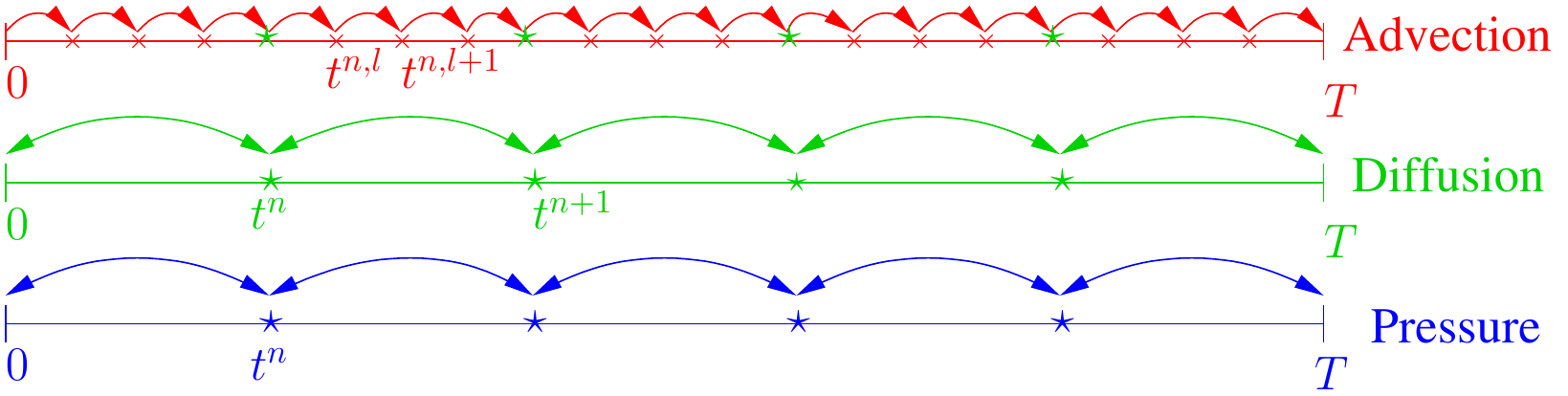}}
{\includegraphics[scale=0.28]{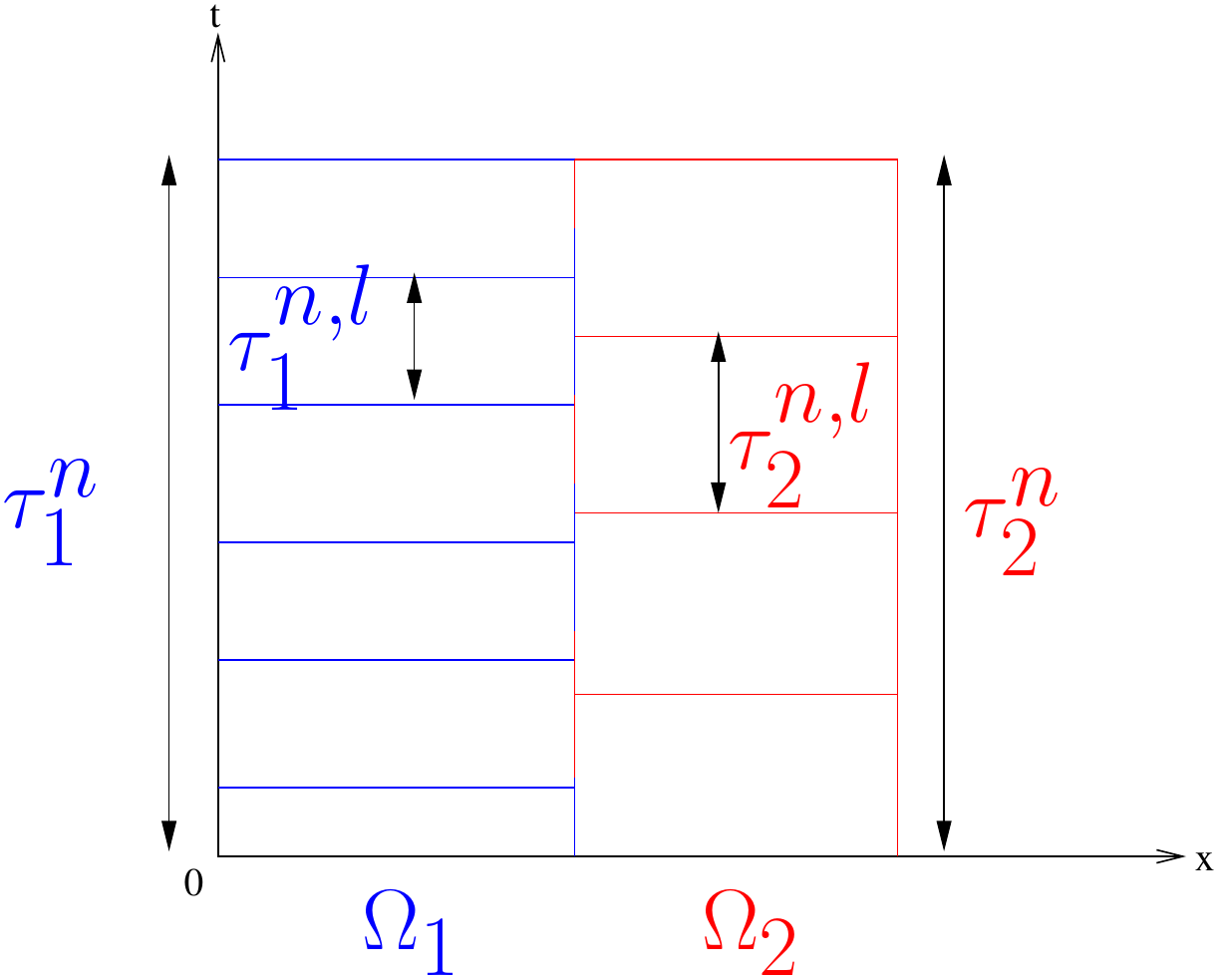}}
\caption{ Time partition with multirate time steps (left) and  nonconforming time grids for advection in the subdomains (right). }
\label{Time_splitting_figure_2}
\end{center}
\end{figure}

For   $1\leq n\leq N$,   let $\mathcal{D}_{1}^{n}$
and $\mathcal{D}_{2}^{n}$ be two possibly different partitions of the time interval $I^{n}$. The finer time step for the advection problem in each interval  $ I^{n}$ is
defined by a sequence  $\left(\tau^{n,l}_{i}\right)_{1\leq l\leq \mathcal{N}^{\textn{a}}_{i}}$, 
where $\mathcal{N}^{\textn{a}}_{i}$ is a positive integer value  such that     $\tau^{n}:=\sum_{i=1}^{\mathcal{N}^{\textn{a}}_{i}} \tau^{n,l}_{i}$. 
We let $t^{n,0}:=t^{n}$, $t^{n,\mathcal{N}^{\textn{a}}_{i}}_{i}:=t^{n+1}$, and   $t^{n,l}_{i}:=\sum_{j=1}^{l}\tau^{n,l}_{j}$, $1 \leq l\leq \mathcal{N}^{\textn{a}}_{i}$,  be the discrete times for the advection in  
 $ I^{n}$. We let
$I^{n,l}_{i}=(t_{i}^{n,l-1},t^{n,l}_{i}], \, 1 \leq l\leq\mathcal{N}^{\textn{a}}_{i}$ (see Figure~\ref{Time_splitting_figure_2}). 

 We denote by $\mathbb{P}_{0}(\mathcal{D}_{j}^{n},\mathrm{\Lambda}_{h})$ the space of piecewise constant functions in time on grid $\mathcal{D}_{j}^{n}$ with values in $\mathrm{\Lambda}_{h}$. To perform the exchange of the data  between the fine time grids
 $\mathcal{D}_{1}^{n}$
and $\mathcal{D}_{2}^{n}$, $1\leq n\leq N$,~(see Figure~\ref{Time_splitting_figure_2} (right)), we
define the  $L^{2}$-projection $\mathcal{P}_{ij}$ from $\mathbb{P}_{0}(\mathcal{D}_{j}^{n},\mathrm{\Lambda}_{h})$ to
$\mathbb{P}_{0}(\mathcal{D}_{i}^{n},\mathrm{\Lambda}_{h})$, as follows~\cite{MR3144798,ahmed2018global}: for   all $0\leq n\leq N$, for $\phi\in \mathbb{P}_{0}(\mathcal{D}_{j}^{n},\mathrm{\Lambda}_{h})$, 
\begin{equation}\label{projection_oper}
\left( \mathcal{P}_{ij}\mathrm{\phi} \right)_{| I^{n,l}_{i}}:=\dfrac{1}{ \tau^{n,l}_{i}}\sum_{m =1}^{\mathcal{N}^{\textn{a}}_{j}} \int_{I^{n,l}_{i}\cap I^{n,m}_{j}}\mathrm{\phi} \, \textnormal{d}t,\quad \forall  l \in \{1, \dotsc, \mathcal{N}^{\textn{a}}_{i} \}.
\end{equation}
We refer the reader to~\cite{MR2235750} for details on 
the implementation of this projection  using an optimal  algorithm  with no additional grid.
\subsection{The module \texttt{PRESSURE}\_\texttt{SOLVER}}\label{sec:interface_pressue}
Let $0\leq n\leq N$ be fixed. The module~\texttt{PRESSURE}\_\texttt{SOLVER} takes as input the saturation  at the present time to linearize the nonlinear coefficients.  So, for  known  
 saturations $s^{n}_{h,i}$, $i\in\{1,2\}$, we let   $\mathcal{A}_{i}^{n}=\mathbf{M}_{i}(s^{n}_{h,i})^{-1}$, and     introduce  a new unknown $\lambda:=(p^{n}_{h,i}-\beta_i(s^{n}_{h,i}))|_{\mathrm{\Gamma}}$,  such that   
 $(\vu^{n}_{h,i},p^{n}_{h,i},\lambda)\in \mathbf{W}_{h,i}\times M_{h,i}\times\mathrm{\Lambda}_{h}$ solves
\bse\label{Chap2_Varia_GP}
\begin{alignat}{3}
\nonumber\dsp- ( p^{n}_{h,i},\div \vv)_{\mathrm{\Omega}_i}+ (\mathcal{A}_{i}^{n}\vu^{n}_{h,i},\vv)_{\mathrm{\Omega}_i}
&+\int_{\Gamma}(\lambda+\beta_i(s^{n}_{h,i})) \vv\cdot \vn_{i}\\
&=\dsp(\rho_{i}(s^{n}_{h,i})\vu_{g}, \vv)_{\mathrm{\Omega}_i},\quad &&\forall \vv\in \mathbf{W}_{h,i}, \\
\dsp( \div \vu^{n}_{h,i}, q)_{\mathrm{\Omega}_i}&=0,\quad &&\forall q\in M_{h,i}.
\end{alignat}\ese 
This system is closed by the weak interface condition
\begin{eqnarray}\label{Chap2_Varia_GP_interface}
\begin{array}{lll}
\dsp\int_{\Gamma}  \jump{\vu^{n}_{h}\cdot \vn}\mu=0,&\quad \forall\mu\in \mathrm{\Lambda}_{h}.
\end{array}
\end{eqnarray}
To efficiently solve  the above multi-domain pressure problem, this system will  
 be reduced   to an interface problem posed   only on $\Gamma$~\cite{ahmed2018multiscale}.  The reduced problem is then solved by an iterative procedure, 
which requires solving subdomain pressure problems at each iteration. For $i\in\{1,2\}$, we let 
\begin{alignat}{2}\label{varia_decom}
 p^{n}_{h,i}:=\overline{p}^{n}_{h,i}(\lambda)+\mathring{p}^{n}_{h,i}\quad\textnormal{and}\quad \vu^{n}_{h,i}:=\overline{\vu}^{n}_{h,i}(\lambda)+\mathring{\vu}^{n}_{h,i},
\end{alignat}
where for $\lambda\in \mathrm{\Lambda}_{h}$,  $(\overline{\vu}^{n}_{h,i},\overline{p}^{n}_{h,i})\in \mathbf{W}_{h,i}\times M_{h,i}$ solves
 \bse\label{Chap2_Varia_GP_0}
 \begin{alignat}{3}
\dsp -(\overline{p}^{n}_{h,i}, \div \vv)_{\mathrm{\Omega}_i}+(\mathcal{A}_{i}^{n} \overline{\vu}^{n}_{h,i},\vv)_{\mathrm{\Omega}_i}&=-\int_{\Gamma}\lambda \vv\cdot \vn_{i},\quad&&\forall \vv\in \mathbf{W}_{h,i},\\ 
( \div \overline{\vu}^{n}_{h,i}, q)_{\mathrm{\Omega}_i} &=0,\quad&&\forall q\in M_{h,i},
\end{alignat}\ese
and where  $(\mathring{\vu}^{n}_{h,i},\mathring{p}^{n}_{h,i})\in \mathbf{W}_{h,i}\times M_{h,i}$  solves
\bse\label{Chap2_Varia_GP_*}
 \begin{alignat}{3}
\nonumber\dsp -\left(\mathring{p}^{n}_{h,i}, \div \vv\right)_{\mathrm{\Omega}_i}+\left(\mathcal{A}_{i}^{n}\mathring{\vu}^{n}_{h,i}, \vv\right)_{\mathrm{\Omega}_i}&+\int_{\Gamma}\beta_i(s_{h,i}^{n}) \vv\cdot \vn_{i}&&\\
&=\left( \rho_{i}(s^{n}_{h,i})\vu_{g}, \vv\right)_{\mathrm{\Omega}_i},\quad&&\forall \vv\in \mathbf{W}_{h,i},\\ 
\dsp \left(\div \mathring{\vu}^{n}_{h,i}, q\right)_{\mathrm{\Omega}_i} &=0,\quad&&\forall q\in M_{h,i}.
\end{alignat}\ese
Define the forms
$s^{n}_{\Gamma,i}:\mathrm{\Lambda}_{h}\times \mathrm{\Lambda}_{h}\rightarrow\RR$, $i\in\{1,2\}$,
$s^{n}_{\Gamma}:\mathrm{\Lambda}_{h}\times \mathrm{\Lambda}_{h}\rightarrow\RR$, and
$g^{n}_{\Gamma}: \mathrm{\Lambda}_{h}\rightarrow\RR$,
\bse\label{bilinear_forms_discrete}\begin{align}
s^{n}_{\Gamma,i}(\lambda,\mu) & :=
 \int_{\Gamma}\mathcal{S}^{\textn{DtN}}_{\Gamma,i,n}(\lambda)\mu:= -\int_{\Gamma}\overline{\vu}^{n}_{h,i}(\lambda)\cdot\vecn_{i}\mu,\\
s^{n}_{\Gamma}(\lambda,\mu) & :=\int_{\Gamma}\mathcal{S}^{\textn{DtN}}_{\Gamma,n}(\lambda)\mu:=\sum_{i=1}^{2}s^{n}_{\Gamma,i}(\lambda,\mu),\\
g^{n}_{\Gamma}(\mu) & :=\int_{\Gamma}g_{\Gamma}^{n}\mu:= \sum_{i=1}^{2}\int_{\Gamma}\mathring{\vu}^{n}_{h,i}\cdot\vecn_{i}\mu,
\end{align}\ese
where $\mathcal{S}^{\textn{DtN}}_{\Gamma,i,n}:\mathrm{\Lambda}_{h}
\rightarrow \mathrm{\Lambda}_{h}$, $1\leq i\leq 2$, and
$\mathcal{S}^{\textn{DtN}}_{\Gamma,n}:=\sum_{i=1}^{2}\mathcal{S}^{\textn{DtN}}_{\Gamma,i,n}$
are   Dirichlet-to-Neumann or Steklov-Poincar\'{e} (SP) type operators. Obviously, the  operator
$\mathcal{S}^{\textn{DtN}}_{\Gamma,i,n}$
 is linear. 
It is now easy  to verify that the pressure problem~\eqref{Chap2_Varia_GP}--\eqref{Chap2_Varia_GP_interface} is reduced to an interface problem:
 \begin{df}[Steklov-Poincar\'{e}  problem]\label{def:reduced scheme}
Find $\lambda\in  \mathrm{\Lambda}_{h}$
such that
\begin{alignat}{4}
 \label{Chap1_interface_pressure}&s^{n}_{\Gamma}(\lambda,\mu)=g^{n}_{\Gamma}(\mu),\quad\forall \mu\in \mathrm{\Lambda}_{h}.
 \end{alignat}
\end{df}
One can show that the Dirichlet-to Neumann operator $\mathcal{S}^{\textn{DtN}}_{\Gamma,n}$ associated to~\eqref{Chap1_interface_pressure}  is a symmetric positive definite one so we  can use   a 
 conjugate gradient method (CG) to calculate $\lambda$. The subdomain solutions are then retrieved via~\eqref{varia_decom}.
 \begin{rem}[Equivalence with the multidomain pressure problem]
In  Definition~\ref{def:reduced scheme}, the 
continuity condition $\jump{p^{n}_{h}-\beta(s^{n}_{h})}=0$   on $\Gamma$ is imposed in space by fixing $\lambda:=(p^{n}_{h,i}-\beta_i(s^{n}_{h,i}))|_{\mathrm{\Gamma}}$, $\forall i\in\{1,2\}$, while the continuity of the total flux $\jump{\vu^{n}_{h} \cdot \vn}=0$ is imposed weakly via~\eqref{Chap1_interface_pressure}.
\end{rem}
 In order to speed up the convergence of~the CG iterations, we may use a 
 \textit{Neumann--Neumann preconditioner} with weighted matrices (cf.~\cite{MR3144798}):
 for  the parameters $a_{i}\in\R$  such that $a_{1}+a_{2}=1$, we introduce the Neumann-to-Dirichlet operators $\mathcal{S}^{\textn{NtD}}_{\Gamma,i,n}$, $i\in\{1,2\}$, given by
  \begin{alignat*}{2}
   && \mathcal{S}^{\textn{NtD}}_{\Gamma,i,n}\,\colon\, \mathrm{\Lambda}_{h}\rightarrow \mathrm{\Lambda}_{h},\,\mathrm{\varphi}\longmapsto p^{n}_{h,i}(\mathrm{\varphi}),\quad i\in\{1,2\},
\end{alignat*} 
 where $(\vu^{n}_{h,i}(\mathrm{\varphi}),p^{n}_{h,i}(\mathrm{\varphi}))$ is such that $\vu^{n}_{h,i}\cdot\vecn_{i}=\mathrm{\varphi}$ on $\mathrm{\Gamma}$ 
and  solves
\bse\label{precond_subprobl}
 \begin{alignat}{3}
\label{precond_subprobl1}\dsp -\left(p^{n}_{h,i}, \div \vv\right)_{\mathrm{\Omega}_i}+\left(\mathcal{A}_{i}^{n} \vu^{n}_{h,i},\vv\right)_{\mathrm{\Omega}_i}&=0,\quad&&\forall \vv\in \mathbf{W}_{h,i},\\
\label{precond_subprobl2}\left( \div \vu^{n}_{h,i}, q\right)_{\mathrm{\Omega}_i} &=0,\quad&&\forall q\in M_{h,i}.
\end{alignat}\ese
The  preconditioned  version (in strong form) of~\eqref{Chap1_interface_pressure} is given by 
\begin{df}[Preconditioned   problem]\label{def:prec_reduced scheme}
   Find $\lambda\in  \mathrm{\Lambda}_{h}$ 
   \begin{equation}
\label{Chap1_interface_pressure_prec} \left(a_{1}\mathcal{S}^{\textn{NtD}}_{\Gamma,1,n}+ a_{2}\mathcal{S}^{\textn{NtD}}_{\Gamma,2,n}\right)\mathcal{S}^{\textn{DtN}}_{\Gamma,n} (\lambda)= g^{n}_{h}.
\end{equation}
  \end{df}
Following~\cite{MR3144798}, an alternative to the Steklov-Poincar\'{e} formulation of Definition~\ref{def:reduced scheme}, can be obtained through an equivalent formulation of the pressure model problem \eqref{Chap2_System_Conserv_tota}-\eqref{Chap2_System_Darcy_total} with the   interface conditions~\eqref{Chap1_Matching_pressure}. Precisely, one can solve equations~\eqref{Chap2_System_Conserv_tota}-\eqref{Chap2_System_Darcy_total}, for  $i\in\{1,2\}$, together with
 \textit{Robin transmission conditions}
  \bse\label{Robin_conditions_press}\begin{alignat}{3} 
&-\vu^{n}_{h,1} \cdot \vn_{1}+\gamma_{1}(p^{n}_{h,1}-\beta(s^{n}_{h,1}))=\vu^{n}_{h,2} \cdot \vn_{2}+\gamma_{1}(p^{n}_{h,2}-\beta(s^{n}_{h,2})),\quad&&\textnormal{on }\Gamma,\\
&-\vu^{n}_{h,2} \cdot \vn_{2}+\gamma_{2}(p^{n}_{h,2}-\beta(s^{n}_{h,2}))=\vu^{n}_{h,1} \cdot \vn_{1}+\gamma_{2}(p^{n}_{h,1}-\beta(s^{n}_{h,1})),\quad&&\textnormal{on }\Gamma,
\end{alignat}\ese
where  $\mathrm{\gamma}_{1}$ and  $\mathrm{\gamma}_{2}$ are strictly positive constants.  Note that~\eqref{Robin_conditions_press} are 
proper linear combinations of the original conditions within \texttt{PRESSURE}\_\texttt{SOLVER} module given by~$\jump{\vu^{n}_{h} \cdot \vn}=0$ and $\jump{p^{n}_{h}-\beta(s^{n}_{h})}=0$ on $\Gamma$.    Therefore, adopting ideas from~\cite{ahmed:hal-01540956,MR3144798,HOANG2017366,seus2017lineartwophase}, one can use 
an \textit{Optimized Schwarz Waveform Relaxation }(OSWR) iterative method  for solving this  problem. This method can be written as follows: 
\begin{algo}[The linear OSWR]~\label{linear_oswr}
{
\setlist[enumerate]{topsep=0pt,itemsep=-1ex,partopsep=1ex,parsep=1ex,leftmargin=1.5\parindent,font=\upshape}
\begin{enumerate}
\item Choose an initial approximation $\mathrm{\eta}^{-1}_{i}\in \mathrm{\Lambda}_{h}$ of $\mathrm{\eta}_{i}:=\vu^{n}_{h,j} \cdot \vn_{j}+\gamma_{i}(p^{n}_{h,j}-\beta(s^{n}_{h,j}))$, for $j=(3-i)$, and a tolerance $\epsilon>0$. Set $k:=-1$.
\item \textn{\textbf{Do}}
\begin{enumerate}
 \item Increase $k:= k+1$.
\item  Compute  $(\vu^{n}_{h,i},p^{n}_{h,i})\in \mathbf{W}_{h,i}\times M_{h,i}$, $i\in\{1,2\}$, such that 
\bse\label{mathematical_problem_RRiter}\begin{alignat}{3} 
\nonumber\dsp- ( p^{k,n}_{h,i},\div \vv)_{\mathrm{\Omega}_i}+ (\mathcal{A}_{i}^{n}\vu^{k,n}_{h,i},\vv)_{\mathrm{\Omega}_i} &&&\\
\nonumber\qquad\qquad\qquad +\dfrac{1}{\gamma_{i}}\dsp\int_{\Gamma}(\vu^{k,n}_{h,i}\cdot\vecn_{i}) (\vv\cdot\vecn_{i})
&=-\dfrac{1}{\gamma_{i}}\int_{\Gamma}(\mathrm{\eta}^{k-1}_{i}+\beta_i(s^{n}_{h,i})) \vv\cdot \vn_{i}&&\\
&\quad+\dsp(\rho_{i}(s^{n}_{h,i})\vu_{g}, \vv)_{\mathrm{\Omega}_i},\quad &&\forall \vv\in \mathbf{W}_{h,i}, \\
\dsp( \div \vu^{k,n}_{h,i}, q)_{\mathrm{\Omega}_i}&=0,\quad &&\forall q\in M_{h,i}.
\end{alignat}\ese
\item Set $\mathrm{\eta}^{k}_{i}:= \vu^{k,n}_{h,j} \cdot \vn_{j}+\gamma_{i}(p^{k,n}_{h,j}-\beta(s^{n}_{h,j}))$, for $j=(3-i)$.
\end{enumerate}
\textn{\textbf{While}} $\dfrac{\|(\mathrm{\eta}^{k}_{1},\mathrm{\eta}^{k}_{2}) - (\mathrm{\eta}^{k-1}_{1},\mathrm{\eta}^{k-1}_{2})\|_{L^{2}(\Gamma)}}{\| (\mathrm{\eta}^{k-1}_{1},\mathrm{\eta}^{k-1}_{2})\|_{L^{2}(\Gamma)}}\geq \epsilon$.
\end{enumerate}}
\end{algo}
The  OSWR method in the context of MFE methods has been  studied in~\cite{HOANG2017366}. Note that the   parameters  $(\mathrm{\gamma}_{1},\mathrm{\gamma}_{2})$ can be optimized to give an improved
convergence rate of this algorithm.  This can be done   by numerically minimizing the convergence factor
corresponding to  the Darcy problem~(cf.~\cite{ahmed:hal-01540956,HOANG2017366} for more
details).
\subsection{The module~\texttt{SPLIT}\_\texttt{SATURATION}\_\texttt{SOLVER}}
Here, we detail the ingredients of the module \texttt{SPLIT}\_\texttt{SATURATION}\_\texttt{SOLVER}. As pointed previously,    the saturation-advection within \texttt{ADVECTION}\_\texttt{SOLVER} module takes multiple finer time steps against one  saturation-diffusion time step  within \texttt{DIFFUSION}\_\texttt{SOLVER}.  The module \texttt{ADVECTION}\_\texttt{SOLVER} employs different time grids between subdomains to adapt the different time scales in the different 
rocks~\cite{HOANG2017366}.

\subsubsection{STEP~1:~\texttt{ADVECTION}\_\texttt{SOLVER}}
\begin{figure}\centering
{\includegraphics[scale=0.3]{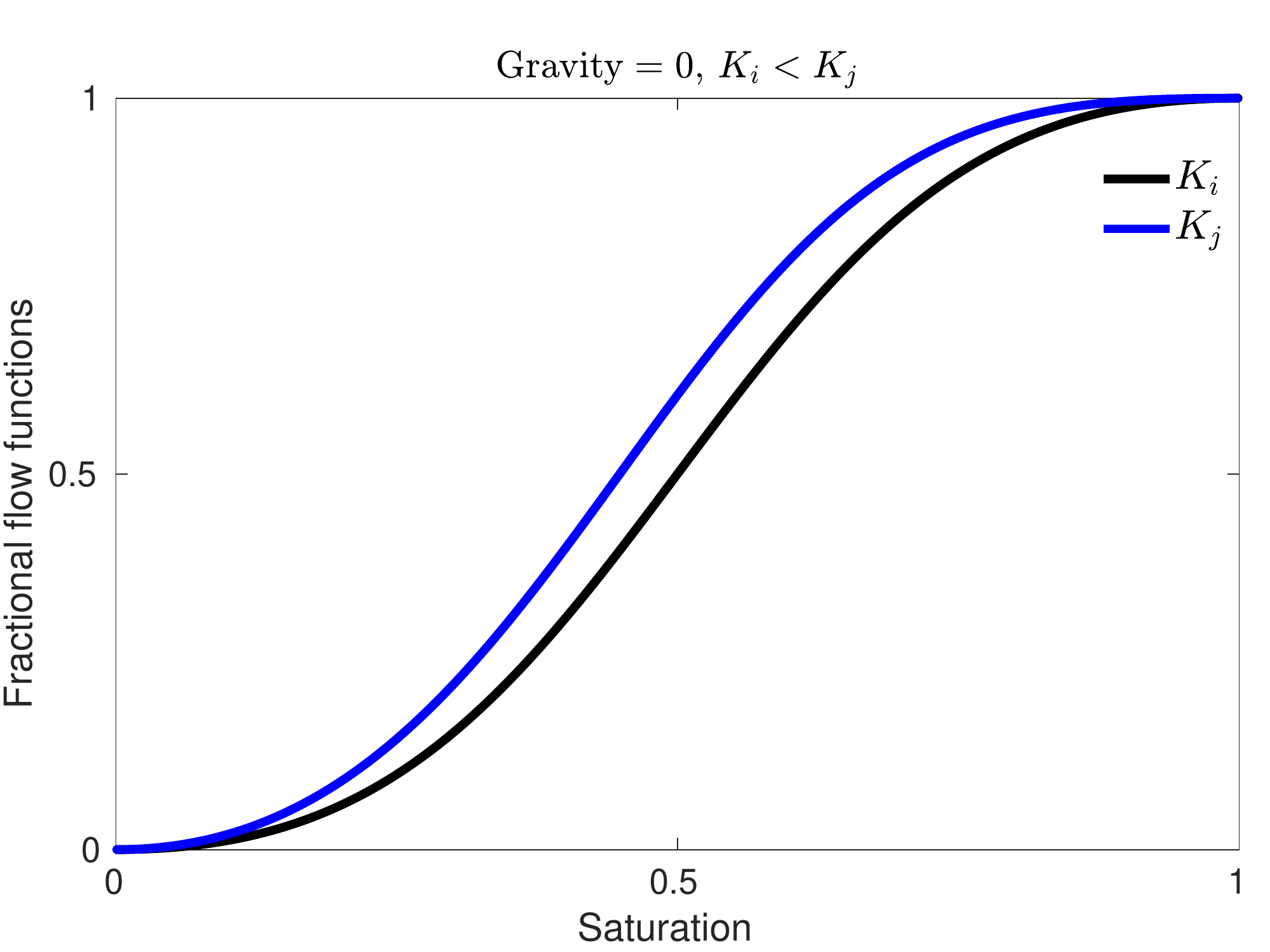}}
{\includegraphics[scale=0.3]{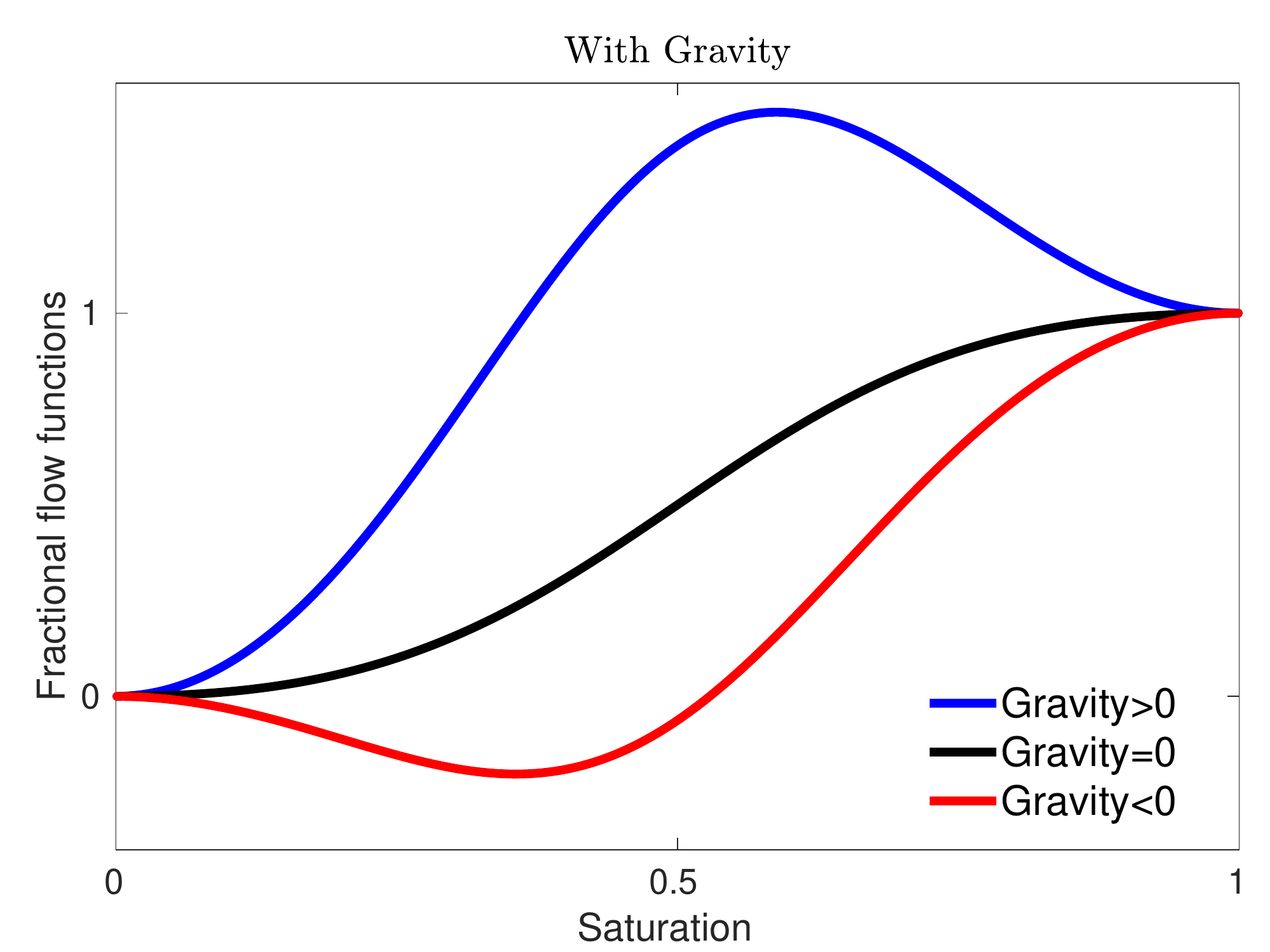}}
\caption{Typical behavior of the fractional flow function for water phase saturation with free-gravity (left),
with gravity effects and dominated advection  (right).}
\label{Advection_Gravity_depends1}
\end{figure}

\begin{figure}[hbtp]
\centering
\includegraphics[scale=0.5]{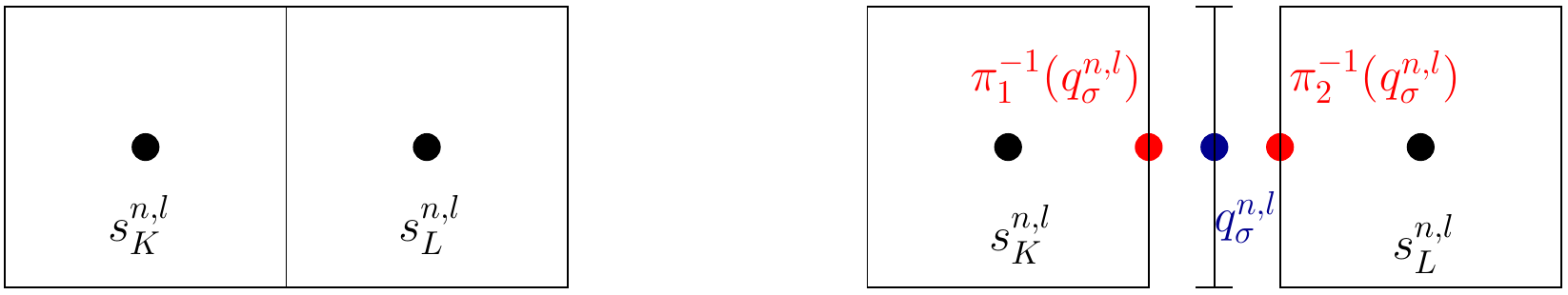}
\caption{Degrees of freedom within one rock type (left) and between two rock types with conforming-in-time grids (right).}
\label{varia_advec}
\end{figure}

The  scheme in~\texttt{ADVECTION}\_\texttt{SOLVER}, is based on conservation of mass element-by-element;
 the subdomain saturation $s^{n,\textnormal{a}}_{h,i}$  is piecewise constant and calculated using first order
cell-centered finite volume method. 
 The initial data $s_{i}^{0}$, $i\in\{1,2\}$, are discretized
as follows:
\begin{equation}\label{init_adv}
 s^{0}_{K,i}=\dfrac{1}{|K|}\int_{K}s_{i}^{0} \text{d}\vecx,\quad\forall K\in \gammaTih.
\end{equation}
For the following time steps, we compute  an intermediate saturation $s^{n,l}_{h,i}$ for all $l\in\{0,\dotsc,\mathcal{N}^{\textn{a}}_{i}-1\}$, by
 \begin{alignat}{2} \label{Chap2_advection_discrete} 
&\displaystyle\sum_{K\in \gammaTih}\int_{K} \mathrm{\Phi}_{K} \dfrac{(s^{n,l+1}_{h,i}-s^{n,l}_{h,i})}{\tau^{n,l+1}_{i}} \,\text{d}\vecx+ \displaystyle\sum_{\sigma\in \mathcal{E}_{K}} |\sigma|\varphi_{K,\sigma}^{n,l}=0,\quad\forall K\in \gammaTih,
\end{alignat}
where $|\sigma|\varphi_{K,\sigma}^{n,l}$ is an approximation of the advection flux through the face $\sigma\in \mathcal{E}_{K}$, $\int_{\mathrm{\sigma}}f_{a,i}^{n}$, with
$f_{a,i}^{n}=\vf_{i}(s)\cdot \vn_{K,\sigma}=(f_{i}(s)\vu_{i} +  f_{gi}(s)\vu_{g})\cdot \vn_{K,\sigma}$ where $\vn_{K,\sigma}$ 
denotes the outward normal to $\sigma$ with respect to $K$. The function $\varphi_{K,\sigma}^{n,l}$ will be defined as a function of the two values of the saturation on the two sides of $\sigma$.
 Let us first suppose the case of conforming time grids between the two rocks, i.e.,  $\mathcal{D}_{1}^{n}=\mathcal{D}_{2}^{n}$.    The numerical flux $ \varphi^{n,l}_{K,\sigma}$ in that case is calculated by 
  \begin{eqnarray}\label{Chap2_advection_flux_definition}
\varphi_{K,\sigma}^{n,l}=\left\{ 
\begin{array}{lll}
 \varphi_{i}^{n}(s_{K}^{n,l},s_{L}^{n,l}),&\quad\textnormal{if } \sigma=K|L \in\mathcal{E}_{h,i}^{\textnormal{int}},\\
 \varphi_{i}^{n}(s_{K}^{n,l},\theta^{n,l}_{K,\sigma}),&\quad\textnormal{if }  \sigma\in\mathcal{E}_{K}\cap\mathcal{E}_{h}^{\mathrm{\Gamma}},\\
 0,&\quad\textnormal{if } \sigma \in\mathcal{E}_{K}\cap\mathcal{E}_{h}^{\textnormal{N}},
\end{array}\right.
\end{eqnarray}
 where   $\theta^{n,l}_{K,\sigma}$ and $\theta^{n,l}_{L,\sigma}$  are two additional unknowns on $\mathrm{\sigma}=K|L\in\mathcal{E}_{h}^{\mathrm{\Gamma}}$, such that $K\in \mathcal{T}_{h,1}$ and $L\in \mathcal{T}_{h,2}$, 
 chosen  to satisfy   the following conditions 
\bse\label{Chap2_advection_flux_definit}
\begin{alignat}{1}\label{Chap2_advection_flux_definit_1}
&\varphi_{1}^{n}(s_{K}^{n,l},\theta^{n,l}_{K,\sigma})+\varphi_{2}^{n}(s_{L}^{n,l},\theta^{n,l}_{L,\sigma})=0,\\
 \label{Chap2_advection_flux_definit_2}&\pi_1(\theta^{n,l}_{K,\sigma})-\pi_2(\theta^{n,l}_{L,\sigma})=0.
\end{alignat}\ese
 We give now the definition of the flux function $\varphi_{i}$. Following \cite{MR2051062,Mishra2017479}, and    taking  advantage that in our case $f_{ai}$ has a particular  bell-shape, 
 (see Figure~\ref{Advection_Gravity_depends1}),
which has either one global maximum and no other local maxima or one global minimum and no other local minima, the function 
$\varphi_{i}(a,b)$, for $(a,b)\in[0,1]^2$ and  $i\in\{1,2\}$, is given by: \\
 When $f_{ai}$ has one maximum:
  \begin{equation*}\begin{array}{lll}
  \varphi_{i}(a,b)= \max\{f_{ai}(\max\{a, \xi_{fi} \}), f_{ai}(\min\{\xi_{fi}, b\})\},\\
  \xi_{fi}=\argmax{f_{ai}}.
 \end{array}\end{equation*}
 When $f_{ai}$ has one minimum:
   \begin{equation*}\begin{array}{lll}
  \varphi_{i}(a,b)= \min\{f_{ai}(\min\{a, \xi_{fi} \}), f_{ai}(\max\{\xi_{fi}, b\})\},\\
  \xi_{fi}=\argmin{f_{ai}}.
 \end{array}\end{equation*}
 In view of the above definition, the flux function  $\varphi_{i}$    is exactly  the Godunov numerical
flux with respect to $f_{ai}$~\cite{MR2051062,MR3208750,Mishra2017479}.
   \begin{rem}[Phase-by-phase function]
  One can also use   the phase-by-phase  upstream flux function instead of the Godunov function, which is simpler for implementation,
 and   whose expression relies on a particular structure  of the flux function $\vf_{i}$ (cf.~\cite{MR3208750,aziz1979petroleum}).
  \end{rem}
To solve numerically the  advection system~\eqref{Chap2_advection_flux_definit} it is desirable  to reduce the number of unknowns by eliminating the  equation~\eqref{Chap2_advection_flux_definit_2}. 
Following~\cite{MR3208750}, finding $(\theta^{n,l}_{K,\sigma},\theta^{n,l}_{L,\sigma})$ solution of~\eqref{Chap2_advection_flux_definit} can be reduced  to the
following problem.
\begin{df}[Conforming-in-time advection problem]
At each time step $t^{n,l}$, find $q_{\sigma}^{n,l}\in\R$ 
such that
\begin{align}\label{Chap2_advection_flux_definit2}
\mathrm{\Psi}^{l}_{\textnormal{a}}(q^{n,l}_{\sigma}):=\varphi_{1}^{n}(s_{K}^{n,l},\pi_1^{-1}(q_{\sigma}^{n,l}))
+\varphi_{2}^{n}(s_{L}^{n,l},\pi_2^{-1}(q_{\sigma}^{n,l}))=0.
\end{align}
\end{df}
This interface problem is with one implicit unknown per interface face (see Figure~\ref{varia_advec}) and can be solved using Newton's method or simply 
using a scalar root finder, e.g.  Regula Falsi method (cf.~\cite{MR3208750}).

 Now, it remains to extend  the above  setting to the nonconforming time grids. To  this aim,   
 we introduce the space-time variables  $\theta^{n,a}_{K,\sigma}:=(\theta^{n,l}_{K,\sigma})_{0\leq l\leq \mathcal{N}^{\textn{a}}_{1}}$ and $\theta^{n,a}_{L,\sigma}:=(\theta^{n,l}_{L,\sigma})_{0\leq l\leq \mathcal{N}^{\textn{a}}_{2}}$, and similarly for $s^{n,a}_{K}$ and $s^{n,a}_{L}$. Therefore,  based on the adoption of Robin interface conditions, i.e., proper linear combinations of the coupling conditions~\eqref{Chap2_advection_flux_definit}, we define 
the interface operator \begin{align}
\label{robin_to_robin_advec}&\mathrm{\Psi}^{l,m}_{\textnormal{a},2\textnormal{L}}(\theta^{n,a}_{K,\sigma},\theta^{n,a}_{L,\sigma}):=
\left\{ 
\begin{array}{lll}\displaystyle\int_{I^{n,l}_{1}}\left(-\varphi_{1}^{n}(s_{K}^{n,a},\theta_{K,\sigma}^{n,a})+\kappa_{1}\pi_{1}(\theta_{K,\sigma}^{n,a})-\mathcal{P}_{12}(\varphi_{2}^{n}(s_{L}^{n,a},\theta_{L,\sigma}^{n,a})+\kappa_{1}\pi_{2}(\theta_{L,\sigma}^{n,a}))\right)\textnormal{d}t,\\
\displaystyle\int_{I^{n,m}_{2}}\left(-\varphi_{2}^{n}(s_{L}^{n,a},\theta_{L,\sigma}^{n,a})+\kappa_{2}\pi_{2}(\theta_{L,\sigma}^{n,a})-\mathcal{P}_{21}(\varphi_{1}^{n}(s_{K}^{n,a},\theta_{K,\sigma}^{n,a})+\kappa_{2}\pi_{1}(\theta_{K,\sigma}^{n,a}))\right)\textnormal{d}t,
\end{array}\right.
\end{align}
where  $\mathrm{\kappa}_{1}$ and  $\mathrm{\kappa}_{2}$ are strictly positive constants, and  $\mathcal{P}_{12}$ and $\mathcal{P}_{21}$,   are  projection-in-time operators defined by~\eqref{projection_oper} (see Subsection~\ref{seb:projection} for more details). The nonconforming-in-time counterpart of the flux continuity condition~\eqref{Chap2_advection_flux_definit2}  is then 
replaced  with a space-time formulation. 
\begin{df}[Nonconforming-in-time advection problem]
 Find $(\theta^{n,a}_{K,\sigma},\theta^{n,a}_{L,\sigma})\in \R^{\mathcal{N}^{\textn{a}}_{1}}\times \R^{\mathcal{N}^{\textn{a}}_{2}}$ such  that
\begin{align}\label{Chap2_advection_flux_definit2_Robin}
&\mathrm{\Psi}^{l,m}_{\textnormal{a},2\textnormal{L}}(\theta^{n,a}_{K,\sigma},\theta^{n,a}_{L,\sigma}):=0,\quad\forall\sigma=K|L\in\mathcal{E}_{h}^{\mathrm{\Gamma}},\quad\forall 0\leq l \leq \mathcal{N}^{\textn{a}}_{1},\,\forall0 \leq m \leq \mathcal{N}^{\textn{a}}_{2}.
\end{align}
\end{df}
This  interface system is solved iteratively with Newton's method~(cf.~\cite{ahmed:hal-01540956,HOANG2017366} for a similar framework), where  the interface unknowns $\theta^{n,a}_{K,\sigma}$ and $\theta^{n,a}_{L,\sigma}$   cooperate to solve the interface  problem~\eqref{Chap2_advection_flux_definit2_Robin}. Note that the operator within~\eqref{Chap2_advection_flux_definit2} is continuous and strictly monotone,
 therefore, this scheme  admits a unique solution (cf.~\cite{MR3208750}). The space-time formulation~\eqref{Chap2_advection_flux_definit2_Robin} is not studied, but we  can at least claim that the operator~\eqref{robin_to_robin_advec} inherits the  monotonicity  in function of each of the interface variables,  as both the projection $\mathcal{P}_{ij}$ and intergration preserve that property. The stability of~\eqref{Chap2_advection_flux_definit2_Robin}  is only verified numerically~(see Section~\ref{section:numerics}). 
 \begin{rem}[Interface data]
 The couple $(\theta^{n,a}_{K,\sigma},\theta^{n,a}_{L,\sigma})$ on  $\sigma=K|L\in\mathcal{E}_{h}^{\mathrm{\Gamma}}$ is  the  common  saturation values, expressing the continuity of the capillary pressure across $\sigma$, while ensuring  the flux transmission. 
\end{rem}
\subsubsection{STEP~2:~\texttt{DIFFUSION}\_\texttt{SOLVER}}\label{diffusion_interface}
We inject the intermediate saturations   $s^{n,*}_{h,i}:=s^{n,\mathcal{N}^{\textn{a}}_{i}}_{h,i}$ calculated through~\texttt{ADVECTION}\_\texttt{SOLVER}  in the module~\texttt{DIFFUSION}\_\texttt{SOLVER}, to calculate the saturations at the coarse time step $n+1$. The scheme of this module is based on MFE scheme for the subdomain diffusion problems and the nonlinear Steklov-Poincar\'{e} and  OSWR methods for the DD approaches, as used in the~\texttt{PRESSURE}\_\texttt{SOLVER} module. The reader is then invited to compare the two  modules  as the ideas and notations are analogous. First,  we define for  a specified capillary pressure $\lambda_{d}:=\pi_{i}(s^{n+1}_{h,i})|_{\Gamma}\in \mathrm{\Lambda}_{h}$, $i\in\{1,2\}$,   the \textit{non-linear  Dirichlet-to-Neumann} operator: 
 \begin{alignat}{2}\label{op:dir_to_neuma}
   && \mathcal{Z}^{\textn{DtN}}_{\Gamma,i,n+1}\,:\, \mathrm{\Lambda}_{h}\times M_{h,i}\, \longrightarrow\,  \mathrm{\Lambda}_{h}, \, (\lambda_{d},s^{n,*}_{h,i})&\longmapsto\vr^{n+1}_{h,i} \cdot \vn_{i},
\end{alignat} 
where  $(\vr^{n+1}_{h,i},s^{n+1}_{h,i})\in \mathbf{W}_{h,i}\times M_{h,i}$, $i\in\{1,2\}$,  is calculated by solving 
the diffusion problem  inside each subdomain together with $\pi_{i}(s^{n+1}_{h,i})=\lambda_{d}$ on $\Gamma$, using MFE method (cf.~\cite{MR3392446,robertsJean});
 \bse\label{Chap2_Varia_Diffusion_theta}
 \begin{alignat}{3}
-( \alpha_{i}(s^{n+1}_{h,i}),\nabla\cdot\vv)_{\mathrm{\Omega}_{i}}+(\vK_{i}^{-1}\vr^{n+1}_{h,i},\vv)_{\mathrm{\Omega}_{i}}&=-\int_{\Gamma}\alpha_{i}(\pi_{i}^{-1}(\lambda_{d})) \vv\cdot \vn_{i},\quad&&\forall \vv\in \mathbf{W}_{h,i},\\
\dsp (\mathrm{\Phi}_{i} \displaystyle\dfrac{s^{n+1}_{h,i}-s^{n,*}_{h,i}}{\tau^{n+1}_{i}},\mu)_{\mathrm{\Omega}_{i}}+(\nabla\cdot\vr^{n+1}_{h,i},\mu)_{\mathrm{\Omega}_{i}}&=0,\quad&&\forall \mu\in M_{h,i}. 
\end{alignat}\ese
 We then define the forms
$d^{n+1}_{\Gamma,i}:\mathrm{\Lambda}_{h}\times \mathrm{\Lambda}_{h}\rightarrow\RR$, $i\in\{1,2\}$, and
$d^{n+1}_{\Gamma}:\mathrm{\Lambda}_{h}\times \mathrm{\Lambda}_{h}\rightarrow\RR$:
\bse\label{nonlinear_bilinear_forms_discrete}\begin{align}
d^{n+1}_{\Gamma,i}(\lambda_{d},\mu) & :=
 \int_{\Gamma}\mathcal{Z}^{\textn{DtN}}_{\Gamma,i,n+1}(\lambda_{d})\mu:= \int_{\Gamma}\vr^{n+1}_{h,i}(\lambda_{d})\cdot\vecn_{i}\mu,\\
d^{n+1}_{\Gamma}(\lambda,\mu) & :=\int_{\Gamma}\mathcal{Z}^{\textn{DtN}}_{\Gamma,n+1}(\lambda)\mu:=\sum_{i=1}^{2}d^{n+1}_{\Gamma,i}(\lambda_{d},\mu),
\end{align}\ese
where $\mathcal{Z}^{\textn{DtN}}_{\Gamma,i,n+1}:\mathrm{\Lambda}_{h}
\rightarrow \mathrm{\Lambda}_{h}$, $1\leq i\leq 2$, and
$\mathcal{Z}^{\textn{DtN}}_{\Gamma,n+1}:=\sum_{i=1}^{2}\mathcal{Z}^{\textn{DtN}}_{\Gamma,i,n+1}$ are non-linear  Dirichlet-to-Neumann operators. 
Now, we are ready to rewrite the equations of this module as  an interface problem posed only on $\Gamma$.
 \begin{df}[Nonlinear Steklov-Poincar\'{e} problem]\label{def:nonlinear_reduced scheme}
Find $\lambda_{d}\in  \mathrm{\Lambda}_{h}$
such that
\begin{alignat}{4}
 \label{Chap1_interface_diffusion}&d^{n+1}_{\Gamma}(\lambda_{d},\mu)=0,\quad\forall\mu\in \mathrm{\Lambda}_{h}.
 \end{alignat}
\end{df}
\begin{rem}[Equivalence with the multidomain diffusion problem]
In  Definition~\ref{def:nonlinear_reduced scheme}, the 
continuity of the capillary pressure across $\Gamma$, i.e., $\jump{\pi(s^{n+1}_{h})}=0$,  is imposed in space by fixing $\lambda_{d}=\pi_{i}(s^{n+1}_{h,i})|_{\Gamma}$, $\forall i\in\{1,2\}$, while the continuity of the diffusive flux $\jump{\vr^{n+1}_{h} \cdot \vn}=0$ is retrieved weakly by solving~\eqref{Chap1_interface_diffusion}.
\end{rem}
\begin{rem}[On Newton-Krylov method for~\eqref{Chap1_interface_diffusion}]\label{rem:Newton_P8GMRES}
If we employ the  Newton-GMRes method (cf.~\cite{kelley1995}) for solving~\eqref{Chap1_interface_diffusion}, the  Newton step is defined  by 
$\lambda_{d}^{k+1}=\lambda_{d}^{k}+\eta^{k}$. Each  step $\eta^{k}$ is computed by a forward difference  GMRes iteration by solving the linear interface-Jacobian problem~$\left(\mathcal{Z}^{\textn{DtN}}_{\Gamma,n+1}\right)^{\prime}(\lambda_{d}^{k})\eta^{k}=-\mathcal{Z}^{\textn{DtN}}_{\Gamma,n+1}(\lambda_{d}^{k})$.  
 A well-known drawback of the GMRes algorithm for solving such interface problem is that the number of iterations depends essentially on  the
number of subdomain solves, therefore depends strongly on the subdomain  discretization. A  preconditioner
is usually needed to reduce the number of iterations to a reasonable level. A left preconditioned
GMRes strategy is based on solving
$P^{-1} \left( \mathcal{Z}^{\textn{DtN}}_{\Gamma,n+1}\right)^{\prime}(\lambda_{d}^{k})\eta^{k}=-P^{-1}\mathcal{Z}^{\textn{DtN}}_{\Gamma,n+1}(\lambda_{d}^{k}),$
where $P$ is an easily invertible approximation to the Jacobian. Physically,  
$P^{-1}$ can be interpreted as a non-linear Neumann-to-Dirichlet operator. Inverting this
non-linear function would lead to a non-linear preconditioner and to resolve that one can construct an approximation of $P$ by solving a linear version of the problem as detailed in \cite{yotov2001interface}. 
\end{rem}

An alternative to the Steklov-Poincar\'{e} formulation~\eqref{Chap1_interface_diffusion}, is to solve the    diffusion problem  in the module~\texttt{DIFFUSION}\_\texttt{SOLVER}  
with the OSWR scheme (see the linear case  in~\texttt{PRESSURE}\_\texttt{SOLVER}) (cf.~\cite{ahmed:hal-01540956,ahmed:hal-02275690}):
\begin{algo}[The nonlinear OSWR]~\label{nonlinear_oswr}
{
\setlist[enumerate]{topsep=0pt,itemsep=-1ex,partopsep=1ex,parsep=1ex,leftmargin=1.5\parindent,font=\upshape}
\begin{enumerate}
\item Choose an initial approximation $\mathrm{\eta}^{-1}_{i}\in \mathrm{\Lambda}_{h}$ of $\mathrm{\eta}_{i}:=\vr^{n+1}_{h,j}\vn_j  +\mathrm{\kappa}_{i}\pi_{j}(s^{n+1}_{h,j})$, with $j=(3-i)$, and a parameter $\mathrm{\kappa}_{i}\in\RR^{+}\setminus\{0\}$. Define a tolerance $\epsilon>0$. Set $k:=-1$.
\item \textn{\textbf{Do}}
\begin{enumerate}
 \item Increase $k:= k+1$.
\item  Compute  $(\vr^{n+1}_{h,i},s^{n+1}_{h,i})\in \mathbf{W}_{h,i}\times M_{h,i}$, $i\in\{1,2\}$, such that 
\bse \label{OSWR_Diffusion_semi_disc}\begin{alignat}{3}
-( \alpha_{i}(s^{k,n+1}_{h,i}),\nabla\cdot\vv)_{\mathrm{\Omega}_{i}}+(\vK_{i}^{-1}\vr^{k,n+1}_{h,i},\vv)_{\mathrm{\Omega}_{i}}&=-\int_{\Gamma}\alpha_{i}(s^{k,n+1}_{h,i}) \vv\cdot \vn_{i},\quad&&\forall \vv\in \mathbf{W}_{h,i},\\
\dsp (\mathrm{\Phi}_{i} \displaystyle\dfrac{s^{k,n+1}_{h,i}-s^{n,*}_{h,i}}{\tau^{n+1}_{i}},q)_{\mathrm{\Omega}_{i}}+(\nabla\cdot\vr^{k,n+1}_{h,i},q)_{\mathrm{\Omega}_{i}}&=0,\quad&&\forall q\in M_{h,i},\\
\int_{\Gamma}\{-\vr^{k,n+1}_{h,i}\cdot\vn_i+\kappa_{i}\pi_{i}(s^{k,n+1}_{h,i})\}\mu&=\int_{\Gamma}\mathrm{\eta}^{k-1}_{i}\mu,\quad&&\forall \mu\in \Lambda_{h}.
\end{alignat}\ese
\item Set $\mathrm{\eta}^{k}_{i}: =\vr^{k,n+1}_{h,j}\cdot\vn_j+\kappa_{i}\pi_{j}(s^{k,n+1}_{h,j})$, for $j=(3-i)$.
\end{enumerate}
\textn{\textbf{While}} $\dfrac{\|(\mathrm{\eta}^{k}_{1},\mathrm{\eta}^{k}_{2}) - (\mathrm{\eta}^{k-1}_{1},\mathrm{\eta}^{k-1}_{2})\|_{L^{2}(\Gamma)}}{\| (\mathrm{\eta}^{k-1}_{1},\mathrm{\eta}^{k-1}_{2})\|_{L^{2}(\Gamma)}}\geq \epsilon$.
\end{enumerate}}
\end{algo}
The  well-posedness of the Robin subdomain problem as well as of the algorithm was shown in~\cite{ahmed:hal-02275690}. Note that  the  proof of convergence can be obtained only for the non-degenerate case, following the techniques of~\cite{MR3144798}. The   free parameters $\mathrm{\kappa}_{i}$, $i\in\{1,2\},$  can be  optimized to ensure  faster  convergence~(cf.~\cite{ahmed2018global,ahmed:hal-01540956,MR3144798}). 
\begin{rem}[At convergence]
 At the convergence of Algorithm~\ref{nonlinear_oswr}, the following
 \textit{Robin transmission conditions } are retrieved
  \bse\label{Robin_conditions_diff}\begin{alignat}{3} 
&-\vr^{n+1}_{h,1} \cdot \vn_{1}+\kappa_{1}\pi_{1}(s^{n+1}_{h,1})=\vr^{n+1}_{h,2} \cdot \vn_{2}+\kappa_{1}\pi_{2}(s^{n+1}_{h,2}),\quad&&\textnormal{on }\Gamma,\\
&-\vr^{n+1}_{h,2} \cdot \vn_{2}+\kappa_{2}\pi_{2}(s^{n+1}_{h,2})=\vr^{n+1}_{h,1} \cdot \vn_{1}+\kappa_{2}\pi_{1}(s^{n+1}_{h,1}),\quad&&\textnormal{on }\Gamma.
\end{alignat}\ese
These conditions are equivalent to the physical ones, i.e., $\jump{\vr^{n+1}_{h} \cdot \vn}=0$ and $\jump{\pi(s^{n+1}_{h})}=0$ on $\Gamma$. 
\end{rem}
 \subsection{The two-stage splitting  algorithm}
 We  provide the   algorithm of Scheme~\ref{flowchart},  by  assembling  
  \texttt{PRESSURE}\_\texttt{SOLVER}  and 
 \texttt{saturation}\_\texttt{SOLVER}, which 
 controls their execution and interaction. 
 \begin{algo}[The overall algorithm of Scheme~\ref{flowchart}]~\label{algo_1}
{
\setlist[enumerate]{topsep=0pt,itemsep=-1ex,partopsep=1ex,parsep=1ex,leftmargin=1.5\parindent,font=\upshape}
\begin{enumerate}
\item Choose  initial saturations $s^{0}_{h,i}\in M_{h,i}$, $i\in\{1,2\}$. Set $n:=-1$.
\item \textn{\textbf{Do}}
\begin{enumerate}
\item  Increase $n:= n+1$. 
\item  Compute $(\vu^{n}_{h,i},p^{n}_{h,i})\in \mathbf{W}_{h,i}\times M_{h,i}$, $i\in\{1,2\}$, at the coarse pressure time step  $t^{n}$, by performing a CG on~\eqref{Chap1_interface_pressure_prec}, or by using Algorithm~\ref{linear_oswr} (OSWR method).
\item Compute  intermediate saturations $s^{n,l}_{h,i}\in  M_{h,i}$, $i\in\{1,2\}$, through advection,  at the finer advection time steps  $t^{n,l}_{i}$, $1\leq l \leq \mathcal{N}^{\textn{a}}_{i}$, by performing Newton's method on~\eqref{Chap2_advection_flux_definit2_Robin}.
\item Set $s^{n,*}_{h,i}:=s^{n,\mathcal{N}^{\textn{a}}_{i}}_{h,i}$, $i\in\{1,2\}$.
 \item Compute the  saturation $s^{n+1}_{h,i}\in  M_{h,i}$, $i\in\{1,2\},$ through diffusion, by performing a Newton-Krylov method on~\eqref{Chap1_interface_diffusion}.
 \end{enumerate}
 \textn{\textbf{While}}  $n\leq N$.
 \end{enumerate}}
\end{algo}
\subsection{An improved IMPES method  with multirate/nonconforming time grids}\label{Scheme2}
In this section, we   review the classical IMPES method, where the  problem is decoupled to pressure and saturation problems and then solved sequentially~\cite{Douglas1983,MR2837398}. 
Thus,  we replace  the module~\texttt{SPLIT}\_\texttt{SATURATION}\_\texttt{SOLVER}  in Scheme~\ref{flowchart} with a specialized solver of the  coupled  advection-diffusion problem (\texttt{ONE}\_\texttt{SATURATION}\_\texttt{SOLVER}). This second approach is called  Scheme~2. The flowchart of Scheme~2 can be drawn with simple modifications on~Scheme~\ref{flowchart}. The drawback of  this   IMPES  based scheme is that  extremely small time steps should be used for overcoming the stability restriction. In order to accelerate the performance of the Scheme~2 and make it capable of solving larger problems  (1) the pressure (fast solution) can cope with a much coarser time step compared to the saturation (slow solution) (2)  the advection-diffusion problem for the saturation  allows   nonconforming  time steps between different rock types.
\subsubsection{The  module \texttt{ONE}\_\texttt{SATURATION}\_\texttt{SOLVER}}
The module~\texttt{ONE}\_\texttt{SATURATION}\_\texttt{SOLVER} takes as input the pressures from the previous coarse time step $n$, and  solves 
the saturation problem (Eqs.~(\ref{Chap2_System_Conserv_wett})--(\ref{Chap2_System_Darcy_wett}) with the IC~\eqref{Chap2_Matching_saturation})
for the finer time steps $t^{n,l}$,  $l \in \{1, \dotsc, \mathcal{N}^{\textn{a}}_{i}\}$, $i\in\{1,2\}$. The scheme uses  an \textit{explicit-in-time finite volume} method for the discretization and the OSWR algorithm  for the domain decomposition. To introduce the FV scheme, we first define  the numerical flux $ F^{n,l}_{K,\sigma}$ over a face $\sigma\in \mathcal{E}_{h}\cap \mathcal{E}_{K}$, $K\in \mathcal{T}_{h,i}$, by 
  \begin{eqnarray}\label{Chap2_advectiondiff_flux_definition}
F_{K,\sigma}^{n,l}=\left\{ 
\begin{array}{lll}
 \varphi_{K,\sigma}^{n,l}+\vK_{i}\dfrac{\alpha_{i}(s_{K}^{n,l})-\alpha_{i}(s_{L}^{n,l})}{d_{K,L}},&\quad\textnormal{if } \sigma=K|L \in\mathcal{E}_{h,i}^{\textnormal{int}},\\
 \varphi_{K,\sigma}^{n,l}+\vK_{i}\dfrac{\alpha_{i}(s_{K}^{n,l})-\alpha_{i}(\theta_{K,\sigma}^{n,l})}{d_{K,\sigma}},&\quad\textnormal{if }  \sigma\in\mathcal{E}_{K}\cap\mathcal{E}_{h}^{\mathrm{\Gamma}},\\
 0,&\quad\textnormal{if } \sigma \in\mathcal{E}_{K}\cap\mathcal{E}_{h}^{\textnormal{N}}.
\end{array}\right.
\end{eqnarray}
 Therein, $\varphi_{K,\sigma}^{n,l}$ is the Godunov flux   function given by~\eqref{Chap2_advection_flux_definition} and where $\theta^{n,l}_{K,\sigma}$ is the unknown face saturation on $\sigma\in\mathcal{E}_{K}\cap\mathcal{E}_{h}^{\mathrm{\Gamma}}$. The numerical flux function $F_{K, \sigma}^{n,l}$ is then  the sum of  a diffusion contribution  due to the capillary effects  and an advection contribution due to the gravity effects and the total flow rate. The FV scheme for the multi-domain saturation is given by
\begin{equation}\label{init_adv_diff}
 s^{0}_{K}=\dfrac{1}{|K|}\int_{K}s^{0} \text{d}\vecx, \quad \forall K\in \mathcal{T}_{h},
\end{equation}
and for  $l\in\{1,\dotsc,\mathcal{N}^{\textn{a}}_{i}\}$, the  discrete saturations  $s^{n,l}_{h,i}$, $i\in\{1,2\}$, is computed by
 \begin{alignat}{2} \label{Chap2_advectiondiff_discrete} 
&\displaystyle\sum_{K\in \gammaTih}\mathrm{\Phi}_{K} \dfrac{s^{n,l+1}_{h,i}-s^{n,l}_{h,i}}{\tau^{n,l+1}_{i}}+ 
\displaystyle\sum_{\sigma\in \mathcal{E}_{K}} |\sigma|F_{K,\sigma}^{n,l}=0,\quad\forall K\in \gammaTih,
\end{alignat}
and such that the saturation traces $\theta^{n,l}_{K,\sigma}$ and $\theta^{n,l}_{L,\sigma}$ on 
both sides of $ \sigma=K|L$,  $K\in \mathcal{T}_{h,1}$, and  $L \in \mathcal{T}_{h,2}$, satisfies the following conditions 
\bse\label{Chap2_advection_flux_definit_NC}
\begin{alignat}{1}\label{Chap2_advection_flux_definit_1NC1}
&-F_{K,\sigma}^{n,l}+\kappa_{1}\pi_1(\theta^{n,l}_{K,\sigma})=\dfrac{1}{\tau^{n,l}_{1}}\int_{I_{1}^{n,l}}\mathcal{P}_{12}(F_{L,\sigma}(t)+\kappa_{1}\pi_2(\theta_{L,\sigma})(t))\textnormal{d}t,\\
\label{Chap2_advection_flux_definit_1NC2}
&-F_{L,\sigma}^{n,l}+\kappa_{2}\pi_2(\theta^{n,l}_{L,\sigma})=\dfrac{1}{\tau^{n,l}_{2}}\int_{I_{2}^{n,l}}\mathcal{P}_{21}(F_{K,\sigma}(t)+\kappa_{2}\pi_1(\theta_{K,\sigma})(t))\textnormal{d}t,
\end{alignat}\ese
where  $\mathrm{\kappa}_{i}$, $i\in\{1,2\},$  are strictly positive constants. The reader is invited to compare the FV scheme~\eqref{init_adv_diff}--\eqref{Chap2_advection_flux_definit_NC}  with that   presented within \texttt{ADVECTION}\_\texttt{SOLVER}, as the ideas are analogous and many of the explanations given there are carried over directly to the present  case (see slo~\cite{ahmed:hal-02275690}). 
\begin{rem}[On~equations~\eqref{Chap2_advection_flux_definit_NC}]\label{rem:nonconform_satur}
 The  equations~\eqref{Chap2_advection_flux_definit_NC} are simply a FV discretization over the time windows $\textnormal{on }\Gamma\times I^{n}$ of the Robin transmission conditions
  \bse\label{Robin_conditions_satur}\begin{alignat}{3} 
&-(\vr_{1}+\vf_{1})\cdot \vn_{1}+\kappa_{1}\pi_{1}(s_{1})=(\vr_{2}+\vf_{2}) \cdot \vn_{2}+\kappa_{1}\pi_{2}(s_{2}),\\
&-(\vr_{2}+\vf_{2}) \cdot \vn_{2}+\kappa_{2}\pi_{2}(s_{2})=(\vr_{1}+\vf_{1}) \cdot \vn_{1}+\kappa_{2}\pi_{1}(s_{1}).
\end{alignat}\ese
\end{rem}
\subsubsection{The OSWR Method}\label{sec:OSWR_satur}
The solution of the FV scheme~\eqref{init_adv_diff}--\eqref{Chap2_advection_flux_definit_NC}  is computed using the OSWR algorithm~\cite{ahmed:hal-01540956,MR3144798}: At the iteration
$k \geq 1$, for $i\in\{1,2\}$,  let $\left(\theta_{h,i}^{0,n,l}\right)_{0\leq l\leq \mathcal{N}^{\textn{a}}_{i}}$  be given,  we then compute the couple 
$(s^{k,n,l}_{h,i},\theta_{h,i}^{k,n,l})$ for all  $ l \in \{1, \dotsc, \mathcal{N}^{\textn{a}}_{i} \}$, by solving the following  problem:  for $i\in\{1,2\}$,  
\begin{equation}
  \label{Interior_discretization_multi-domain_OSWR}
\begin{aligned}
\displaystyle\sum_{K\in \gammaTih}\mathrm{\Phi}_{K} \dfrac{s^{k,n,l+1}_{h,i}-s^{k,n,l}_{h,i}}{\tau^{n,l+1}_{i}}+ \displaystyle\sum_{\sigma\in \mathcal{E}_{K}} |\sigma|F_{K,\sigma}^{k,n,l}=0,\quad\forall K\in \gammaTih,
\end{aligned}
\end{equation}
with the transmission conditions
\begin{alignat}{1}
\label{Exterior_discretization_OSWR_1}
 -F_{K,\sigma}^{k,n,l}+\kappa_{i}\pi_i(\theta^{k,n,l}_{K,\sigma})=\dfrac{1}{\tau^{n,l}_{i}}\int_{I_{i}^{n,l}}\mathcal{P}_{ij}(F_{L,\sigma}^{k-1}(t)+\kappa_{i}\pi_j(\theta^{k-1}_{L,\sigma})(t))\textnormal{d}t,\quad\forall \sigma \in\mathcal{E}_{K}\cap\mathcal{E}_{h}^{\mathrm{\Gamma}},
\end{alignat}
 for $j =(3-i)$,  and  $L \in \mathcal{T}_{h,j}$. For an overview and further details on the OSWR scheme for an explicit-in-time scheme, we refer the reader to~\cite{HOANG2017366}.
 \subsubsection{The modified IMPES}
The  algorithm of Scheme~2  assembles \texttt{PRESSURE}\_\texttt{SOLVER}  and \texttt{ONE}\_\texttt{SATURATION}\_\texttt{SOLVER}   in a coupling algorithm.
 \begin{algo}[The overall algorithm of Scheme~2]~\label{algo_2}
{
\setlist[enumerate]{topsep=0pt,itemsep=-1ex,partopsep=1ex,parsep=1ex,leftmargin=1.5\parindent,font=\upshape}
\begin{enumerate}
\item Choose  initial saturations $s^{0}_{h,i}\in M_{h,i}$, $i\in\{1,2\}$. Set $n:=-1$.
\item \textn{\textbf{Do}}
\begin{enumerate}
\item  Increase $n:= n+1$. 
\item  Compute $(\vu^{n}_{h,i},p^{n}_{h,i})\in \mathbf{W}_{h,i}\times M_{h,i}$, $i\in\{1,2\}$, at the coarse pressure time step  $t^{n}$, by performing a CG method on~\eqref{Chap1_interface_pressure_prec}, or by using Algorithm~\ref{linear_oswr} (OSWR method).
\item Compute the saturation $s^{n,l}_{h,i}\in  M_{h,i}$, $i\in\{1,2\}$, at the finer saturation time steps  $t^{n,l}_{i}$, $l\in\{1,\cdots,\mathcal{N}^{\textn{a}}_{i}\}$, by performing the OSWR algorithm~\eqref{Interior_discretization_multi-domain_OSWR}-\eqref{Exterior_discretization_OSWR_1}.
 \item Set $s^{n+1}_{h,i}:=s^{n,\mathcal{N}^{\textn{a}}_{i}}_{h,i}$. 
 \end{enumerate}
 \textn{\textbf{While}}  $n\leq N$.
 \end{enumerate}}
\end{algo}
 \section{Numerical results}\label{section:numerics}
 In this section, we discuss numerical solutions to the
incompressible two-phase flow problem in three space
dimensions using Algorithms~\ref{algo_1} and~\ref{algo_2}.   In many  examples, we use the 
tolerances listed in Table~\ref{Tolerance_table}. In all the experiments presented here the absolute permeability tensor
$\vK$ is actually a scalar value $K$, and the Mualem--Van Genuchten model is used for the relative permeability and capillary pressure curves~\cite{chavent1986mathematical}:
\begin{alignat}{3}
 \nonumber&\pi(s)=p_{\pi}\big((1-s)^{-1/m}-1\big)^{1/n},\\
 \label{vgcaprel}&k_{w}(s)=\sqrt{s}[1-(1-s^{n})^{m}]^{2},\\
 \nonumber&k_{n}(s)=(1-s)^{2}(1-s^{n})^{2m},
\end{alignat}
with $n=2.8$, $m=1-1/n$ and  $p_{\pi}=\mathrm{\alpha}\sqrt{\mathrm{\Phi}/K}$,  
 where $\alpha$ is a proportionality constant or the Van Genuchten factor. 
Note that a small value of the Van  Genuchten factor $\alpha$ indicates an advection-dominated problem. All of our tests include gravity effects.  It is also worth noting that in all tests,  the advective flux due to   total flow rate dominates the one due to the gravity effects.

  \begin{table}[h]
  \centering
\begin{tabular}{|l || *{3}{>{\centering}p{1cm}|}c|}
\hline
 Interf. solver
    & CG    &   OSWR  &   Newton  &   GMRes                  \\
    \hline
     Tol.
    & 1E-6    &   1E-6  &   1E-4  &   1E-3\\                    
    \hline \hline
    Subd. solver
    & \multicolumn{2}{c|}{Newton} 
        & \multicolumn{2}{c|}{GMRes}\\ 
                                                   
    \hline
    Tol
    & \multicolumn{2}{c|}{1E-4} 
        & \multicolumn{2}{c|}{1E-4} \\                                     
    \hline
\end{tabular}
\caption{Relative tolerances.}
\label{Tolerance_table}
\end{table}
\begin{rem}[Adaptive time stepping]
 In Algorithm~\ref{algo_1} and~\ref{algo_2}, one can use an adaptive  time stepping  strategy  in order  to  systematize  the  determination  of  adequate coarse and fine time steps for each problem. This strategy can be build based on a posteriori error estimates and distinguishing the various error components as used in the multirate scheme in~\cite{MR3983155}.
\end{rem}
\begin{rem}[On the implementation] The two algorithms   are implemented in  the Matlab Reservoir Simulation Toolbox~\cite{MRST-Lie2014,lie_2019}. The meshes are produced using  the  three-dimensional surface meshing software BLSURF interfaced with GHS3D~\cite{FreyGeorge2008} software for the  three-dimensional volumetric meshes.  The implementation  uses the automatic differentiation feature of the toolbox to compute the Jacobian matrices for solving the nonlinear diffusion subdomain problems by the Newton method~(cf.~\cite{lie_2019}).
\end{rem}

\subsection{Test case~1: saturation-diffusion problem between two rock types}
\begin{figure}[h]
\begin{center}
\includegraphics[scale=0.35]{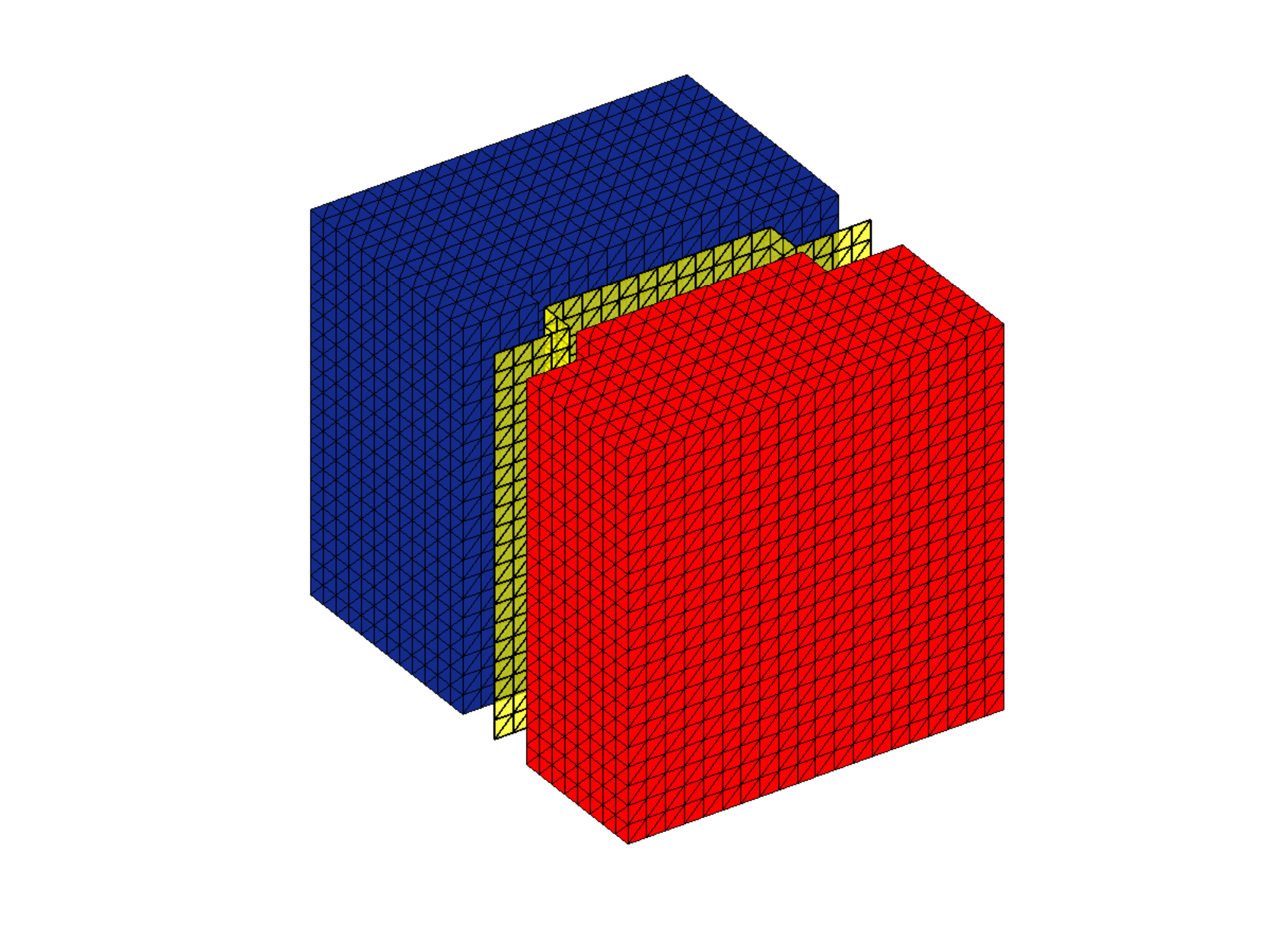}
\includegraphics[scale=0.305]{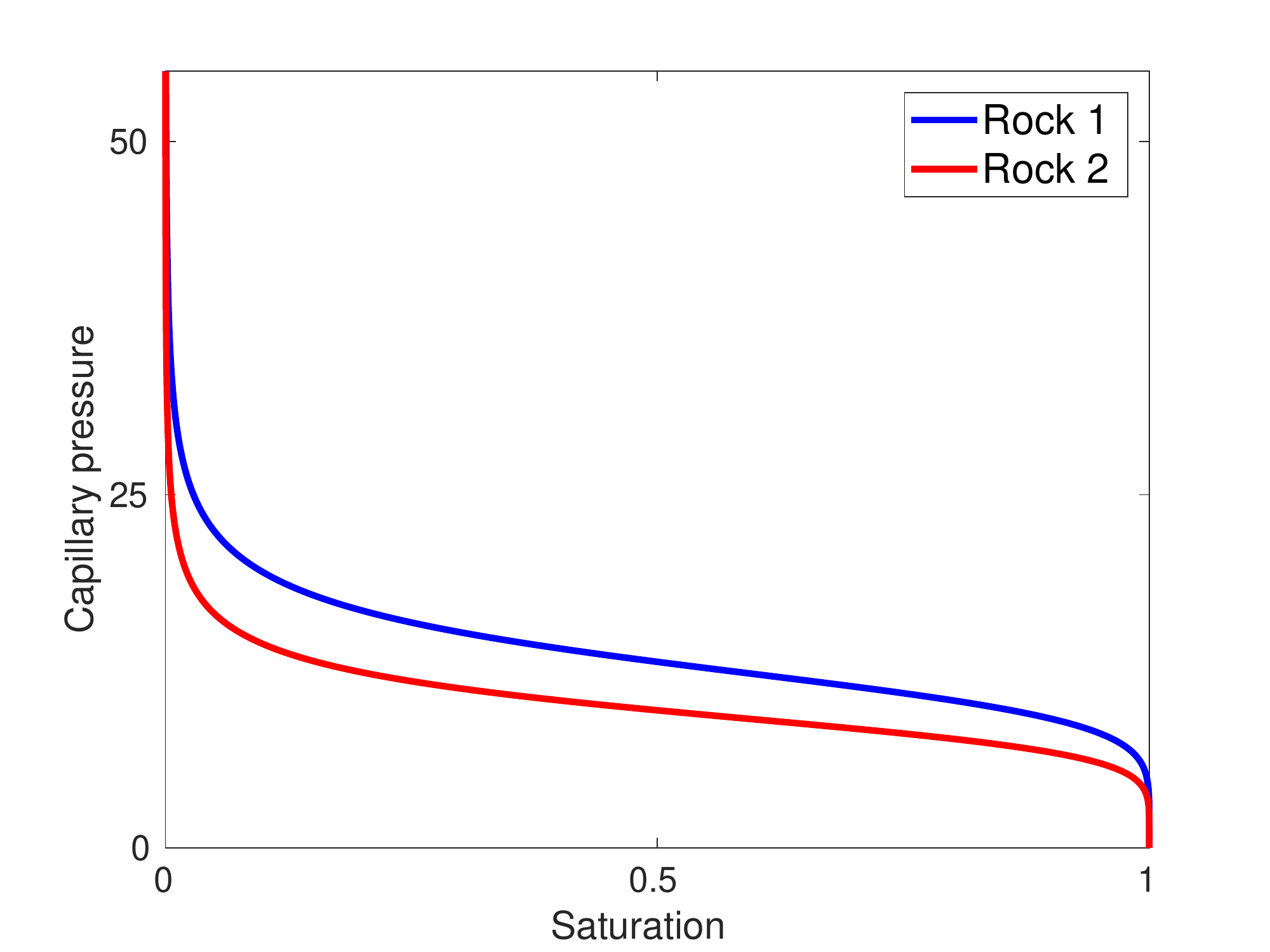}
\caption{Test case~1: an exploded view of the division into two subdomains  (left) and  associated capillary pressure curves (right). }
\label{Two-domains_diffusion}
\end{center}
\end{figure}

In the first case, we  particularize the two-phase flow model  presented above for the sole saturation-diffusion equation (the total velocity being neglected)~\cite{ahmed:hal-02275690}. The goal of this test case is to assess and validate the module~\texttt{DIFFUSION}\_\texttt{SOLVER}. We study the convergence behavior of the Newton-GMRes method applied on~\eqref{Chap1_interface_diffusion} and the OSWR method given by  Algorithm~\ref{nonlinear_oswr}. 
We fix $T=250$ s  and let $\mathrm{\Omega}=[0,1]^{3}$ be decomposed into two subdomains with two rock types
(see Figure~\ref{Two-domains_diffusion}). 

For the spatial discretization, we use a uniform tetrahedral mesh with 48000
elements.  In~\eqref{vgcaprel}, we let $\alpha=100$  for Rock~1 and so is $p_{\pi,1}=4.47$ psi, and  then  $\alpha=140$  so that $p_{\pi,2}=3.52$ psi for Rock~2. For the two rocks,
we set $\mathrm{\Phi}=0.25$, and  $K=1$.  The saturation is set equal to $0.97$ on $\mathrm{\Gamma}^\text{in} :=
\{x=0 \}$. On the outflow boundary $\mathrm{\Gamma}^\text{out} :=
\{x=1 \}$, the saturation at time $t^{n+1}$
is set equal to that inside the closest cell at time $t^n$ (cf.~\cite{ahmed:hal-01540956,Ahmed2016}). We assume homogeneous Neumann boundary conditions on
the remaining part of the boundary. The initial condition is taken to be 0.95 in $\{(x,y,z)\in\mathrm{\Omega}\,\, |\,\,x<0.45 \}$
 and zero elsewhere which fits the continuity of the capillary pressure at the interface between the
rocks.  

\begin{figure}\centering{
      \includegraphics[scale=0.4]{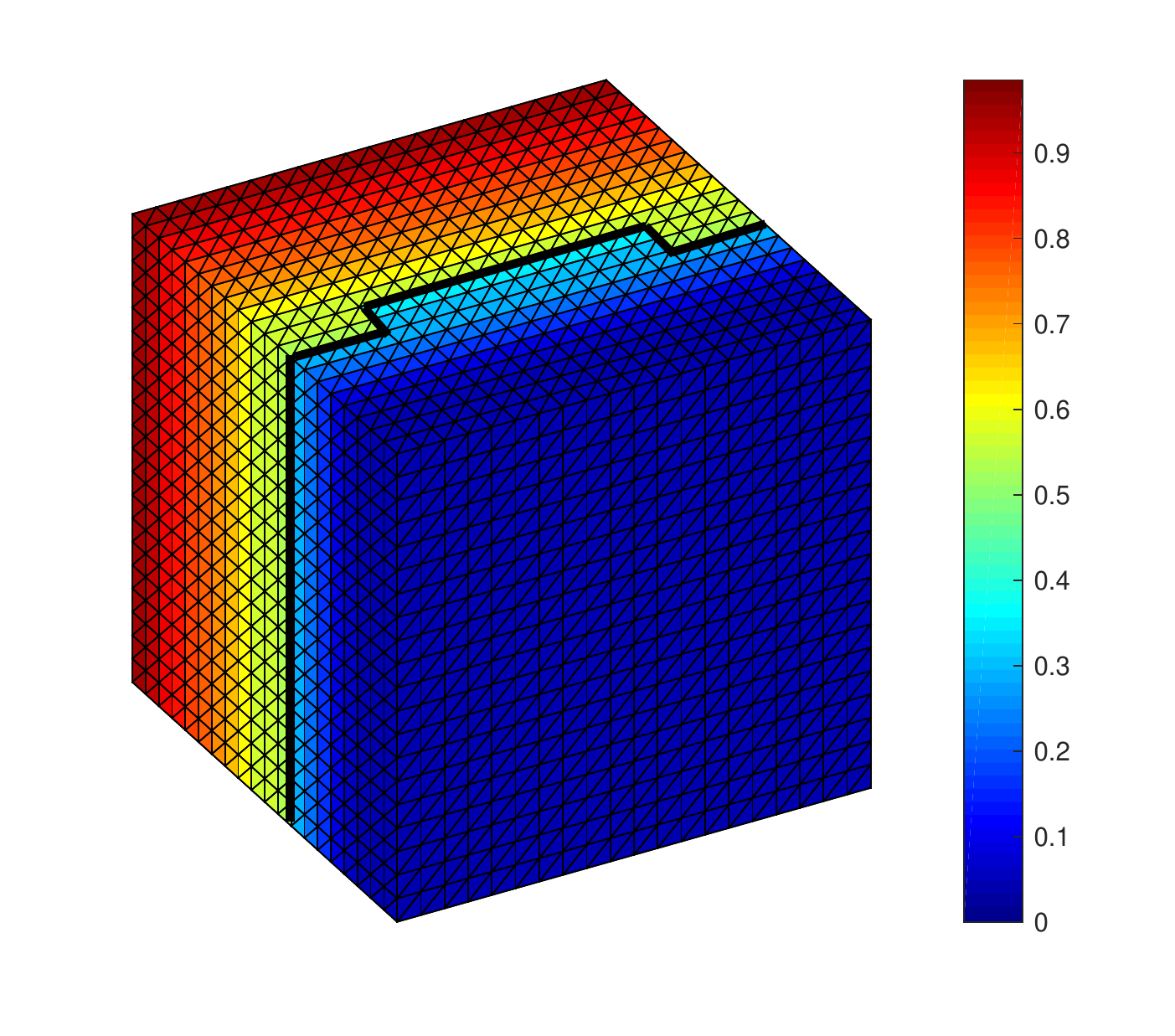}
       \includegraphics[scale=0.4]{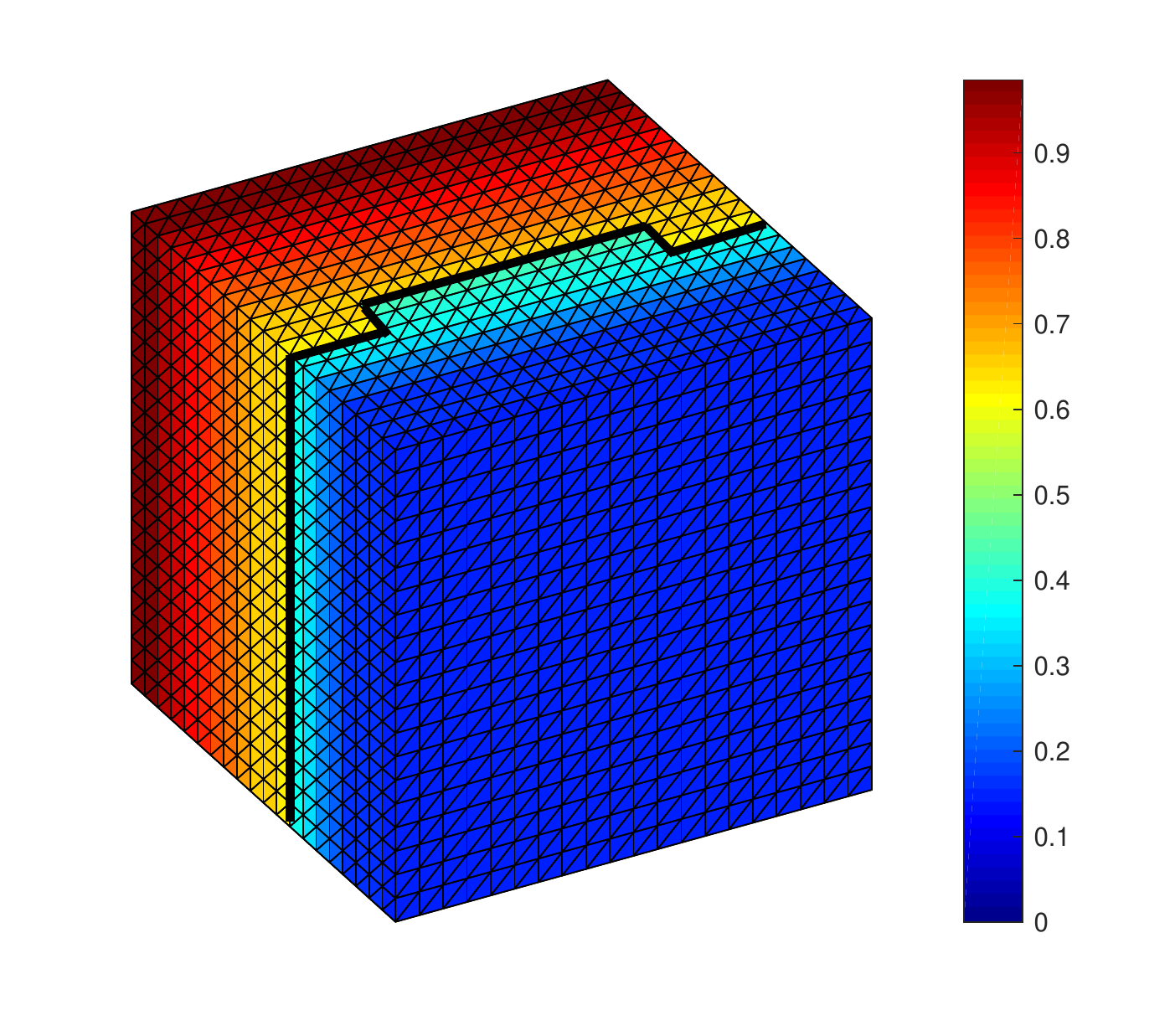}\\
{\includegraphics[scale=0.4]{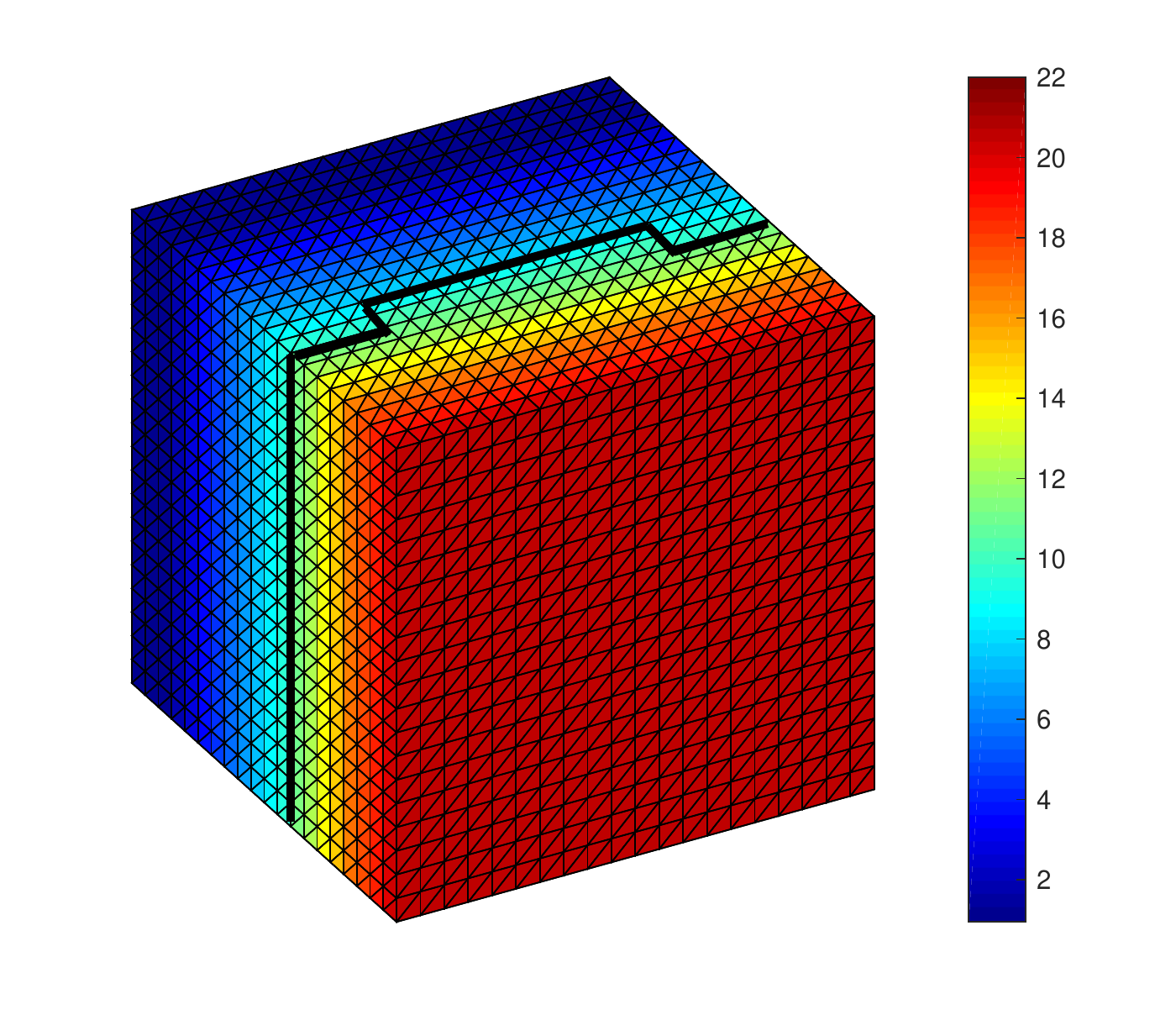}}
{\includegraphics[scale=0.4]{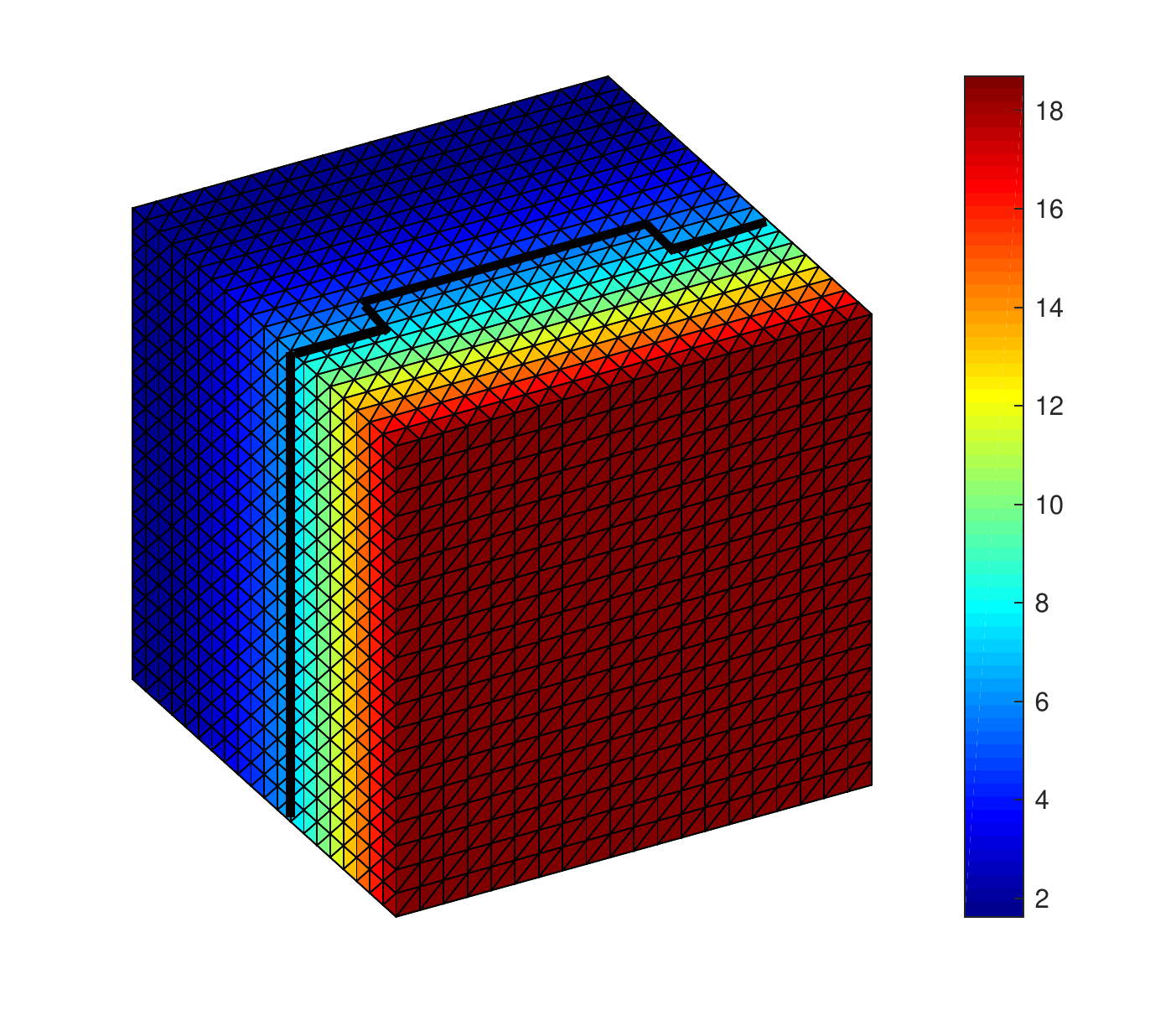}}}
\caption{Test case~1:  saturation $s(t)$ (top) and  capillary pressure $\pi(s(t))$ (bottom) for $t=83$ and $t=113$.}
\label{Two-domains_diffusion_saturation}
\end{figure}
The evolution of the saturation at two time steps is shown in Figure~\ref{Two-domains_diffusion_saturation} (top). We remark 
a very sharp, discontinuous change in saturation at the interface between the two rocks. 
Clearly, the gas cannot penetrate to the domain $\mathrm{\Omega}_{2}$ with the same intensity 
due to the change in the capillary pressure function. Figure~\ref{Two-domains_diffusion_saturation} (bottom) shows the   capillary pressure at two time steps.
We see that, unlike  the saturation, the capillary pressure is continuous  at the
interface between the two rocks, highlighting   numerically that the  
transmission conditions are satisfied.

We come now to the  convergence analysis of the iterative methods involved by \texttt{DIFFUSION}\_\texttt{SOLVER} (Newton-GMRes vs OSWR). First,  we  compare in Table~\ref{Newton_table} the results obtained with the Newton-GMRes  method without and  with preconditioner (see Remark~\ref{rem:Newton_P8GMRES}). We then  present 
in Table~\ref{NewvsOSWR_table} the results obtained using the OSWR method.
\begin{table}[H]
 \centering \begin{tabular}{|l|l|l|}
\hline
\multicolumn{3}{|c|}{Newton-GMRes method}\\
\hline
Interf. Newton & Interf. GMRes& Subd. Newton\\
\hline
\multicolumn{3}{|c|}{without precond.}\\
\hline
Tot.\quad Avg.&Tot.\quad Avg.&Tot.\quad Avg.\\
\hline
318 \quad 6.36&835\quad\;16.72&1660\quad13.23\\
\hline
\multicolumn{3}{|c|}{with precond.}\\
\hline
Tot.\quad Avg.&Tot.\quad Avg.&Tot.\quad Avg.\\
\hline
206 \quad 4.11&705\quad\;11.2&1510\quad9.21\\
\hline
\end{tabular}
\caption{Test case~1: computational iterations of the Newton-GMRes method.}
\label{Newton_table}
\end{table}

 We see   that the Newton-GMRes method with  preconditioner  improves the convergence rate compared to the case with no preconditioner, i.e.,  
 the average number of  Newton iterations  is reduced from 6.36 to 4.11.  That  of the GMRes method,  the average number of  iterations  is reduced from 16.72 to 11.2  iterations. In terms of number of subdomain solves,  
 the  results shows that the   method without a preconditioner needs more subdomain solves. However, we note that each  preconditioned GMRes iteration  costs 
 twice as much as one iteration  of the  method without preconditioning (to construct the preconditioner $P$).  
 \begin{table}[H]
 \centering \begin{tabular}{|l|l|}
 \hline
\multicolumn{2}{|c|}{OSWR method}\\
  \hline
Interf. OSWR& Subd. Newton\\
\hline
Tot.\quad Avg.&Tot.\quad Avg.\\
403\quad \;8.15&1695\quad 13.19\\
\hline
\end{tabular}
\caption{Test case~1: computational iterations of the OSWR method.}
\label{NewvsOSWR_table}
\end{table}
 For the  results obtained with the OSWR method,  the method is  slower than  the Newton-GMRes method and the average number of  OSWR iterations is 8.15 iterations. This can be explained  by the fact the free parameters in~\eqref{OSWR_Diffusion_semi_disc} are not the optimal ones.  In terms of number of 
 subdomain solves, the results are comparable with the method without preconditioning, even though OSWR method needs  different subdomain 
 solves (with nonlinear robin boundary condition).

\subsection{Test case~2: fully two-phase flow  between  two rock types}
We consider a numerical experiment from \cite{alboin2000domain} given by the displacement of a non-wetting fluid by a wetting fluid in a domain $\mathrm{\Omega}=[0,10]^{3}$ 
made up of  two  different rock types and for  $T=5\cdot 10^6$s.
We let Rock~1 have an absolute  permeability equal to $1$  millidarcies and five  times larger than that on Rock~2. The porosity is fixed  to $0.3$ in Rock~1 and $0.7$ in Rock~2. The injection boundary is taken orthogonal  to the interface $\mathrm{\Gamma}$ 
between the two rock types. The initial saturation is set equal to $0.05$ in Rock~1 while in Rock~2 it is set to satisfy  
the equality of the capillary pressure on the interface. 
\subsubsection{Computational performance of Algorithm~\ref{algo_1}}
 We evaluate here the computational performance of Algorithm~\ref{algo_1}. The coarse time steps within the modules \texttt{DIFFUSION}\_\texttt{SOLVER} and \texttt{PRESSURE}\_\texttt{SOLVER} are  fixed in size, i.e., $\tau^{n}=1\cdot 10^{4}$s.  The  finer time steps for the advection in the subdomains are also fixed in size and we choose $\tau^{n,l}_{1}=5\cdot10^2$s in Rock~1 and    $\tau^{n,l}_{2}=1\cdot10^3$s in Rock~2.
We choose the tolerances for the various algorithms involved in Algorithm~\ref{algo_1} from~Table~\ref{Tolerance_table}. 

\begin{figure}[H]\centering
{\includegraphics[scale=0.45]{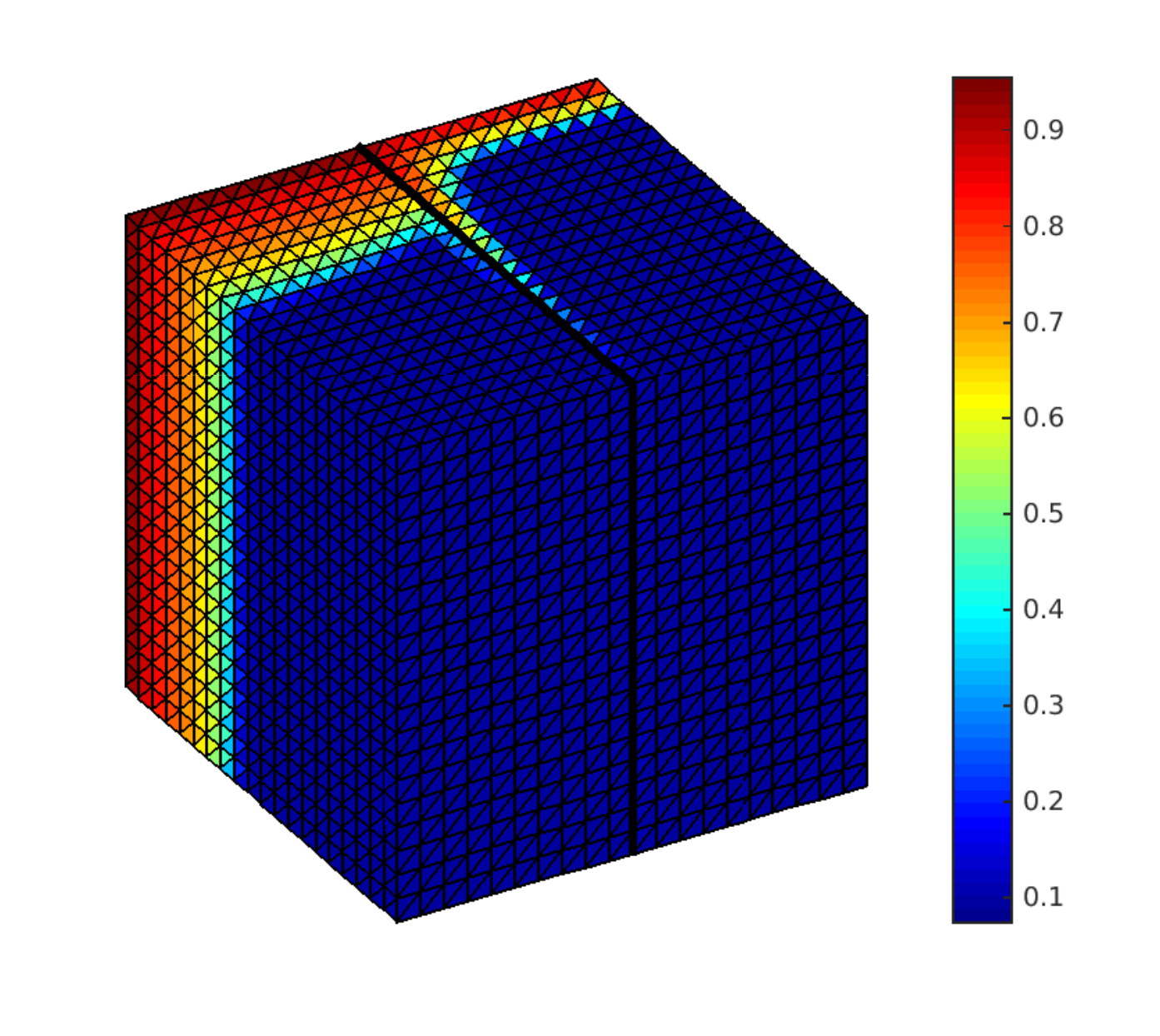}}
{\includegraphics[scale=0.45]{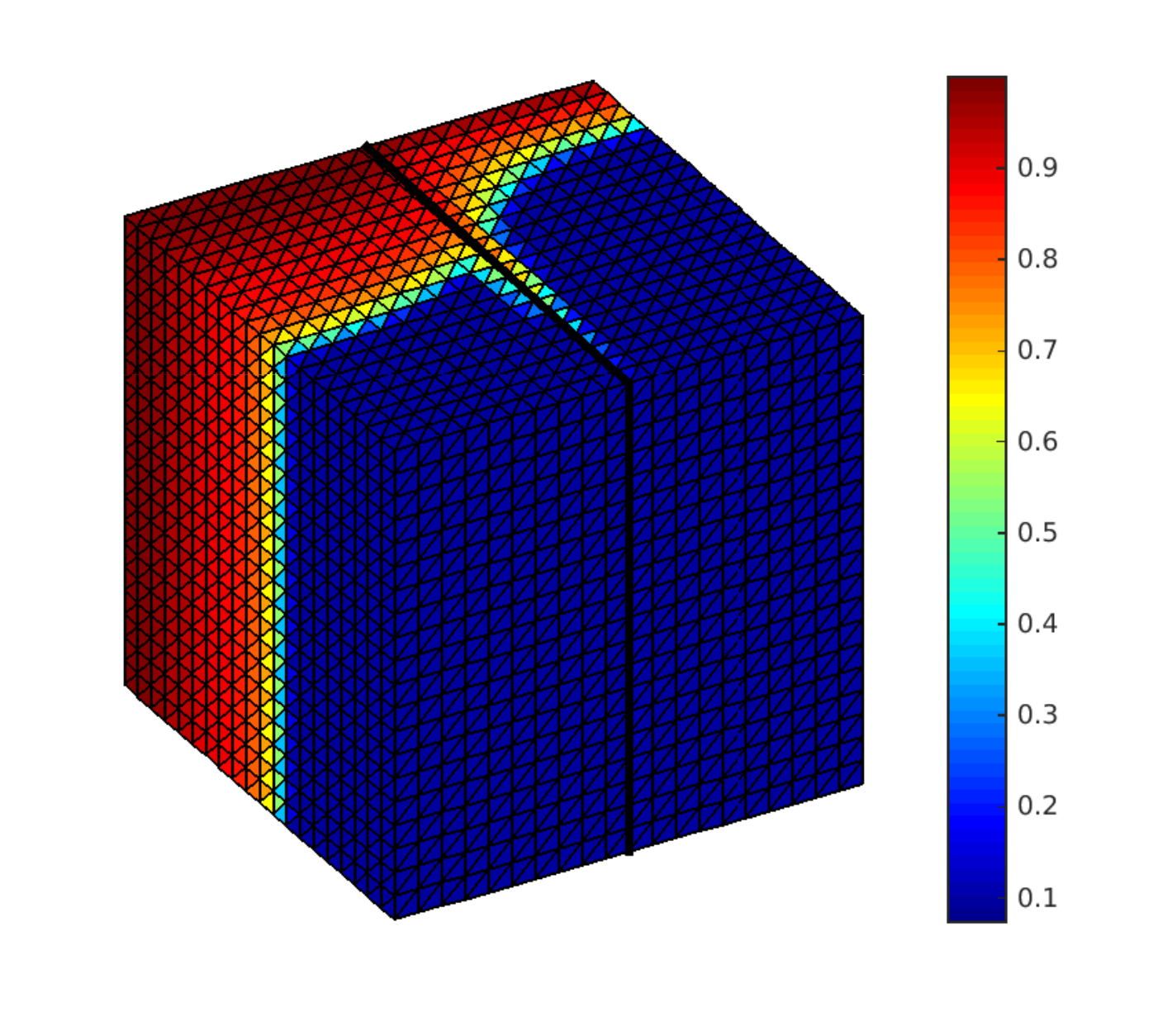}}\\
{\includegraphics[scale=0.45]{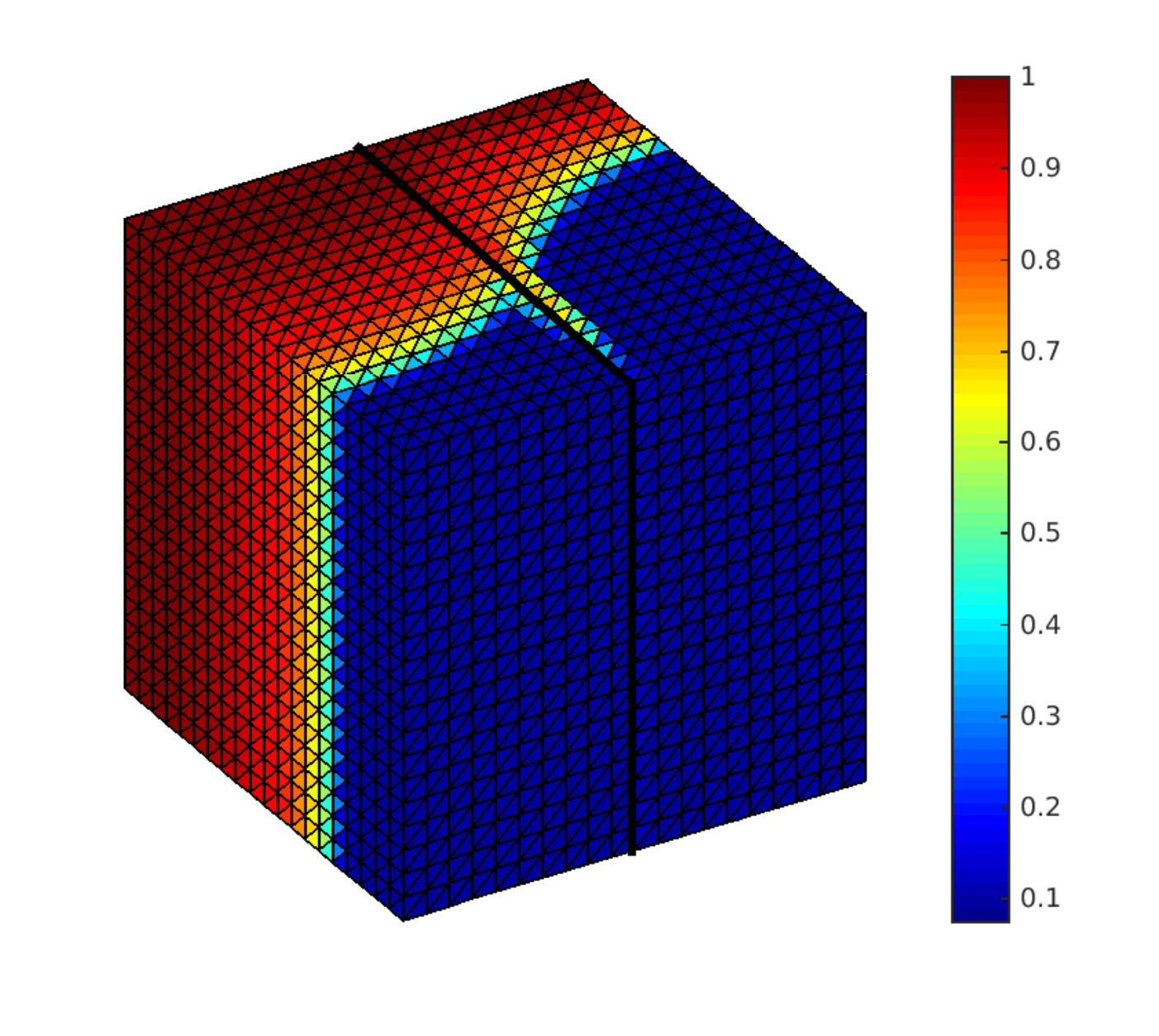}}
{\includegraphics[scale=0.45]{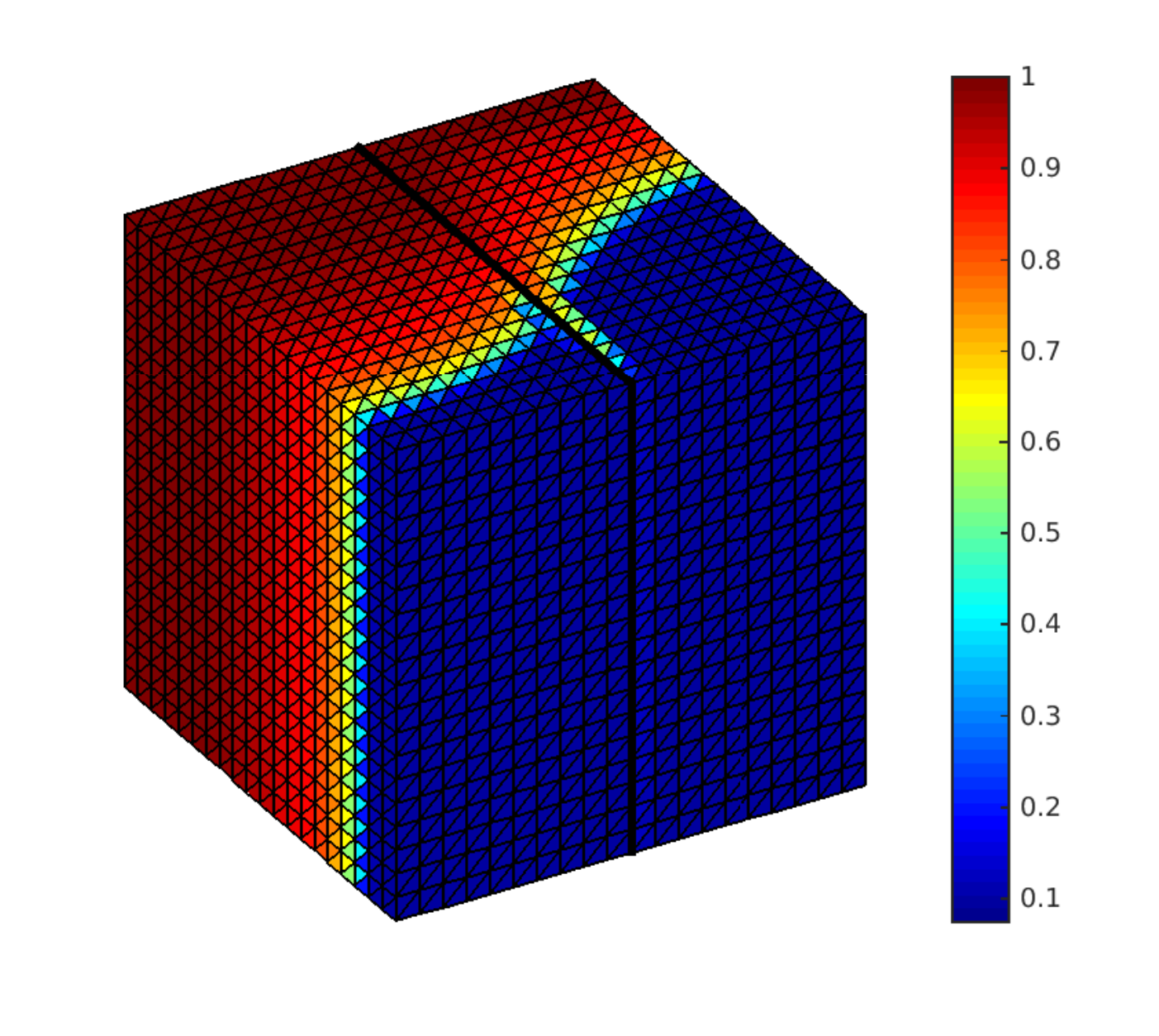}}
\caption{Test case~2: saturation $s(t)$ for $t=5\cdot 10^4$, $t=2\cdot 10^5$, $t=8\cdot 10^5$ and $t=2\cdot 10^6$.}
\label{Two-domain_saturation}
\end{figure}

\begin{figure}[H]\centering
\centering{\includegraphics[scale=0.353]{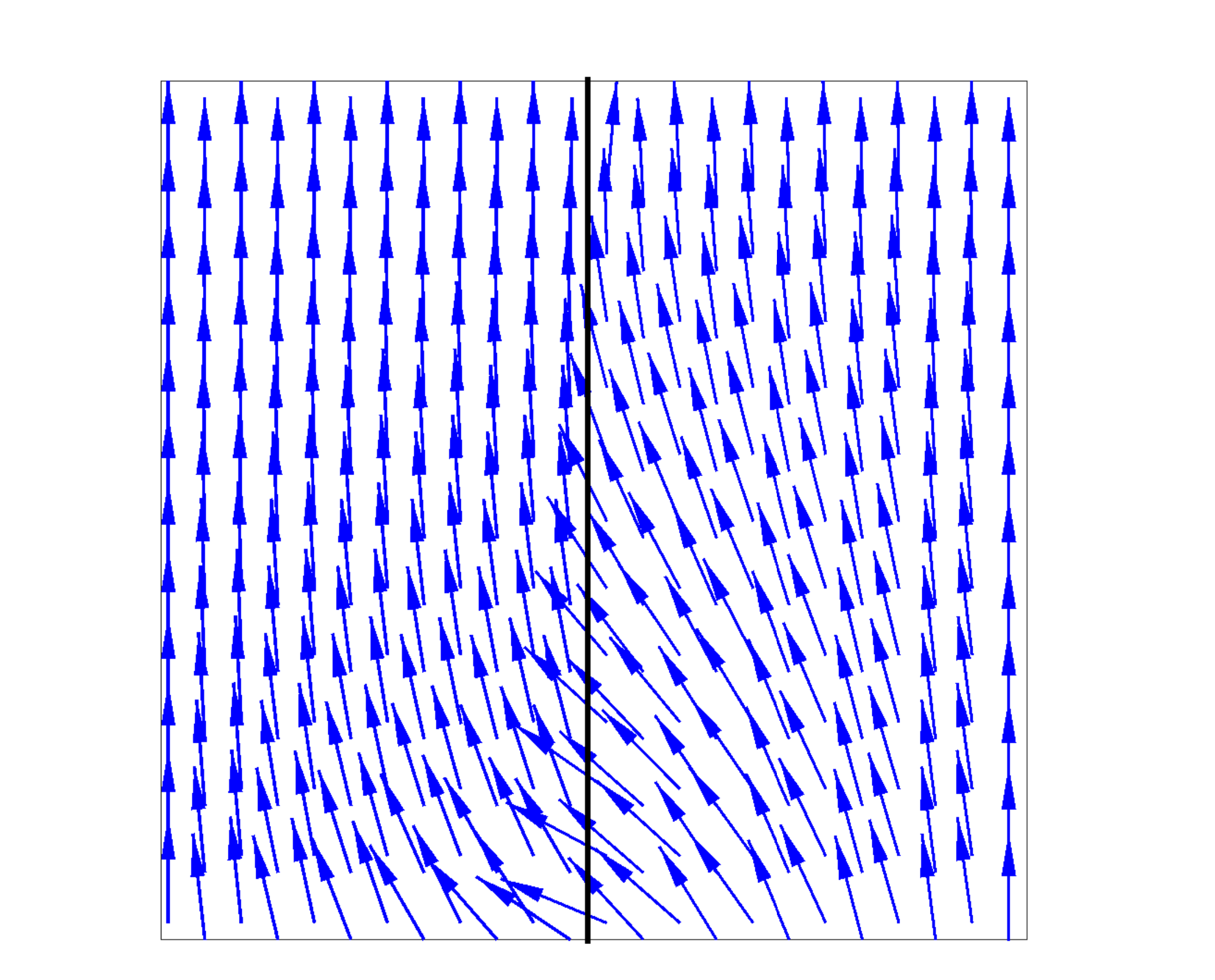}}
\caption{Test case~2: velocity profile between two rock types.} 
\label{two-domain_velocity}
\end{figure}

In Figure~\ref{Two-domain_saturation}, the saturation solution  is depicted for four time steps, 
and we can show that the saturation is discontinuous across the interface so that it  respects the continuity
of the capillary pressure. We can  see also that because of the contrast in the capillary pressure field, the saturation front moves faster around
Rock~2. Thus, the capillary pressure has smoothed out a finger effect and  a spike effect is observed  at the interface. In Figure~\ref{two-domain_velocity}, we show the velocity profile between the rock types. We can see that the fluid inside Rock~2  snakes around the interface  to travel again 
through Rock~1. 

\begin{figure}[h]\centering
{\includegraphics[scale=0.4]{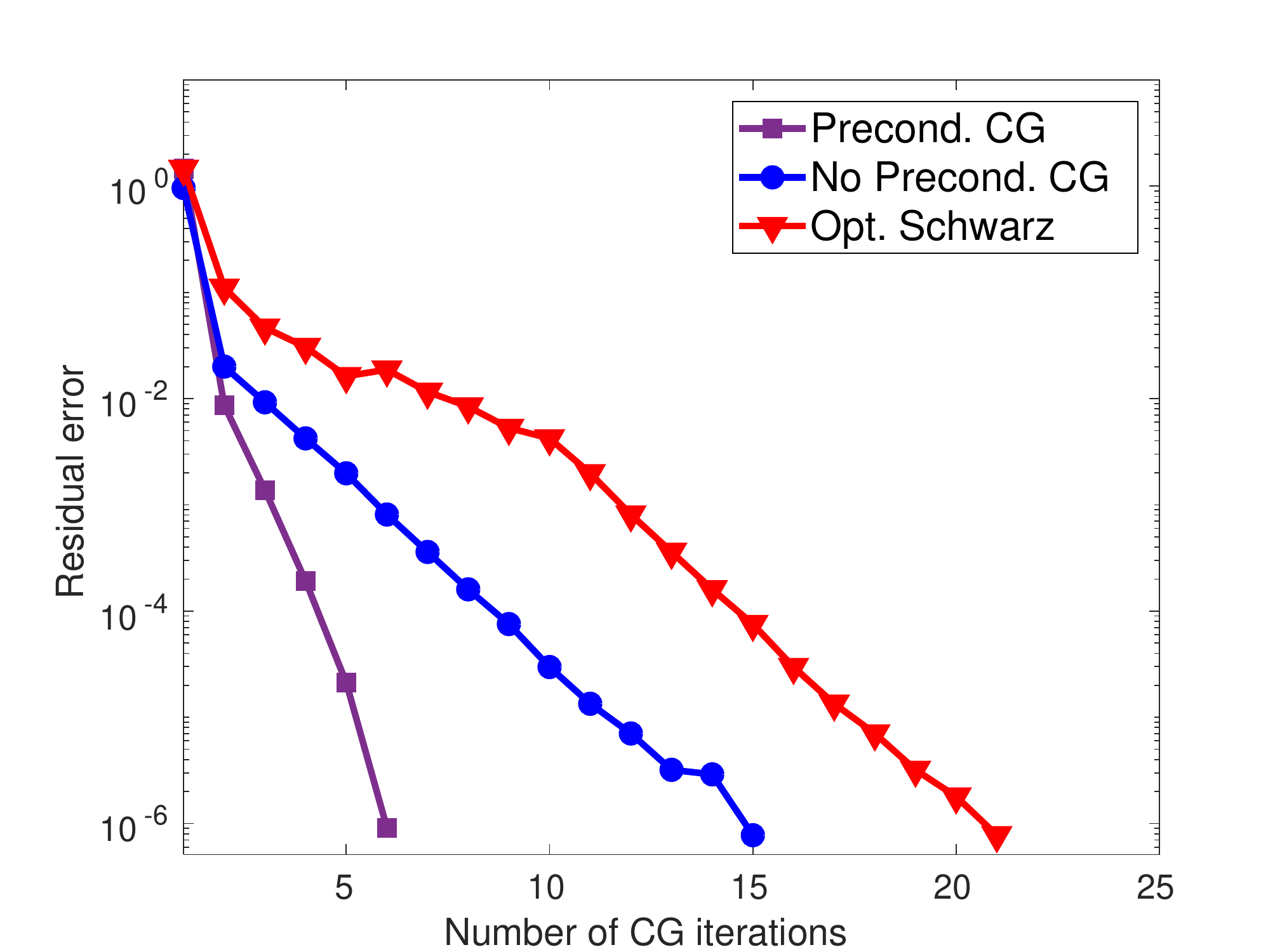}}
{\includegraphics[scale=0.4]{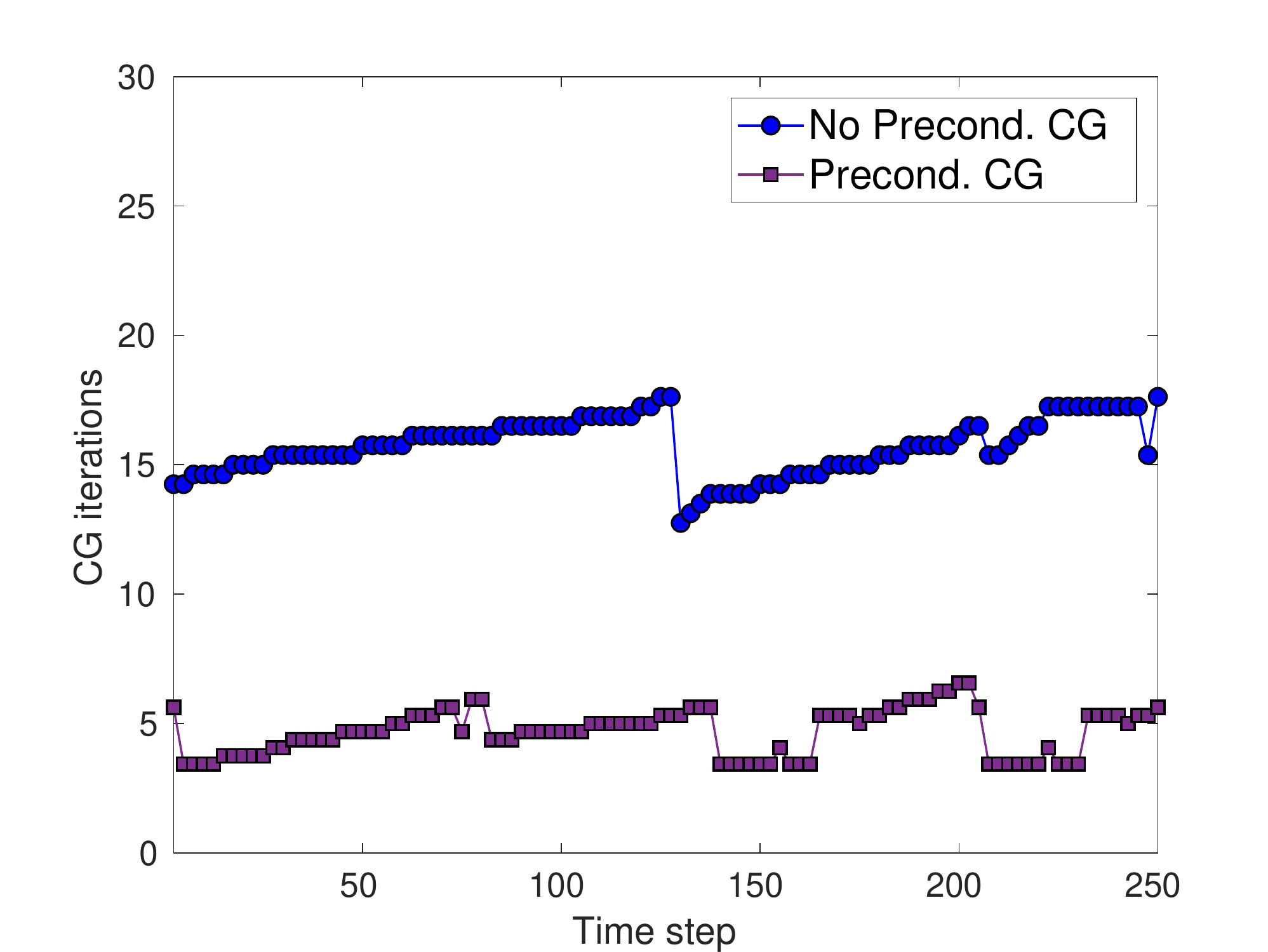}}
\caption{Test case~2: convergence of Interf. CG for $t=2\cdot 10^6$ (left) and  cumulative   number  of Interf. CG iterations (right).}
\label{Convergence_CG_two_rocks}
\end{figure}

\begin{figure}[h]\centering
{\includegraphics[scale=0.38]{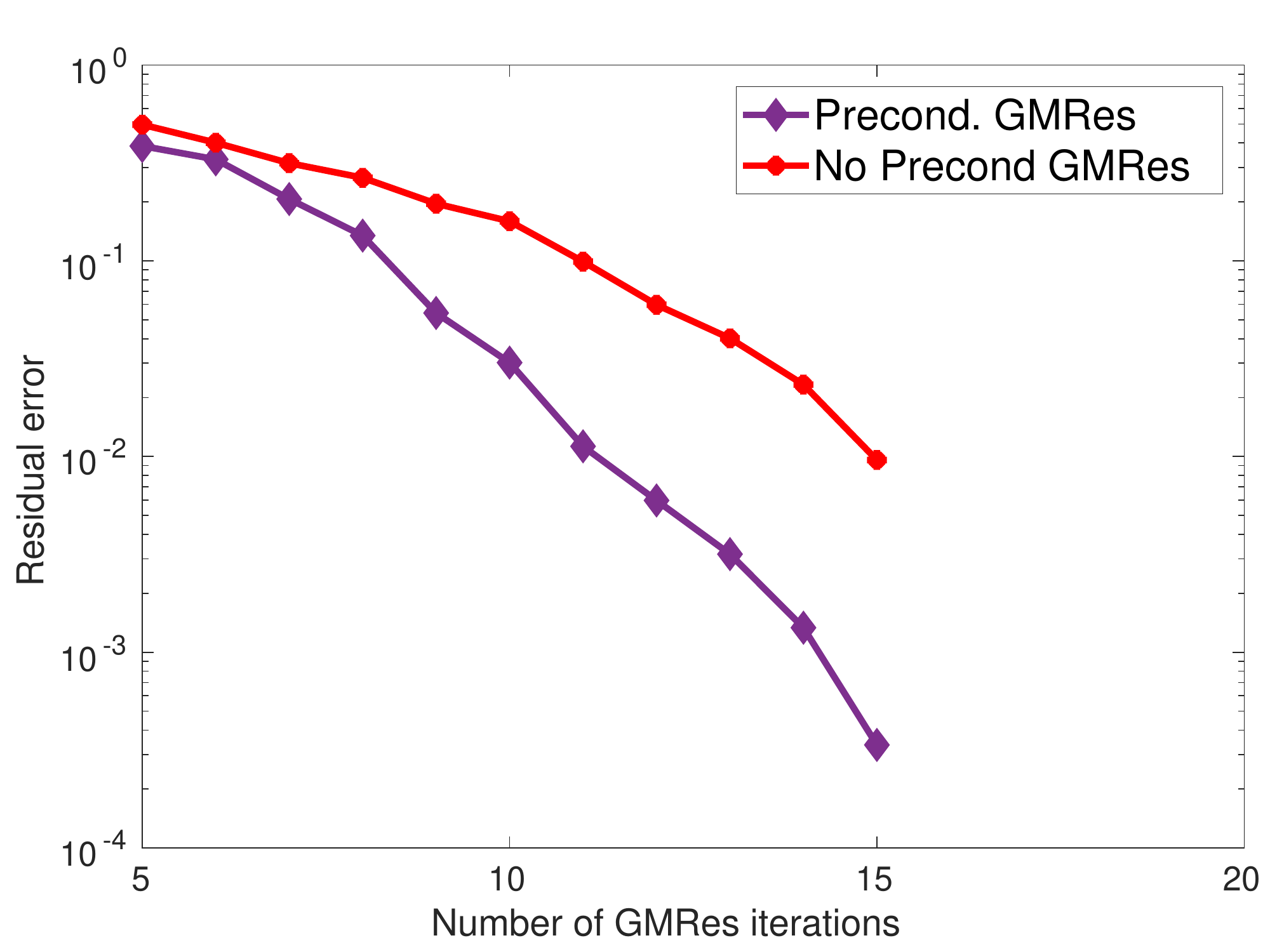}}
{\includegraphics[scale=0.38]{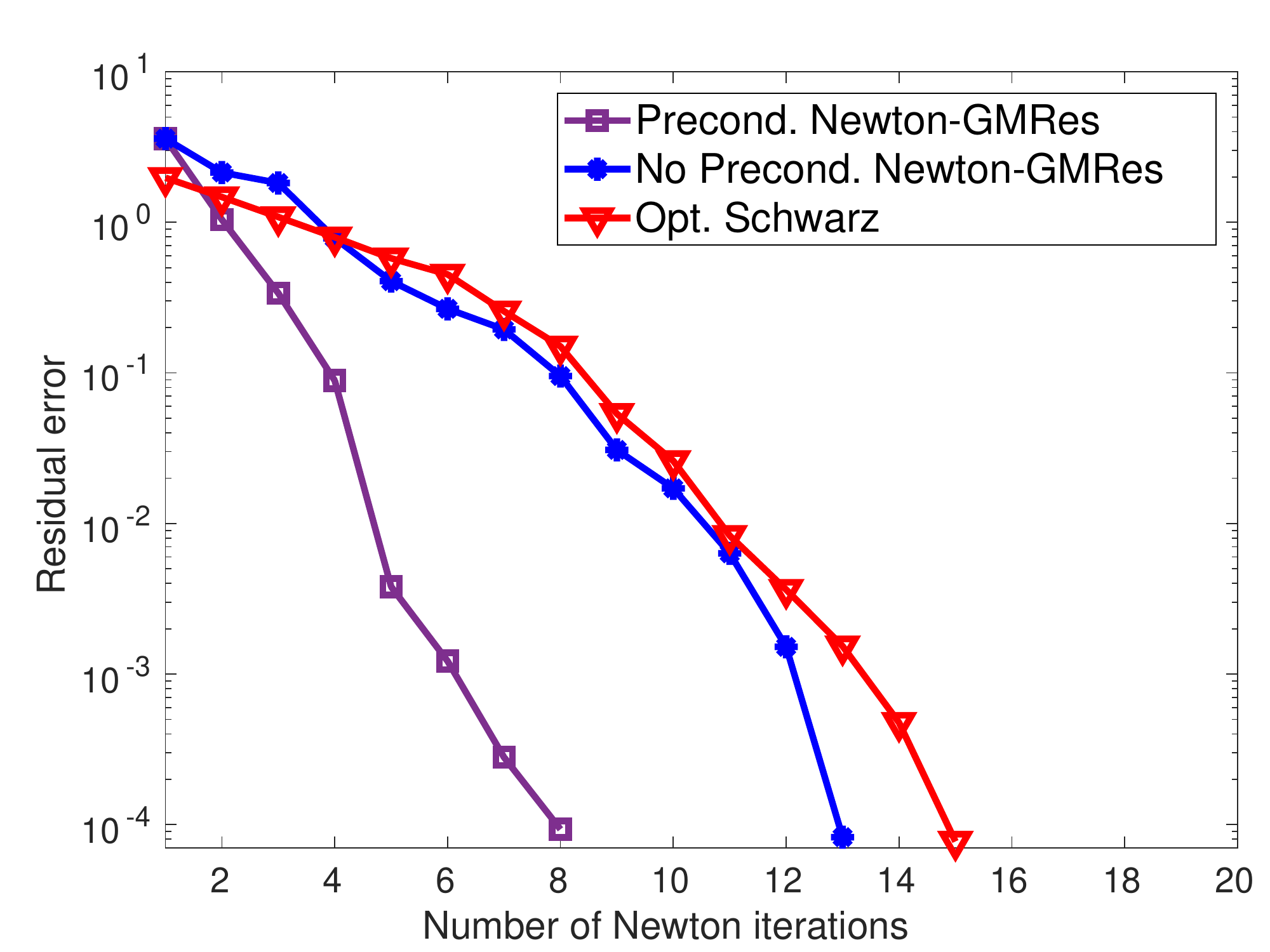}}
\caption{Test case~2: convergence of the Interf. GMRes  (left) and Interf. Newton (right) for $t=2\cdot 10^6$.}
\label{Convergence_Newton_GMRES_two_rocks}
\end{figure}

 We come now to the   convergence analysis  of the iterative solvers involved in   Algorithm~\ref{algo_1}.  We start with the linear ones used in~\texttt{PRESSURE}\_\texttt{SOLVER}.  In  Figure~\ref{Convergence_CG_two_rocks}, we plot  the   residual error for the CG solver  with and without  preconditioning and the OSWR method
 in one fixed time step $t^{n}=2\cdot 10^6$ (left) and  the cumulative  number  of  CG iterations  as  a  function  of  time
(right). One can clearly observe 
that the effect of the  preconditioner on the   convergence of the CG method is significant  and that the  desired residual  tolerance is achieved with  6 iterations. The average number of  CG iterations is reduced from 15.1 to 4.8. That of the OSWR method (Algorithm~\ref{linear_oswr}) is  achieved  after 21 iterations. For the solvers within~the module~\texttt{DIFFUSION}\_\texttt{SOLVER},  the preconditioner of the  GMRes iterations, used in the   Newton-GMRes method,  improves the linear convergence; we obtain the  relative tolerance within 15 iterations (see Figure~\ref{Convergence_Newton_GMRES_two_rocks} (left)).
Note that the reduction in the number of interface iterations is still  small and the average is reduced from 17.81 to 13.87, but 
 faster   Newton convergence is observed during most of the iterations (see Figure~\ref{Convergence_Newton_GMRES_two_rocks} (right)). 
 The OSWR method (Algorithm~\ref{nonlinear_oswr}) for this step  behaves as  in the pressure problem so it converges slower than the Newton-GMRes method. For the  Newton  solver within~\texttt{ADVECTION}\_\texttt{SOLVER}, the cost is negligible as the solution of the interface advection problem~\eqref{Chap2_advection_flux_definit2_Robin} does not require  subdomain solves, which is not the case in the other modules of the algorithm.


\subsubsection{ Comparison of the algorithms:   accuracy-in-time and overall cost}
Due to the nonconformity in time  and the use of different splitting techniques,  we  want to  see  whether 
nonconforming time grids  between subdomains  as well as multirate time steps (within Algorithm~\ref{algo_1} and~\ref{algo_2}) preserve the accuracy in time. 

To study the accuracy of the nonconformity in time between subdomains, we compute with each algorithm  a reference solution on a very fine time grid $(\tau^{n,l}_{i}=\tau^{n}/200=5\cdot10^1$) and a fixed mesh. We then consider four 
initial time grids, refined then  4 times by a factor of 2:
\begin{itemize}
 \item Time grid 1 (Conforming fine): $\tau^{n,l}_{i}=5\cdot10^2$, $i\in\{1,2\}$,
 \item Time grid 2 (Nonconforming, fine in Rock~1): $\tau^{n,l}_{1}=5\cdot10^2$, $\tau^{n,l}_{2}=1\cdot10^3$,
 \item Time grid 3 (Conforming coarse): $\tau^{n,l}_{i}=2\cdot10^3$, $i\in\{1,2\}$,
 \item Time grid 4 (Nonconforming, coarse in Rock~1): $\tau^{n,l}_{1}=2\cdot10^3$, $\tau^{n,l}_{2}=2\cdot10^2$.
\end{itemize}

\begin{figure}[h]\centering
\centering{\includegraphics[scale=0.3]{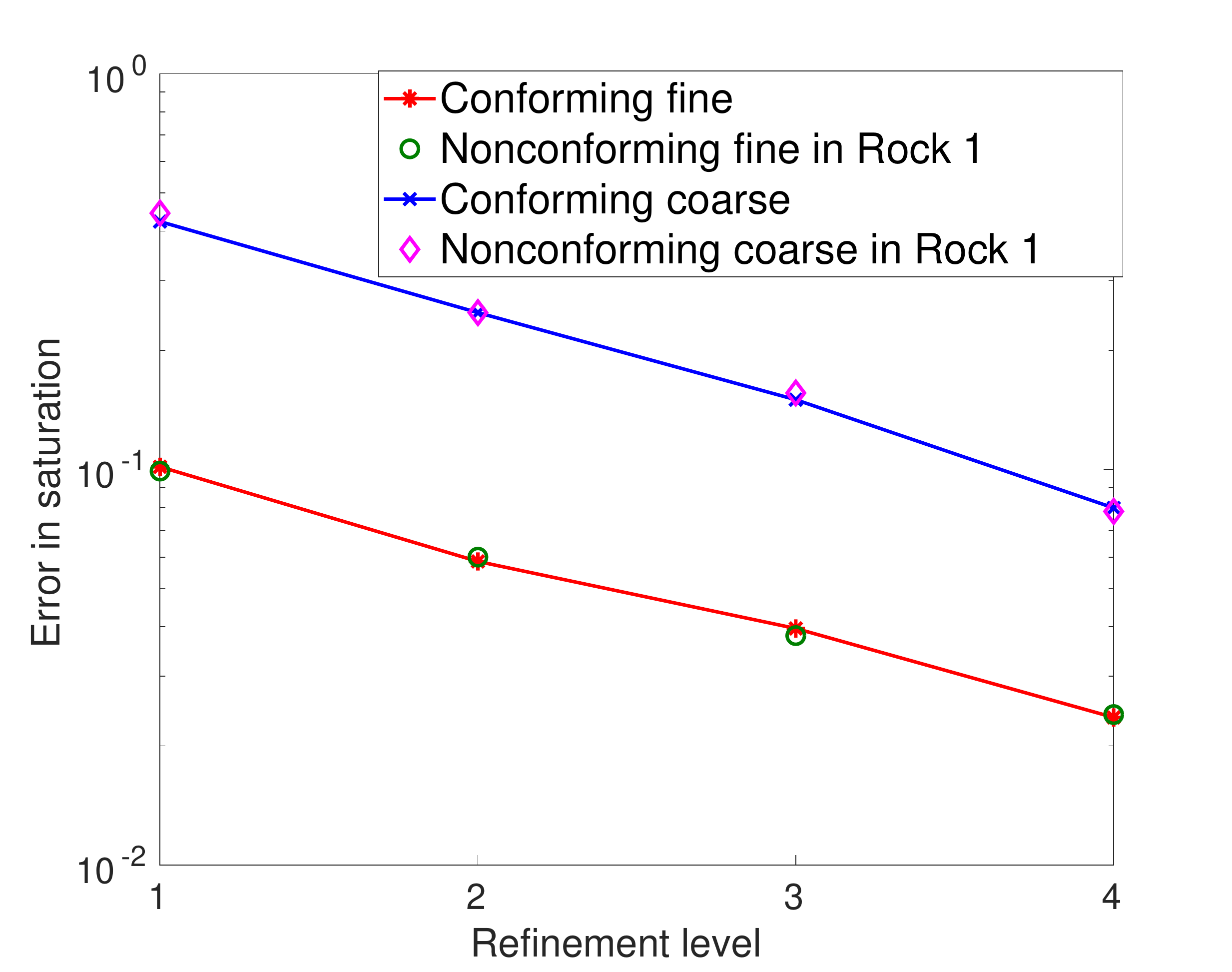}
\includegraphics[scale=0.3]{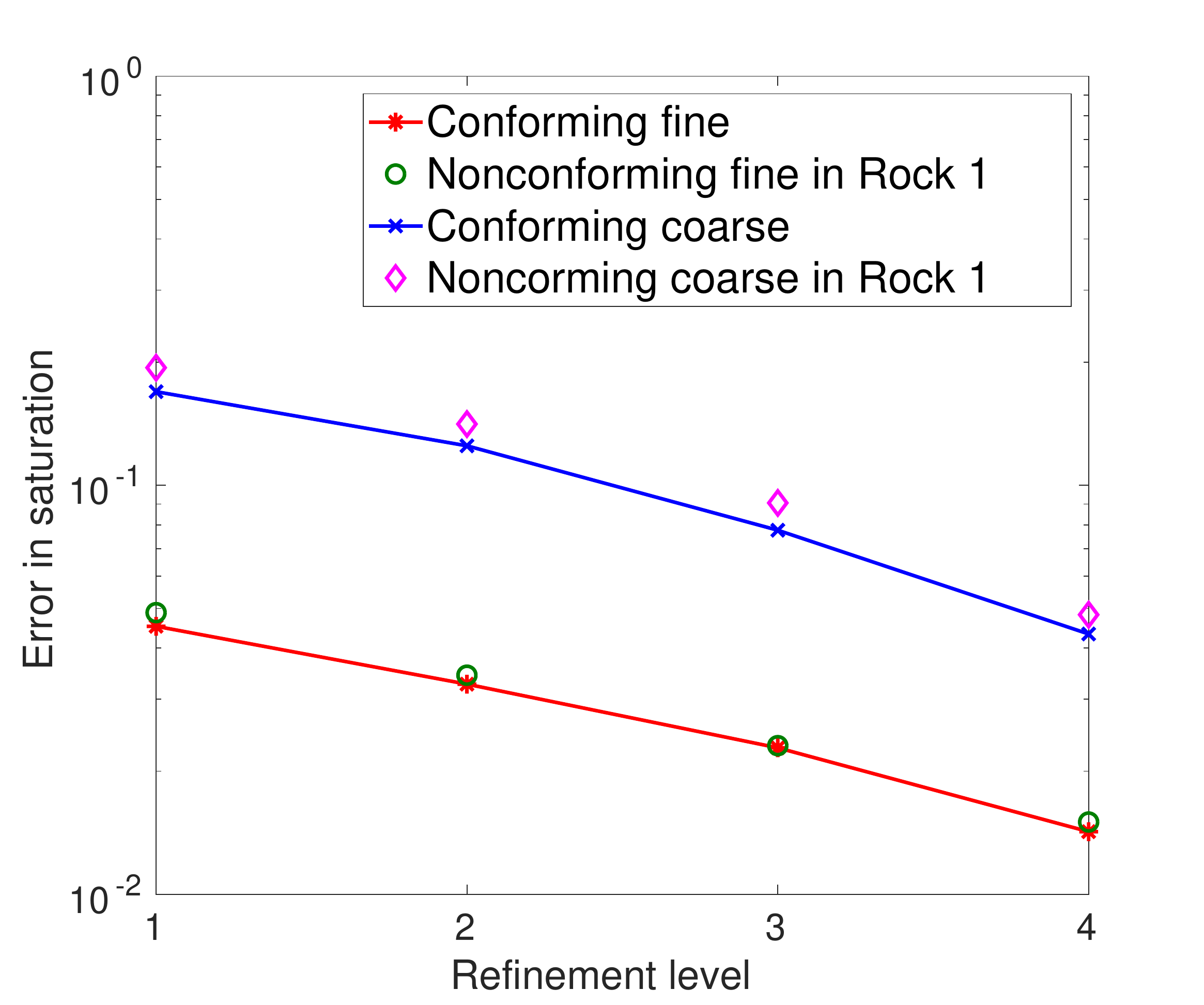}
}
\caption{Test case~2: errors in
 saturation between the reference  and the numerical solutions versus the refinement level for Algorithm~\ref{algo_1} (left) and Algorithm~\ref{algo_2} (right).}

\label{Error_saturation}
\end{figure}

In Figure~\ref{Error_saturation},  the error in the $L^2(0,T;L^2(\Omega))$-norm of the saturation
versus the refinement level is depicted for the two algorithms. Clearly, both algorithms preserve the accuracy in time as the errors 
obtained in the nonconforming case with a fine time step in Rock~1 coincides with those obtained with the finer conforming case. 

 We come now to the comparison of Algorithm~\ref{algo_1} and Algorithm~\ref{algo_2} as well as the assessment of the accuracy in time of the multirate time strategies. In that case, we consider conforming time grids in the two algorithms. We compute using  Algorithm~\ref{algo_2} a reference solution on a very fine time grid $(\tau^{n,l}_{i}=\tau^{n}=2\cdot10^1)$ and a fixed mesh. We then test the two algorithms with
$\tau^{n}=\mathcal{N}^{\textnormal{a}}\tau^{n,l}_{i}$, with a  factor $\mathcal{N}^{\textnormal{a}}=48$ (the number of advection fine time steps within one coarse pressure and diffusion time step), divided then 4 times by a factor of 2. 

\begin{figure}[h!]\centering{
\includegraphics[scale=0.35]{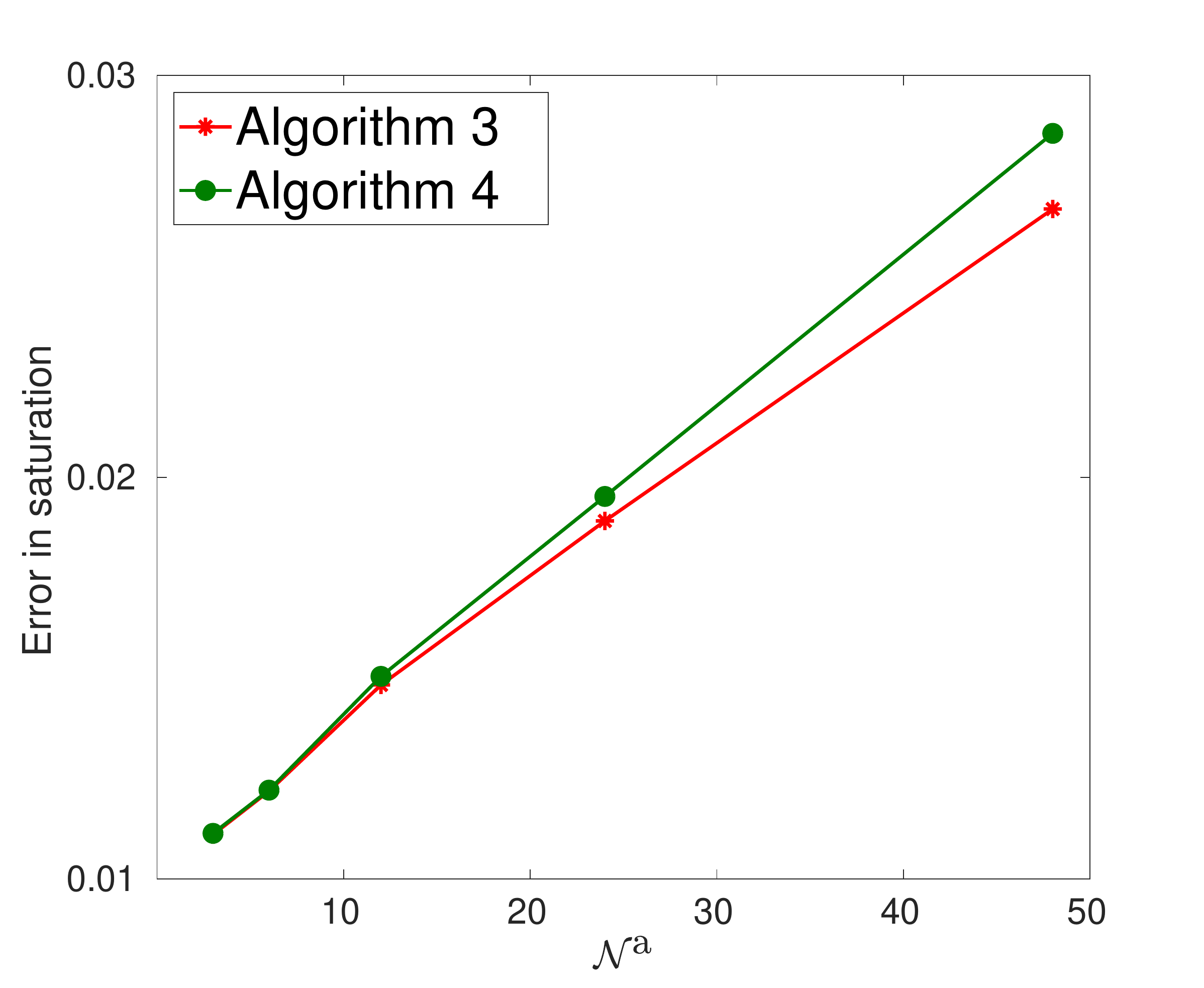}
}
\caption{Test case~2: errors in
 saturation between the reference  and the numerical solutions versus the multirate factor $\mathcal{N}^{\textnormal{a}}$.}
\label{Error_saturation-vsspliiting}
\end{figure}
In Figure~\ref{Error_saturation-vsspliiting},  we show the error in the $L^2(0,T;L^2(\Omega))$-norm of the saturation
versus the splitting level for the two algorithms.
Clearly, an excellent  quality of the solution is obtained from both algorithms even with $\mathcal{N}^{\textnormal{a}}=48$. Particularly, as stated previously Algorithm~\ref{algo_1}  decreases the number of inner time steps and  Newton iterations required for the full saturation  problem in Algorithm~\ref{algo_2}. As a result, the excellent  quality of the solution from  
 Algorithm~\ref{algo_1} and its reduced cost compared to Algorithm~\ref{algo_2}  as well the  reduced exchanged data (also  smaller subdomain solves are required) between the different rocks make Algorithm~\ref{algo_1} an efficient tool to deal with two-phase  flow model in heterogeneous media.
  \subsection{Test case~3: extensions to~multiple subdomains}
Based on the above comparison, we choose here   to apply Algorithm~\ref{algo_1}  on more complex configurations of porous media made-up with different rock types. For the next two cases, we let $\mathrm{\Omega}=[0,10]^{3}$ and $T=5\cdot 10^6$s. 
 \subsubsection{Fractured porous media}
\begin{figure}[h!]
\begin{center}
\includegraphics[scale=0.3]{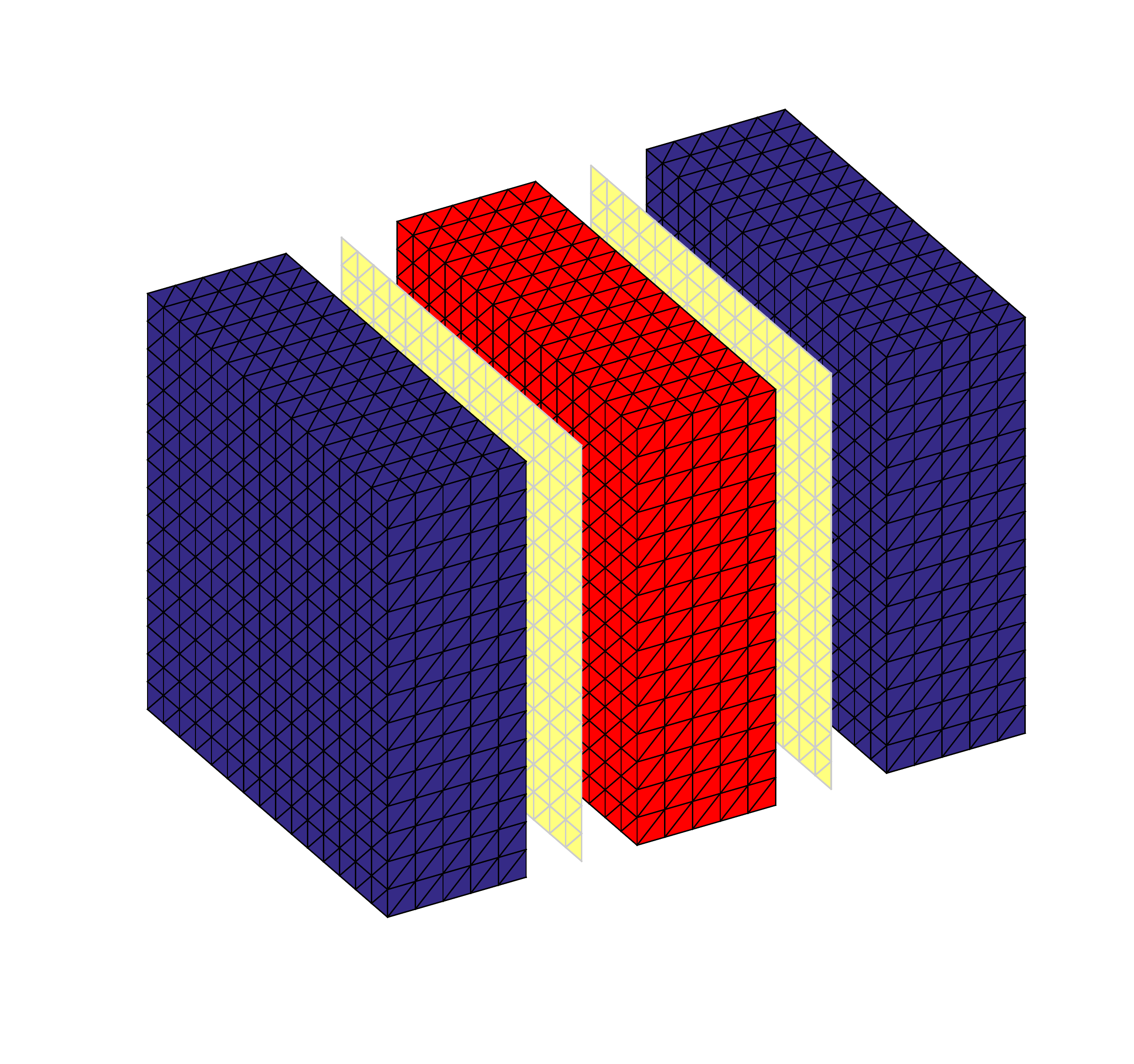}
\caption{Test case~3 (fractured): an exploded view of the division into three subdomains with  two rock types.}
\label{Three-domains_figure}
\end{center}
\end{figure}

In the first configuration, we  consider the  test  case  depicted  in
Figure~\ref{Three-domains_figure} (left)  where   a  three dimensional domain  is  divided  into two  equally  sized  subdomains  by  a  fracture  of  the same size. Precisely, Rock~1 (rock matrix) appears in the left  and right five
layers of elements, separated by  Rock~2 (fracture) in the center
also five layers of elements. 

 In  this example, the properties of
 Rock~1 are the same as the previous test case so that $K_1=1$  and $\phi_{1}=0.3$. In Rock~2,  we fix the absolute permeability  to be  ten times larger than that of Rock~1.  The coarser time steps are of fixed size and 
 $\tau^{n}=1\cdot 10^{3}$s. The time step for the advection in the fracture is taken   five times smaller than  used in  the surrounding rocks, i.e., $\tau^{n,l}_{1}=\tau^{n}/100$ and  $\tau^{n,l}_{2}=\tau^{n}/20$. Water is injected uniformly through one side perpendicular to the fracture~(see~Figure~\ref{Three-domains_figure}). The production boundary is the opposite side. The other boundaries are impermeable.

\begin{figure}[!]\centering
{\includegraphics[scale=0.4]{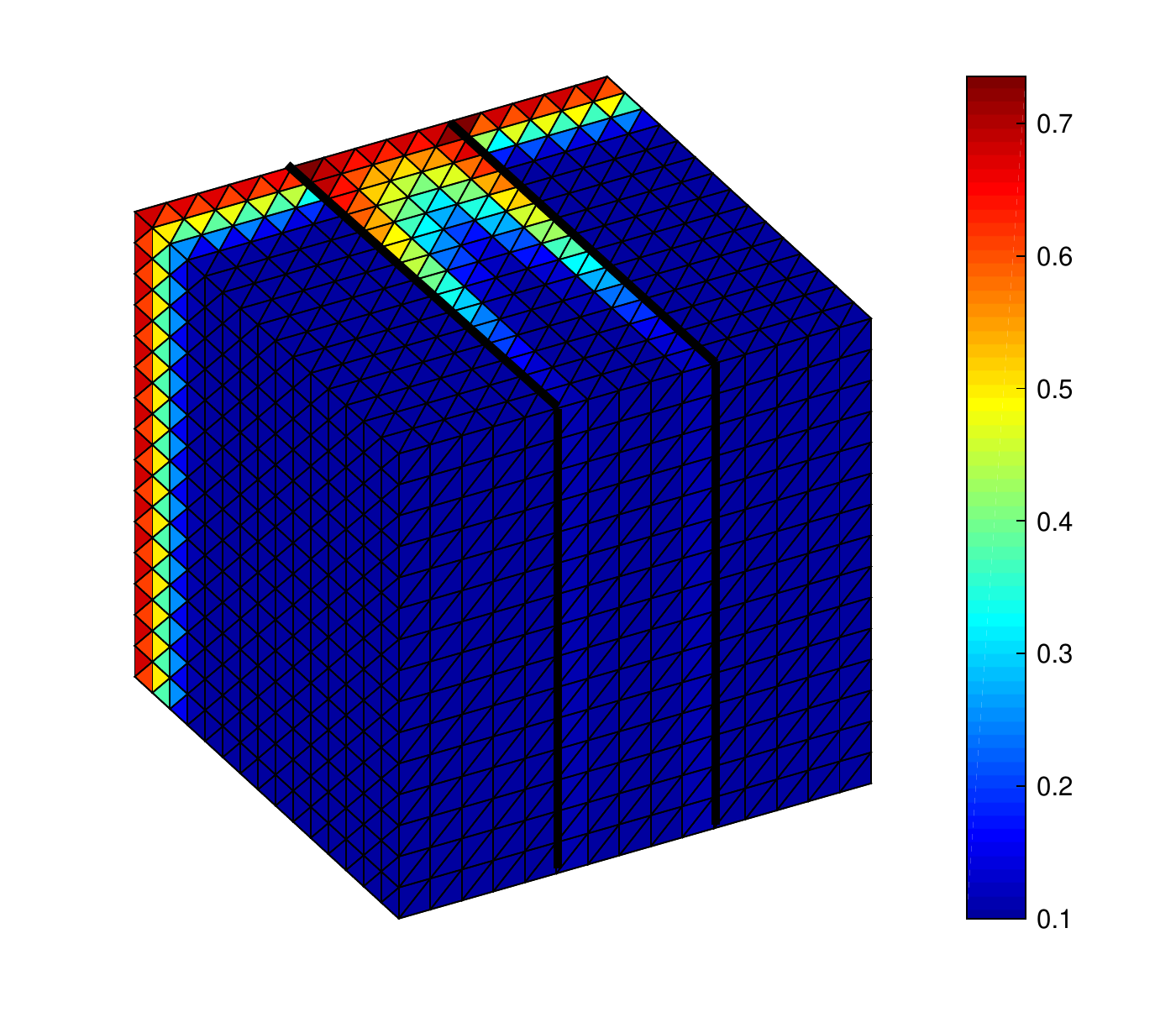}}
{\includegraphics[scale=0.4]{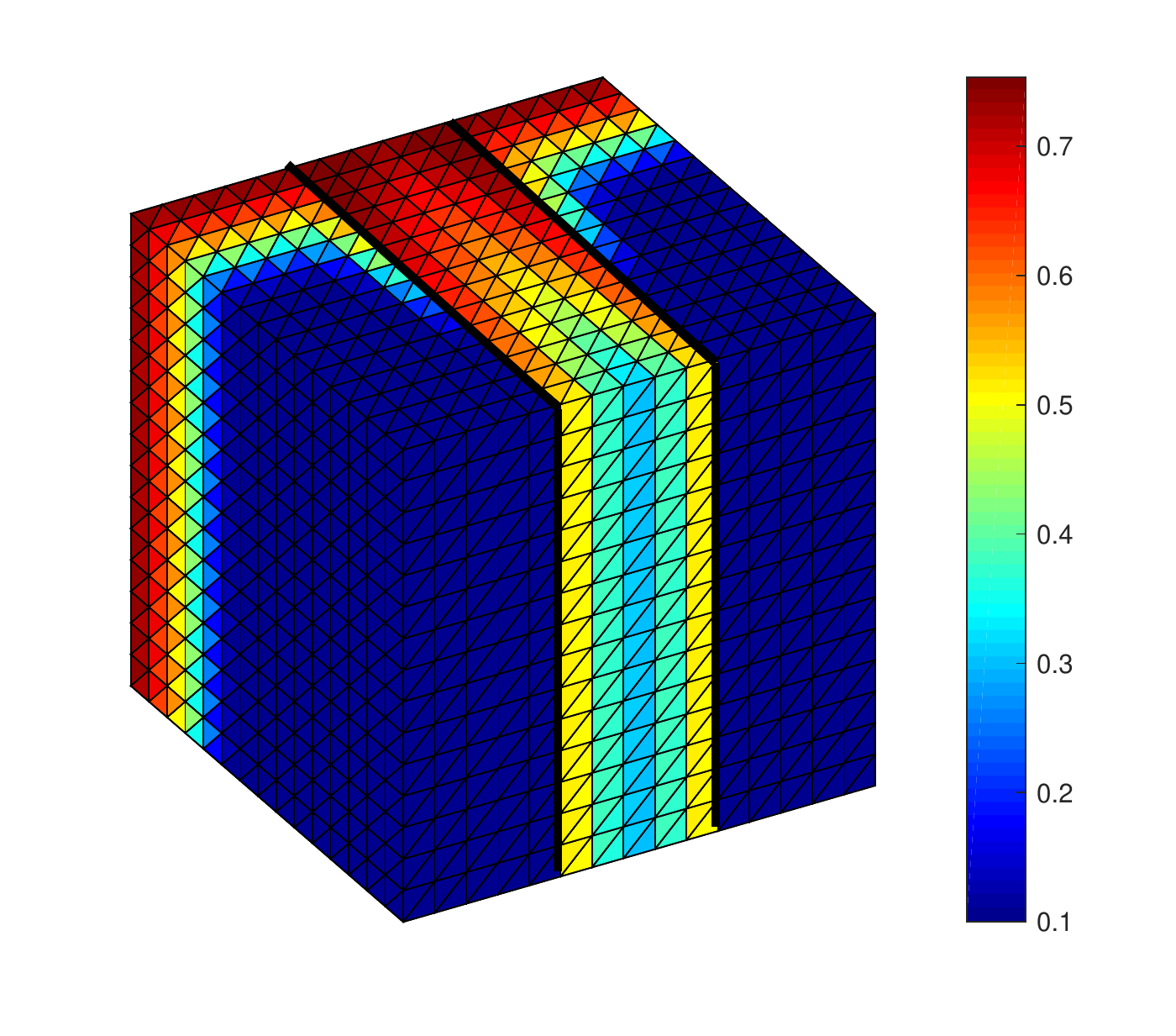}}\\
{\includegraphics[scale=0.4]{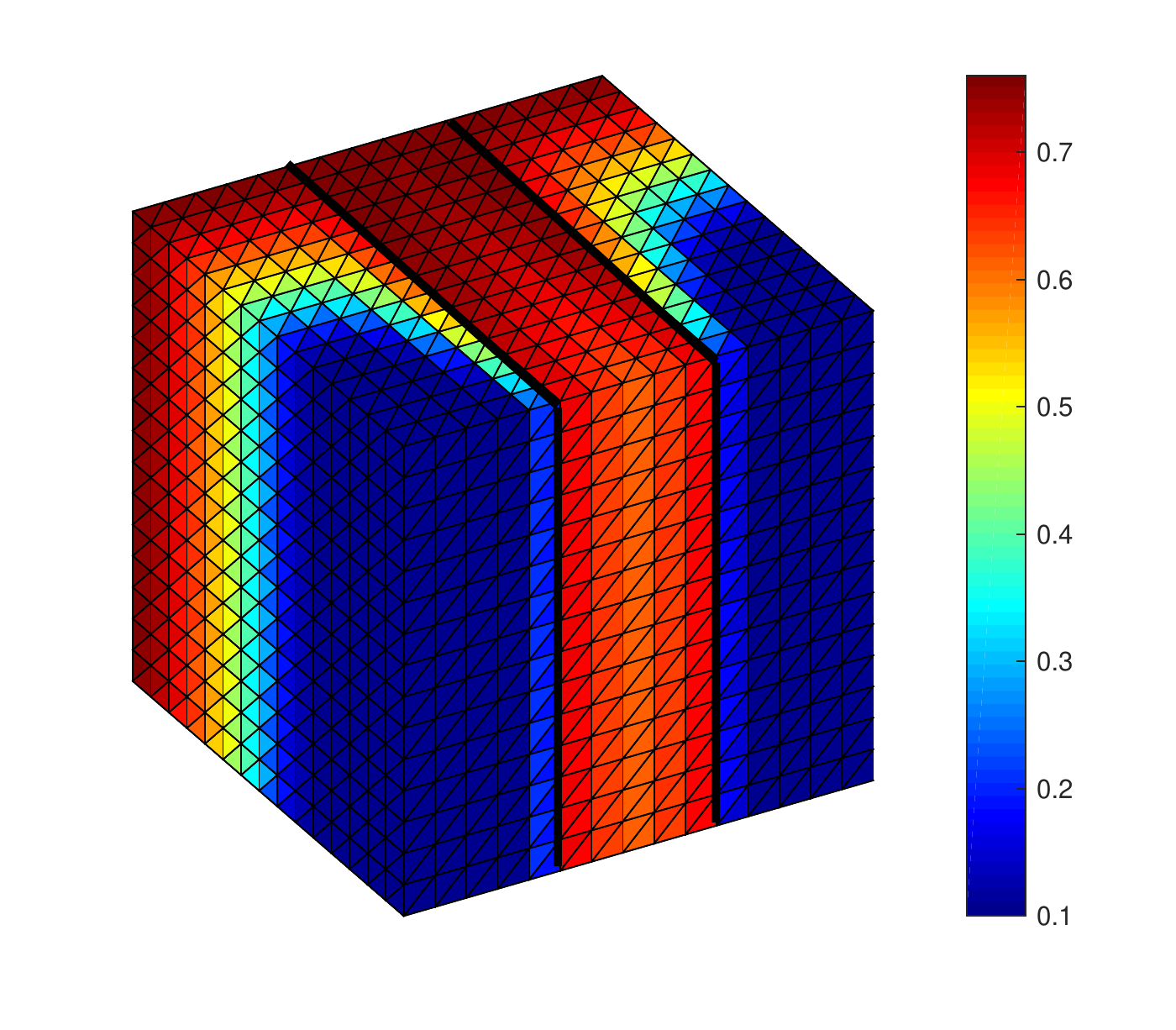}}
{\includegraphics[scale=0.4]{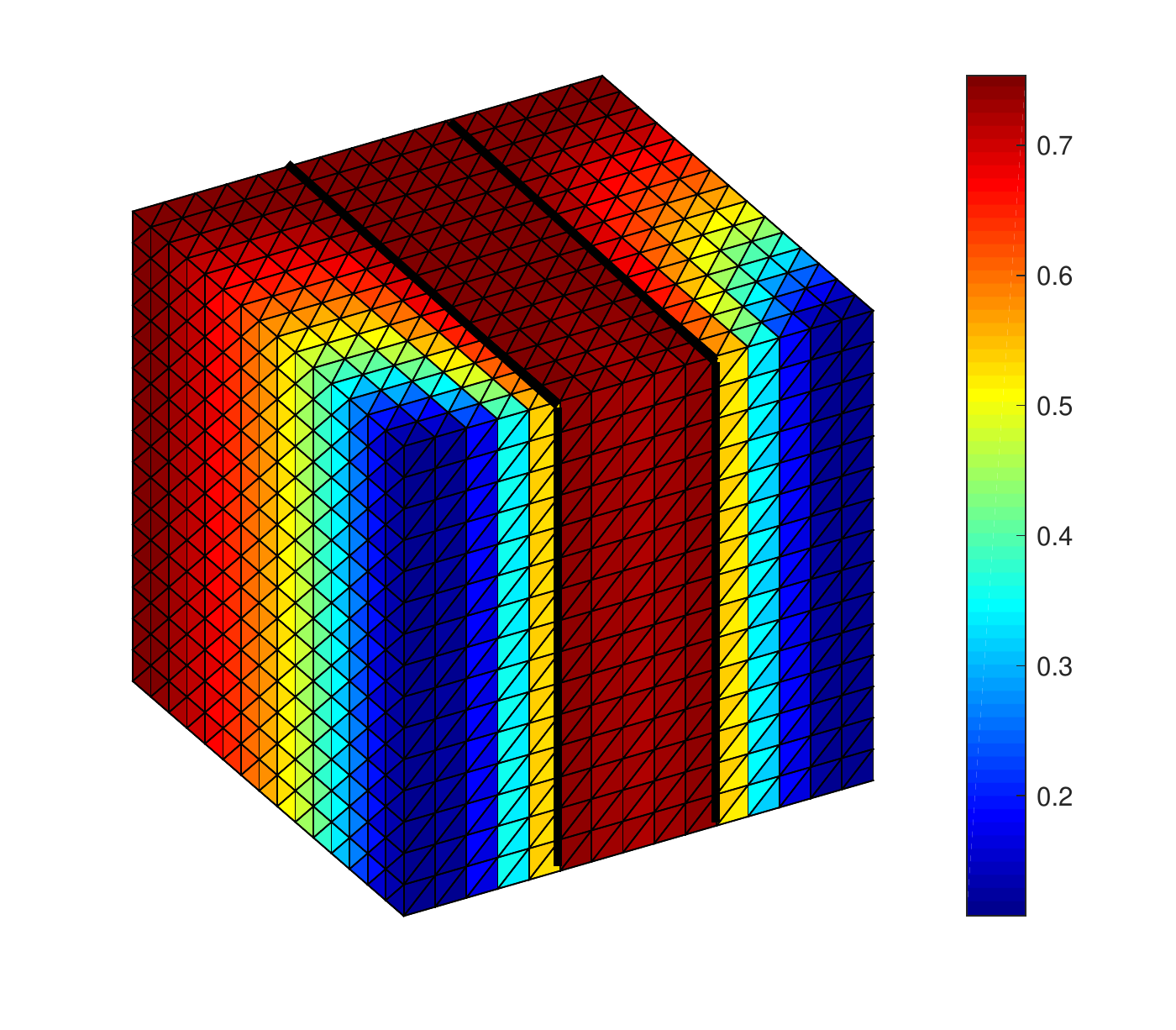}}
\caption{Test case~3 (fractured):  saturation $s(t)$ for $t=2\cdot 10^2$, $t=4\cdot 10^3$, $t=8\cdot 10^3$ and $t=2\cdot 10^4$.}
\label{Three-domain_saturation}
\end{figure}

Plots of the simulation results are shown in Figure~\ref{Three-domain_saturation}. The results illustrate the   discontinuous behavior of the saturation at rock  interfaces due to the strong capillarity effects  present in the fracture.  The water saturation  front  snakes around the fracture to travel through the fracture and then moves from the injection boundary to the  production boundary   of  the fracture. The convergence 
 results (not shown here) are similar to what is  observed for the previous test case, confirming the ability  of Algorithm~\ref{algo_1} to deal with more complex configuration of porous media. Particularly,    the efficiency of the preconditioners  to improve the convergence rate of the solvers compared to the ones without preconditioners is also confirmed in this test case.

\subsubsection{Multiple rock types with strong heterogeneity}
\begin{figure}[h!]
 \begin{center}\includegraphics[scale=0.4]{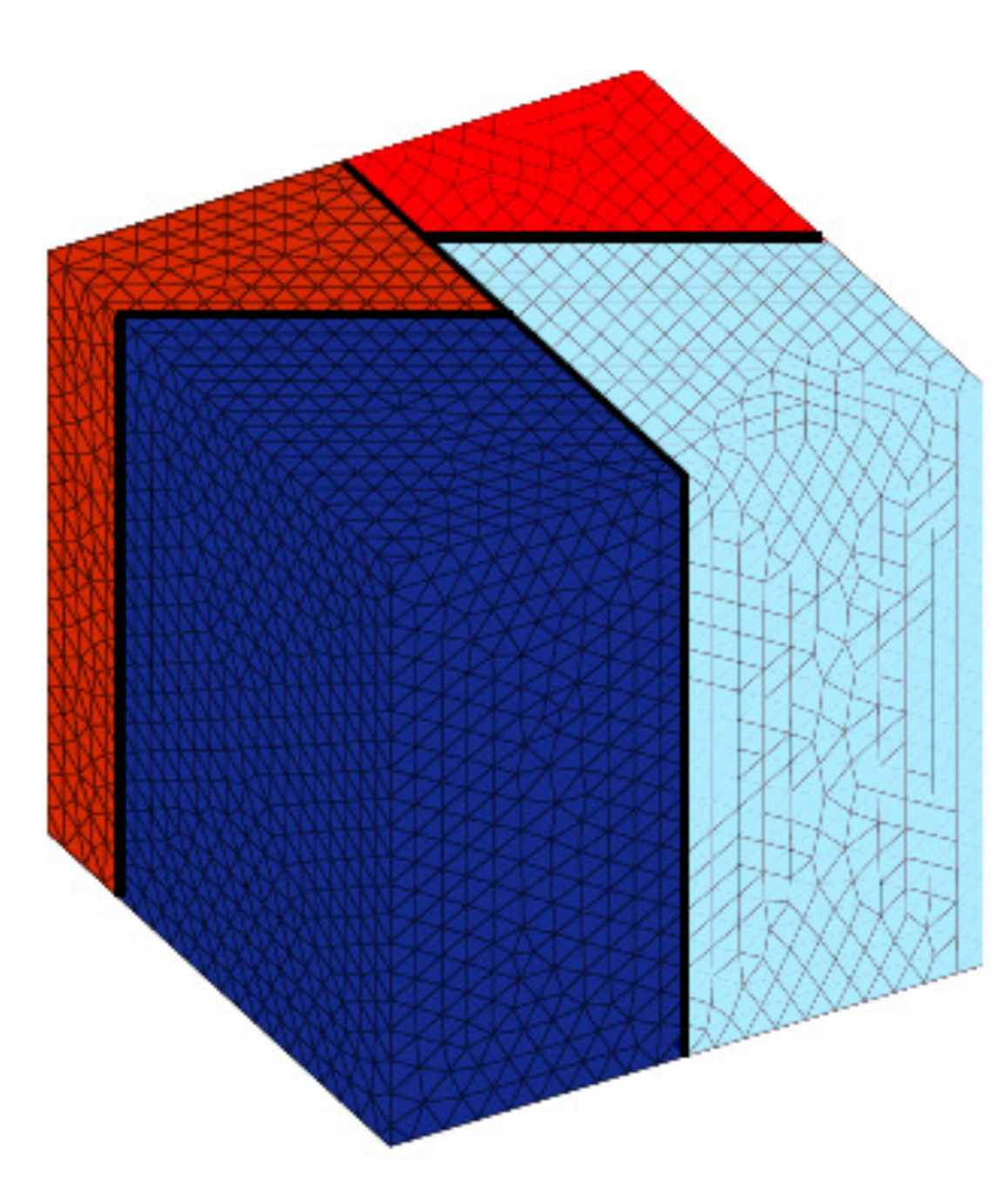}
 \caption{Test case~3 (multiple rocks): The mesh for multiples rock types.}
 \label{four-domains_figure}
 \end{center}
 \end{figure}

In this  test case, we test the capability of  Algorithm~\ref{algo_1} to deal 
with multiple rock types with highly contrasting physical parameters. Precisely, we consider a domain made up with of four rock types (see Figure~\ref{four-domains_figure}).  The rock properties are listed in Table~\ref{table_lasttest}; we can see that Rock~1 is the most permeable and Rock~4 is the one with the lowest permeability.  We also fix $\alpha$ in order to make an advection-dominated problem (with negligible capillarity effects) within Rock~1 and Rock~2,    and  in Rock~3 and Rock~4 we neglect totally the advection effects. In these  rocks, we have only a saturation-diffusion problem as in Test case~1.  

For the initial condition, we assume that the domain  contains some quantity of water situated only within Rock~1 and Rock~2; we set  $s=0.95$ 
in Rock~1 and elsewhere is set to satisfy the continuity of the capillary pressure.   The time scales for the advection 
are taken different in  the different rocks depending on the distribution of the absolute permeability, i.e., $\tau^{n,l}_{1}=\tau^{n}/10$,   $\tau^{n,l}_{2}=\tau^{n}/5$, and $\tau^{n,l}_{3}=\tau^{n,l}_{4}=\tau^{n}/2$   
 with  coarse and fixed  time steps for diffusion and pressure, i.e., $\tau^{n}=\tau^{n}=5\cdot 10^{2}$s. 

\begin{table}[h]
\centering
\begin{tabular}{ | c || c | c |c | } \hline
 & $K$ (md)  & $\phi$ & $\alpha$  \\ \hline
{Rock 1} & 6  &  0.5 & 0.3 \\ \hline
{Rock 2} & 3 & 0.5 & 0.3 \\ \hline
{Rock 3} & 0.6 &0.3  & 5\\ \hline
{Rock 4} & 0.3 &  0.3 & 5\\ \hline
\end{tabular}%
\caption{Test case~3 (multiple rocks): the  physical properties.}\label{table_lasttest}%
\end{table}

\begin{figure}[h]\centering{
\includegraphics[scale=0.4]{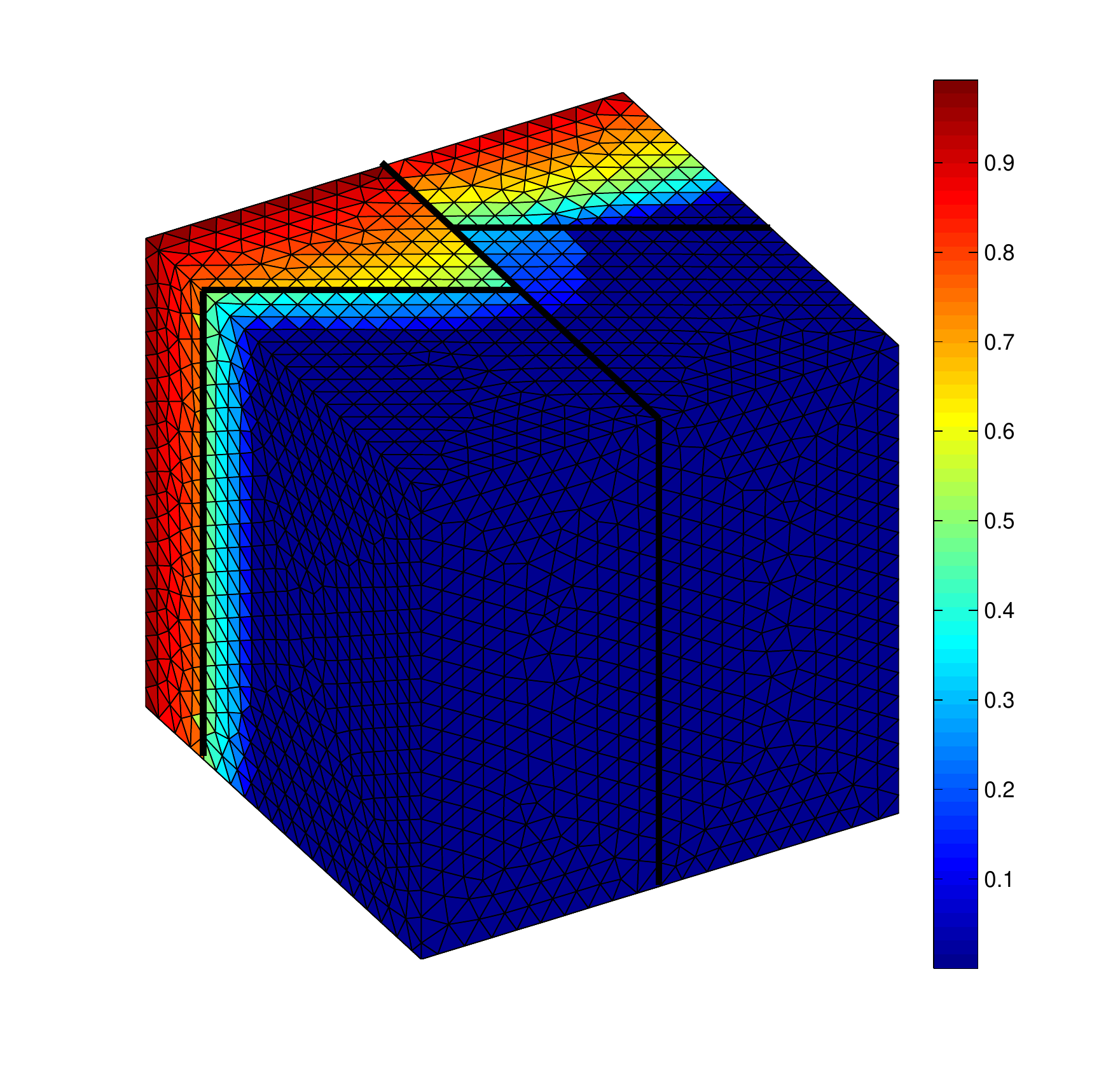}
\includegraphics[scale=0.4]{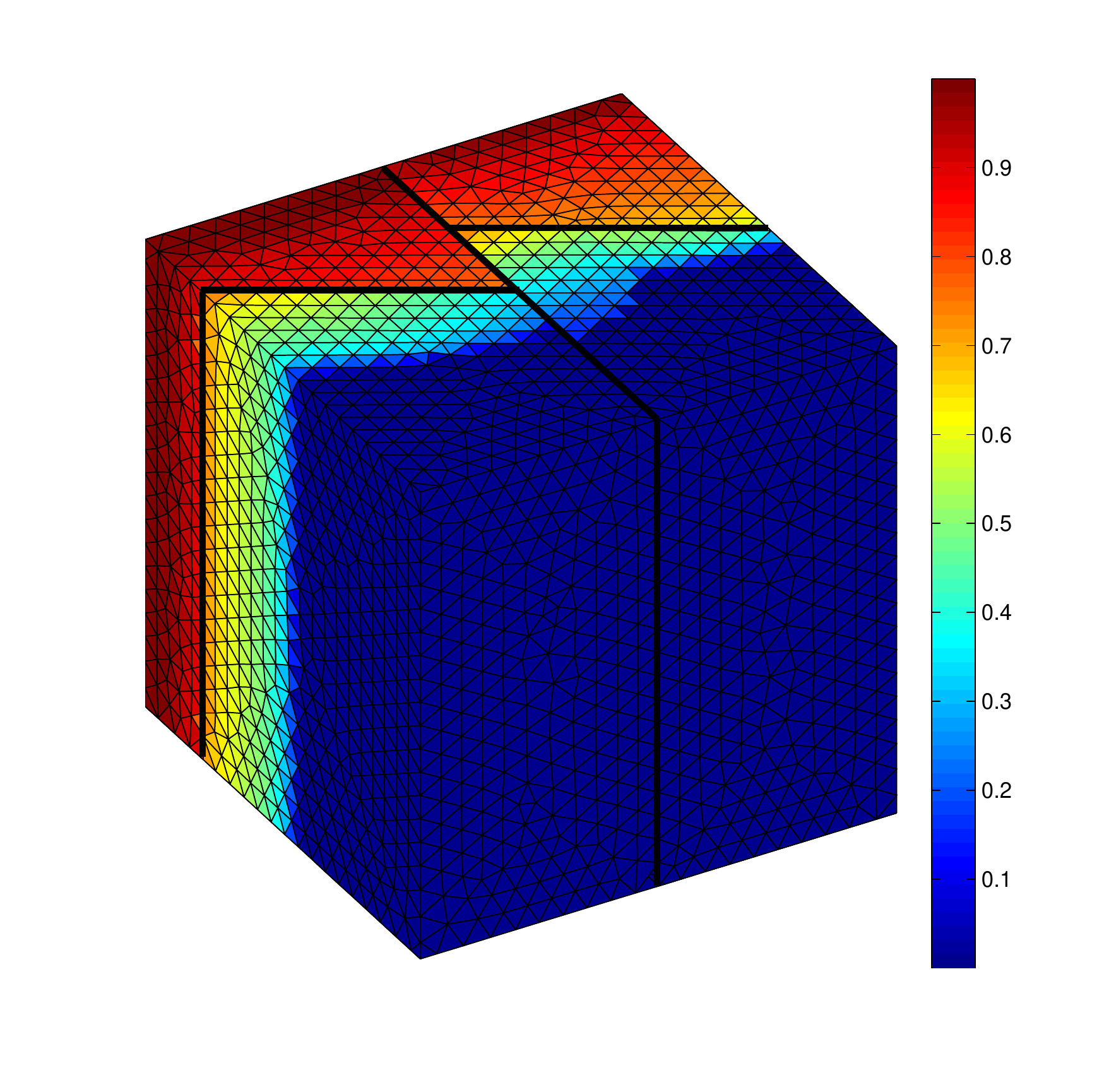}}
\caption{Test case~3 (multiple rocks): saturation $s(t)$ for  $t=4.6\cdot 10^4$ and $t=1.9\cdot 10^5$.}
\label{MRT-domain_saturation}
\end{figure}

We show the results in~Figure~\ref{MRT-domain_saturation} for two time steps.  One remark  that a very sharp and discontinuous change in saturation
at the rock type interfaces happens due to the higher contrast in the capillary pressure. The saturation of the wetting phase in the more permeable rocks (Rock~1 and 2)  is increasing rapidly in these two rocks and   a small quantity of the wetting phase penetrates into the Rocks~3 and 4 in which advection effects are neglected. We  remark in Figure~\ref{MRT-domain_saturation}  that the wetting phase  propagates in these subdomains   with a finite speed due   to  their  low  permeabilities and  due to the only diffusion effects that are present. In terms of computational effort, a good performance of the various inner algorithms is
observed.
\section{Conclusion}\label{section:conclusion}
We propose in this paper
a  splitting-based domain decomposition  methods to simulate two-phase flow 
in a porous medium composed of  different rock types.  The solution is resolved through a sequential  approach that consists 
of splitting the original problem into three (or two) problems, where  (global) pressure, saturation-advection and saturation-diffusion (or full saturation) problems are solved sequentially at each time step.  The resulting  schemes, that  differ by how  we treat the saturation problem (fully coupled or decoupled) provides us  with flexible and efficient ways   to treat the discontinuity of the saturation between rock types,  as we can adapt  the time  scales for the advection and diffusion effects, as well as we can adapt  the time scales for advection with respect to the rock type. 
Numerical experiments including those with several rock types and fractures  showcase the computational efficiency of the methods  and highlight their  flexibility  to handle  complex transmission  conditions between different rock types.  Work underway addresses the stability of the derived schemes and  provides adaptive multirate and local time steps    based  on  a  posteriori error estimates.
\section*{Acknowledgments}
This research was partially funded by the Hydrinv Inria Euro Med 3+3: HYDRINV project.  
It has also received funding from EPIC project (within LIRIMA: http://lirima.inria.fr) and the 
Tunisian Ministry of Higher Education and Scientific Research. The author thanks 
J\'er\^ome Jaffr\'{e} and Jean E. Roberts for helpful
discussions. 
\appendix\normalsize
 \section{Appendix: Application to a reduced fracture model between two rock types}\label{section:appendix}
\begin{figure}[h]
\begin{center}
{\includegraphics[scale=0.5]{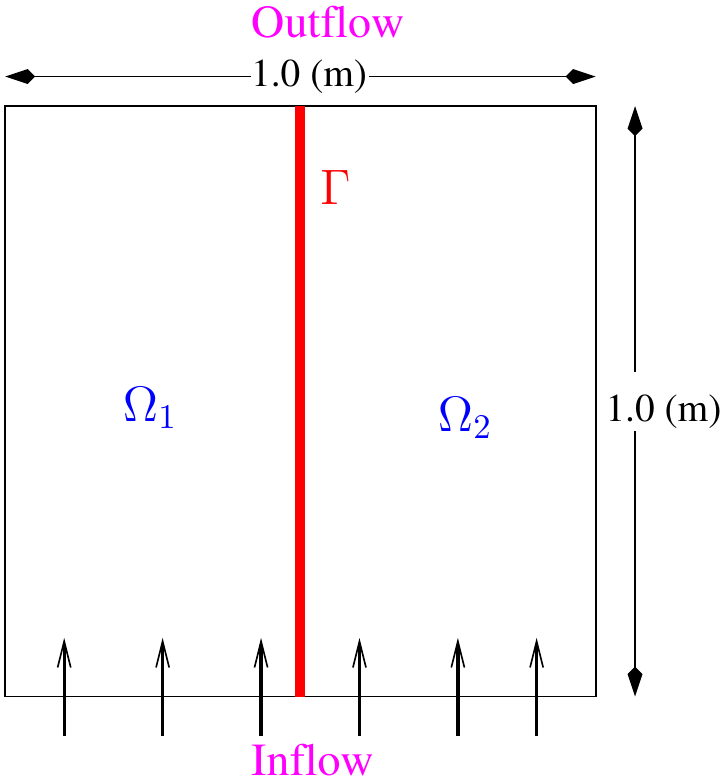}}
{\includegraphics[scale=0.37]{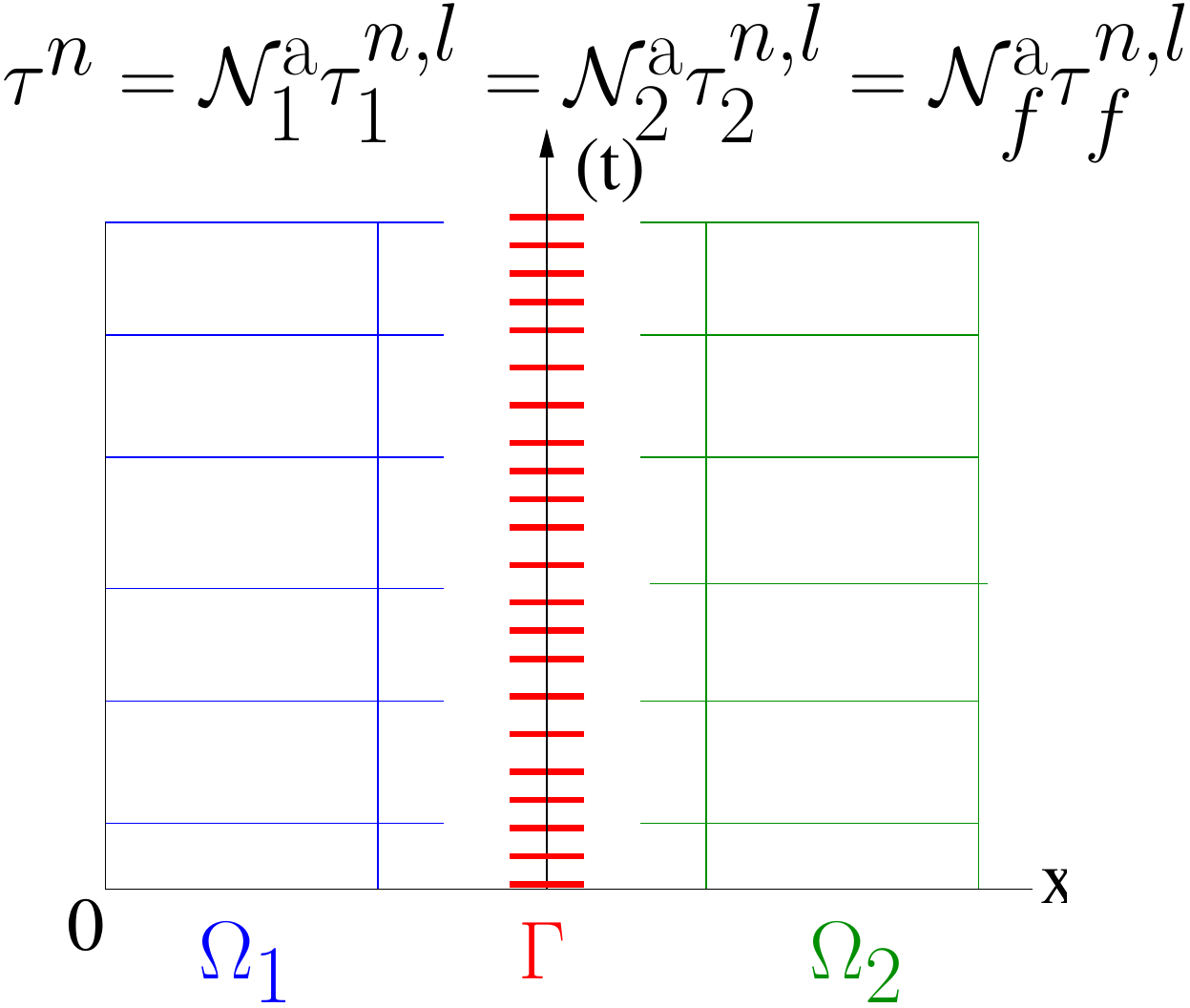}}
\caption{Reduced fracture model: geometry of the test case where the fracture is considered as an interface (left),  and  multirate time steps with nonconforming   grids for advection  in the rock matrix and in the fracture (right). }
\label{geo_fracture}
\end{center}
\end{figure}

\begin{figure}[h]\centering{
\includegraphics[scale=0.35]{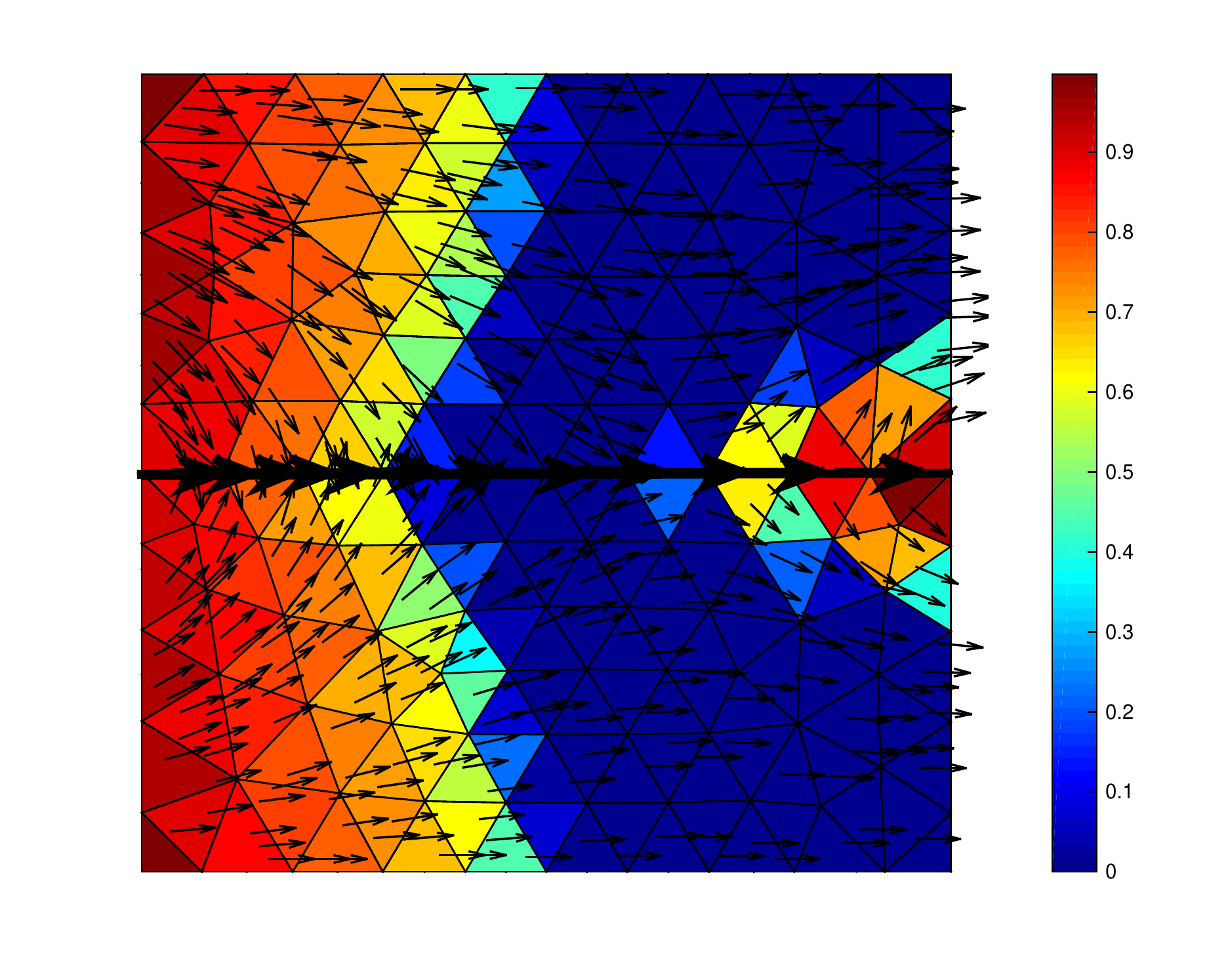}
\includegraphics[scale=0.25]{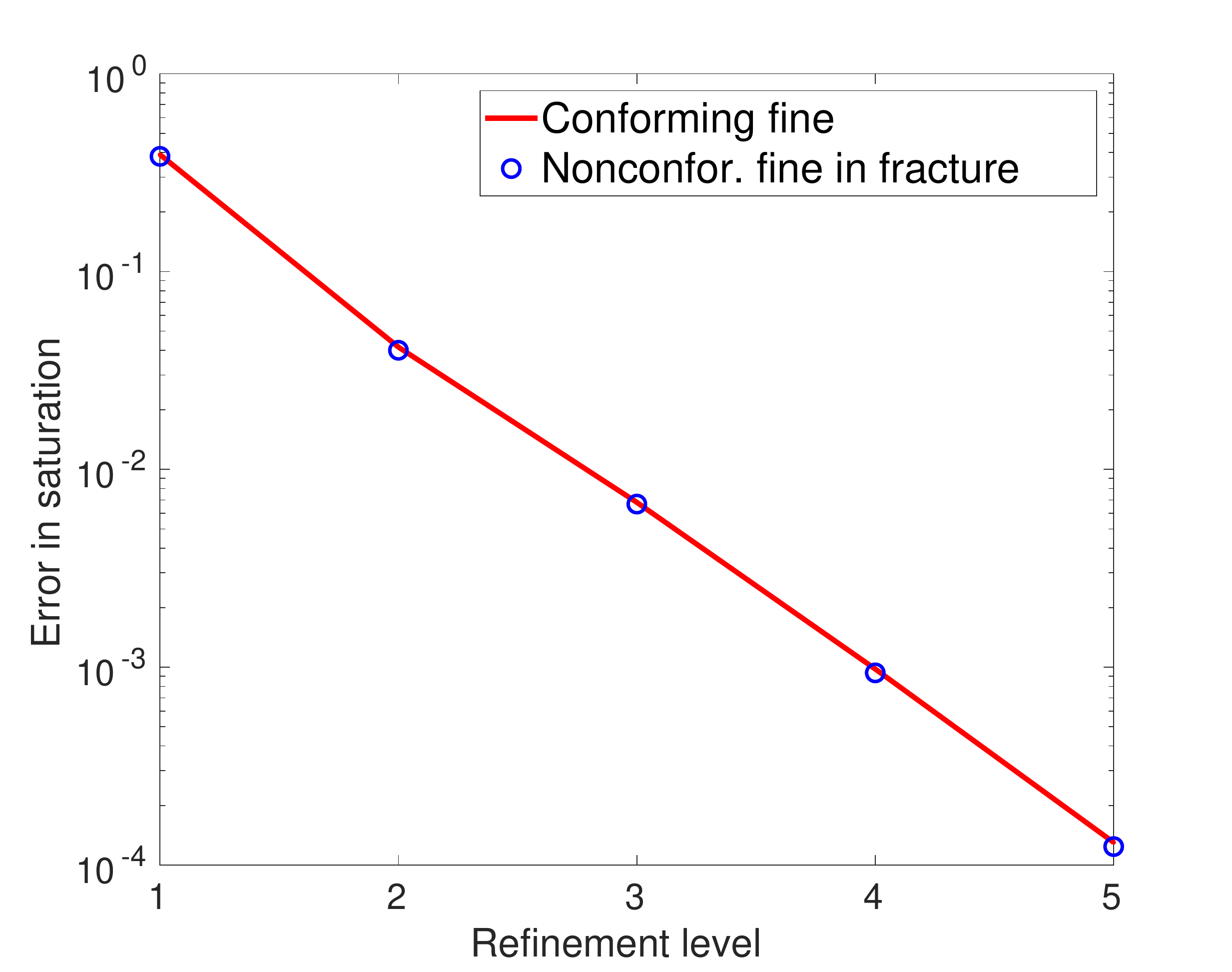}
}
\caption{Reduced fracture model: saturation $s(t)$ and  total velocity $\vu(t)$ for    $t=1\cdot 10^3$ (left) and errors in
 saturation between the reference  and the numerical solutions versus the refinement level (right).}
\label{MRT-domain_saturation_fracture}
\end{figure}

A further important feature of the developed DD approaches is the ability to integrate   models for   reduced fractures  for two-phase flow 
 in a natural way. Precisely, we apply Algorithm~\ref{algo_1} to the  discrete fracture model presented in~\cite{Ahmed2016} in which a fracture is treated as an interface of dimension 1 in a 2-dimensional simulation, with fluid exchange between the 1-dimensional fracture flow and the 2-dimensional flow in the surrounding rock matrix (see Figure~\ref{geo_fracture}). Next, we  give a  short overview of the model and the used discrete scheme.
\\\\
\textbf{The model} as presented in~\cite{Ahmed2016} is given by two-phase model problem~\eqref{mathematical_problem_DD} in each space--time domain $\mathrm{\Omega}_{i}\times(0,T)$ together with the following two phase flow
in  the fracture interface:
\bse\label{mathematical_fracture_DD}\begin{alignat}{2}
\label{Chap2_System_Conserv_wett_f}&\mathrm{\Phi}_{f} \dfrac{\p {\hat{s}}}{\p t} + \nabla_{{\tau}}\cdot\left(\vf_{f}(\hat{s})+\vr_{f}(\hat{s})\right) =\jump{\vu_{w}\cdot\vecn},\quad&\textnormal{ on }\mathrm{\Gamma}\times(0,T),\\
\label{Chap2_System_Conserv_add_f}&\vf_{f}(\hat{s})=f_{f}(\hat{s})\hat{\vu}_{f}+  f_{gf}(\hat{s})\vu_{g\tau},\quad&\textnormal{ on }\mathrm{\Gamma}\times(0,T),\\
\label{Chap2_System_Conserv_dif_f}&\vr_{f}(\hat{s})=- \vK_{f\tau}\nabla_{{\tau}} \alpha_{f}(\hat{s}),\quad&\textnormal{ on }\mathrm{\Gamma}\times(0,T),\\
\label{Chap2_System_Conserv_total_f}&\nabla_{{\tau}}\cdot\hat{\vu}_{f} =\jump{\vu\cdot\vecn},\quad&\textnormal{ on }\mathrm{\Gamma}\times(0,T),\\
\label{Chap2_System_Darcy_total_f}&\hat{\vu}_{f}=-\mathbf{M}(\hat{s})(\nabla_{{\tau}} \hat{p}-\rho_{f}(\hat{s})\vu_{g\tau}),\quad&\textnormal{ on }\mathrm{\Gamma}\times(0,T),
\end{alignat}\ese
where $\nabla_{{\tau}}$ denotes tangential component of the gradient operator, $\vu_{g\tau}$ is the tangential component of $\vu_{g}$
on $\mathrm{\Gamma}$, and where the functions $\mathbf{M}_{f}(\hat{s})$, $f_{gf}(\hat{s})$  and $\alpha_{f}(\hat{s})$
are defined using $\vK_{f\tau}$ the tangential component of the permeability on the fracture $\vK_{f}$ (see~\cite{Ahmed2016,Brenner2017,list2018upscaling} for more details). 
 We impose homogeneous Neumann boundary condition on $\partial\Gamma\times (0,T)$ and we assume  that  the  initial   saturation  is
known. For  model closure, we introduce
the following coupling conditions, for $i\in\{1,2\}$, 
\bse\label{Chap1_Matching_pressuref}\begin{alignat}{2} 
\label{Chap1_Matching_pressure_f}&p_{i}-\beta_{i}(s_{i})=\hat{p}-\beta_{f}(\hat{s}),\quad&\textnormal{ on }\mathrm{\Gamma}\times(0,T),\\
\label{Chap1_Matching_capilarity_f}&\pi_{i}(s_{i})=\pi_{f}(\hat{s}),\quad&\textnormal{on }\mathrm{\Gamma}\times(0,T).
\end{alignat}\ese
The above  systems are also coupled through the source terms appearing in the conservation equations in the 
fracture~(see Eqs~\eqref{Chap2_System_Conserv_wett_f} and~\eqref{Chap2_System_Conserv_total_f}), representing the difference between the fluid entering the fracture from one subdomain and that exiting through the other subdomain.   
\\\\
\textbf{The scheme} is given by  extending  Scheme~\ref{flowchart} to  the above setting, leading to solve sequentially   reduced  pressure,   saturation-advection  and saturation-diffusion interface problems posed only  on the fracture. With the use of the definitions~\eqref{op:dir_to_neuma}--\eqref{nonlinear_bilinear_forms_discrete}, we  promptly arrive to rewrite the pressure problem~\eqref{Chap2_System_Conserv_total_f}-\eqref{Chap2_System_Darcy_total_f} as an interface problem, in which simply we replace  the  pressure problem~\eqref{Chap1_interface_pressure}  by the following:  at each coarse time step $n$,  with known saturations $\hat{s}^{n}_{h}$, we solve for $\hat{p}^{n}_{h} \in \mathrm{\Lambda}_{h}$ such that
\begin{equation*}
\nabla_{{\tau}}\cdot\left(-\mathbf{M}_{f}(\hat{s}^{n}_{h})(\nabla_{{\tau}} \hat{p}^{n}_{h}-\rho_{f}(\hat{s}^{n}_{h})\vu_{g\tau})\right)+\mathcal{S}^{\textn{DtN}}_{\Gamma,n} (\hat{p}^{n}_{h})= g^{n}_{h}.
\end{equation*}     
The   operator associated to this reduced problem is  symmetric  positive definite and  can be solved  using   the  CG method~(cf.~\cite{ahmed2019robust} for details).  Now, we turn to the saturation problem. The  flexibility  in  the
time scales in the subdomains  and  in  the  fracture is a crucial asset in our numerical method,  and it allows to significantly improve the accuracy of the 
scheme when highly permeable fractures are present between different rock types.  The  same notations to obtain~\eqref{Chap2_advection_flux_definit2} allows to take into account the advection effects on  the fracture from the subdomains;  we simply replace~\eqref{Chap2_advection_flux_definit2}  by  solving for $\theta^{n,l}_{f,h}:=\hat{s}^{n,l}_{h}\in\Lambda_{h}$, such that
\begin{alignat*}{2} 
\label{discrete_advec_frac}&\int_{\sigma} \mathrm{\Phi}_{\sigma} \dfrac{\theta^{n,l+1}_{f,h}-\theta^{n,l}_{f,h}}{\tau^{n,l+1}_{f}} \,\text{d}\sigma + 
\displaystyle\sum_{e\in \mathcal{E}_{\sigma} } |e| \varphi^{n,l}_{\sigma,e} = \displaystyle  \int_{I^{n,l}_{f}}|\sigma|\{\mathcal{P}_{f1}(\varphi_{1}(s^{n,a}_{K},\theta^{n,a}_{h,1}))+\mathcal{P}_{f2}(\varphi_{2}(s^{n,a}_{L},\theta^{n,a}_{h,2}))\},\, \forall \sigma=K|L \in \mathcal{E}^{\mathrm{\Gamma}}_{h},\\
&\displaystyle  \int_{I^{n,l}_{f}}\{\pi_{f}(\theta^{n,l}_{f,h})-\mathcal{P}_{fi}(\pi_{i}(\theta^{n,a}_{h,i}))\}=0,\qquad \forall i\in\{1,2\},
\end{alignat*}
for all $l\in\{0,\cdots,\mathcal{N}^{\textnormal{a}}_{f}-1\}$,  where  $\mathcal{P}_{fi}$, $i\in\{1,2\}$, are projection-in-time type operators  
given by~\eqref{projection_oper}. Therein,   $|e| \varphi^{n,l}_{\sigma,e} $ is an approximation of the advection flux through the edge $e$,  $\int_{e} \vf^{n}_{f}(\hat{s})\cdot\vecn_{e}$. Similarly to $\varphi^{n,l}_{K,\sigma}$ in~\eqref{Chap2_advection_flux_definition},  $\varphi^{n,l}_{\sigma,e}$ is a function of the two values of the saturation adjacent to~$e=\sigma^{-}|\sigma^{+}$ 
and   we calculate it  using the introduced  Godunov scheme with $i=f$. Now, it remains to extend the diffusion step.  As for the reduced pressure problem, we replace~\eqref{Chap1_interface_diffusion} by solving 
$\theta^{n+1}_{f,h} \in \mathrm{\Lambda}_{h}$ such that
\begin{equation*}
\label{Chap1_interface_diffusion_f}
 \mathrm{\Phi}_{f}\dfrac{\theta^{n+1}_{f,h}-\theta^{n,\mathcal{N}^{\textnormal{a}}_{f}}_{f,h}}{\tau^{n+1}_{f}} +\nabla_{\tau}\cdot\left(-\vK_{f\tau}\nabla_{\tau}\alpha_{f}(\theta^{n+1}_{f,h})\right)+\mathcal{Z}^{\textn{DtN}}_{\Gamma,n+1}(\pi_{f}( \theta^{n+1}_{f,h}))=0.
\end{equation*}
  This problem can be solved iteratively using fixed point  iterations or the Newton method. The result of this system is then  the saturation at the coarse time step $n+1$, i.e., $\hat{s}^{n+1}_{h}=\theta^{n+1}_{f,h}$.
\begin{rem}[The interface preconditioners]
 To improve the convergence of the our DD method, the inverse of second order operator on the fracture $\nabla_{{\tau}}\cdot\left(-M_{f}(\hat{s}^{n}_{h})\nabla_{{\tau}} \cdot\right)$ 
 is used as preconditioner for the interface pressure problem~\cite{budivsa2019block,Ahmed2016,hoang2016space}. Similarly, a linearized version of the interface operator 
 $\nabla_{\tau}\cdot\left(\vK_{f\tau}\nabla_{\tau}(\alpha_{f}(\cdot)\right)$ is used as preconditioner for GMRes when Newton method is used to solve 
 the interface diffusion problem (cf.~\cite{Ahmed2016,ahmed2019robust} for more details).
\end{rem}
\textbf{The test case} we consider takes  a unit square domain $\Omega$,       divided  into
two  equally  sized  subdomains  by  a  fracture $\mathrm{\Gamma}$ located  in the middle of the domain, perpendicular 
to the injection and production faces of $\mathrm{\Omega}$. 
The  permeability  of  the  matrix  is
 given by $\vK= 1$,  very  low  compared  to  the  permeability  of  the
fracture $\vK_{f}=1\cdot 10^{2}$.  The porosity is equal to $\mathrm{\Phi}=0.1$ in the matrix and to 
$\mathrm{\Phi}=0.3$ in the fracture. 
We consider a triangular mesh with 800 grids. In time, we  fix $T=5\cdot 10^3$, and use uniform coarse time steps  in the fracture and subdomains, i.e., $\tau^{n}=T/100$. For the advection, we choose a finer and nonconforming time steps between the subdomains and the fracture; $\tau_{i}^{n,l}=\tau^{n}/5$, $i\in\{1,2\}$, and $\tau_{f}^{n,l}=\tau^{n}/20$ as depicted in Figure~\ref{geo_fracture} (right). In 
Figure~\ref{MRT-domain_saturation_fracture} (left),  the saturation and 
the total velocity at $t^{n}=T/5$ is shown.
In this experiment, the flow is driven not only by the difference of the permeabilities and the capillary pressure fields, but also by the presence of the fracture. The wetting phase  
moves  immediately through  the fracture and the saturation front in the surrounding matrix  snakes along the matrix-fracture interfaces. The wetting phase accumulates into the fracture near the production boundary 
until it eventually spreads out into the surrounding matrix. To verify the accuracy in time, we consider 
two initial time grids refined then  5 times by a factor of 2:
\begin{itemize}
 \item Time grid 1 (Conforming fine): $\tau^{n,l}_{i}=\tau^{n}/5$, $i\in\{1,2\}$,
 \item Time grid 2 (Nonconforming, fine in fracture): $\tau^{n,l}_{i}=\tau^{n}/5$, $i\in\{1,2\},$ and $\tau^{n,l}_{f}=\tau^{n}/20$.
\end{itemize}
A reference solution is calculated with a very fine and single-rate time step and on fixed mesh with Algorithm~\ref{algo_2}. In Figure~\ref{MRT-domain_saturation_fracture} (right)
the error in the $L^2(0,T;L^2(\Gamma))$-norm of the saturation in the fracture
versus the refinement level is depicted.
As expected,  first order convergence is ensured in the nonconforming case and the errors
obtained in the nonconforming case with a finer  time step in the fracture
are  nearly  the  same  as  in  the  finer  conforming case.   Thus,  the  use  of nonconforming and multirate grids preserves the accuracy in time. Note that excellent results are obtained with a ratio 10 of the fine time step to the coarse  time step. 

\enlargethispage{0.2cm}
\bibliographystyle{siamplain}


\begin{thebibliography}{10}

\bibitem{MR2051062}
{\sc Adimurthi, J.~Jaffr{\'e}, and G.~D. {Veerappa Gowda}}, {\em Godunov-type
  methods for conservation laws with a flux function discontinuous in space},
  SIAM J. Numer. Anal., 42 (2004), pp.~179--208,
  \href{http://dx.doi.org/10.1137/S003614290139562X}{doi:\nolinkurl{10.1137/S003614290139562X}},
  \url{https://doi.org/10.1137/S003614290139562X}.

\bibitem{ahmed:hal-01540956}
{\sc E.~Ahmed, S.~{Ali Hassan}, C.~Japhet, M.~Kern, and M.~Vohral\'{\i}k}, {\em
  A posteriori error estimates and stopping criteria for space-time domain
  decomposition for two-phase flow between different rock types}.
\newblock working paper or preprint, June 2017,
  \url{https://hal.inria.fr/hal-01540956}.

\bibitem{ahmed2018multiscale}
{\sc E.~Ahmed, A.~Fumagalli, and A.~Budi{\v s}a}, {\em A multiscale flux basis
  for mortar mixed discretizations of reduced darcy-forchheimer fracture
  models}, Computer Methods in Applied Mechanics and Engineering,  (2019),
  \href{http://dx.doi.org/10.1016/j.cma.2019.05.034}{doi:\nolinkurl{10.1016/j.cma.2019.05.034}}.

\bibitem{ahmed2019robust}
{\sc E.~Ahmed, A.~Fumagalli, A.~Budi{\v s}a, E.~Keilegavlen, J.~M. Nordbotten,
  and A.~R. Forin}, {\em Robust linear domain decomposition schemes for reduced
  non-linear fracture flow models}, arXiv preprint arXiv:1906.05831,  (2019).

\bibitem{Ahmed2016}
{\sc E.~Ahmed, J.~Jaffr{\'e}, and J.~E. Roberts}, {\em A reduced fracture model
  for two-phase flow with different rock types}, Math. Comput. Simulation, 137
  (2017), pp.~49--70,
  \href{http://dx.doi.org/10.1016/j.matcom.2016.10.005}{doi:\nolinkurl{10.1016/j.matcom.2016.10.005}},
  \url{https://doi.org/10.1016/j.matcom.2016.10.005}.

\bibitem{ahmed2018global}
{\sc E.~Ahmed, C.~Japhet, and M.~Kern}, {\em Global-in-time domain
  decomposition for a nonlinear diffusion problem}, in HAL Preprint 02263280,
  to appear in Domain Decomposition Methods in Science and Engineering XXV,
  Springer International Publishing, 2019.

\bibitem{ahmed:hal-02275690}
{\sc E.~Ahmed, C.~Japhet, and M.~Kern}, {\em Space-time domain decomposition
  for two-phase flow between different rock types}.
\newblock working paper or preprint, Aug. 2019,
  \url{https://hal.inria.fr/hal-02275690}.

\bibitem{MR3983155}
{\sc E.~Ahmed, J.~M. Nordbotten, and F.~A. Radu}, {\em Adaptive asynchronous
  time-stepping, stopping criteria, and a posteriori error estimates for
  fixed-stress iterative schemes for coupled poromechanics problems}, J.
  Comput. Appl. Math., 364 (2020),
  \href{http://dx.doi.org/10.1016/j.cam.2019.06.028}{doi:\nolinkurl{10.1016/j.cam.2019.06.028}},
  \url{https://doi.org/10.1016/j.cam.2019.06.028}.

\bibitem{alboin2000domain}
{\sc C.~Alboin, J.~Jaffr{\'e}, J.~E. Roberts, and C.~Serres}, {\em Modeling
  fractures as interfaces for flow and transport in porous media}, in Fluid
  flow and transport in porous media: mathematical and numerical treatment
  ({S}outh {H}adley, {MA}, 2001), vol.~295 of Contemp. Math., Amer. Math. Soc.,
  Providence, RI, 2002, pp.~13--24,
  \href{http://dx.doi.org/10.1090/conm/295/04999}{doi:\nolinkurl{10.1090/conm/295/04999}}.

\bibitem{alboinDD}
{\sc C.~Alboin, J.~Jaffr{\'e}, J.~E. Roberts, X.~Wang, and C.~Serres}, {\em
  Domain decomposition for some transmission problems in flow in porous media},
  in Numerical Treatment of Multiphase Flows in Porous Media, Z.~Chen, R.~E.
  Ewing, and Z.-C. Shi, eds., Berlin, Heidelberg, 2000, Springer Berlin
  Heidelberg, pp.~22--34.

\bibitem{MR3564686}
{\sc T.~Almani, K.~Kumar, A.~Dogru, G.~Singh, and M.~F. Wheeler}, {\em
  Convergence analysis of multirate fixed-stress split iterative schemes for
  coupling flow with geomechanics}, Comput. Methods Appl. Mech. Engrg., 311
  (2016), pp.~180--207,
  \href{http://dx.doi.org/10.1016/j.cma.2016.07.036}{doi:\nolinkurl{10.1016/j.cma.2016.07.036}},
  \url{https://doi.org/10.1016/j.cma.2016.07.036}.

\bibitem{MR3208750}
{\sc B.~Andreianov and C.~Canc{\`e}s}, {\em A phase-by-phase upstream scheme
  that converges to the vanishing capillarity solution for countercurrent
  two-phase flow in two-rock media}, Comput. Geosci., 18 (2014), pp.~211--226,
  \href{http://dx.doi.org/10.1007/s10596-014-9403-5}{doi:\nolinkurl{10.1007/s10596-014-9403-5}},
  \url{https://doi.org/10.1007/s10596-014-9403-5}.

\bibitem{MR3142429}
{\sc T.~Arbogast, M.~Juntunen, J.~Pool, and M.~F. Wheeler}, {\em A
  discontinuous {G}alerkin method for two-phase flow in a porous medium
  enforcing {$H({\textnormal{div}})$} velocity and continuous capillary
  pressure}, Comput. Geosci., 17 (2013), pp.~1055--1078,
  \href{http://dx.doi.org/10.1007/s10596-013-9374-y}{doi:\nolinkurl{10.1007/s10596-013-9374-y}},
  \url{https://doi.org/10.1007/s10596-013-9374-y}.

\bibitem{aziz1979petroleum}
{\sc K.~Aziz and A.~Settari}, {\em Petroleum Reservoir Simulation}, vol.~476,
  Applied Science Publishers London, 1979.

\bibitem{Brenner2017}
{\sc K.~Brenner, M.~Groza, L.~Jeannin, R.~Masson, and J.~Pellerin}, {\em
  Immiscible two-phase {D}arcy flow model accounting for vanishing and
  discontinuous capillary pressures: application to the flow in fractured
  porous media}, Comput. Geosci., 21 (2017), pp.~1075--1094,
  \href{http://dx.doi.org/10.1007/s10596-017-9675-7}{doi:\nolinkurl{10.1007/s10596-017-9675-7}},
  \url{https://doi.org/10.1007/s10596-017-9675-7}.

\bibitem{brun2019monolithic}
{\sc M.~K. Brun, E.~Ahmed, I.~Berre, J.~M. Nordbotten, and F.~A. Radu}, {\em
  Monolithic and splitting based solution schemes for fully coupled
  quasi-static thermo-poroelasticity with nonlinear convective transport},
  arXiv preprint arXiv:1902.05783,  (2019).

\bibitem{budivsa2019block}
{\sc A.~Budi{\v s}a and X.~Hu}, {\em Block preconditioners for
  mixed-dimensional discretization of flow in fractured porous media}, arXiv
  preprint arXiv:1905.13513,  (2019).

\bibitem{MR2465972}
{\sc C.~Canc{\`e}s}, {\em Nonlinear parabolic equations with spatial
  discontinuities}, NoDEA Nonlinear Differential Equations Appl., 15 (2008),
  pp.~427--456,
  \href{http://dx.doi.org/10.1007/s00030-008-6030-7}{doi:\nolinkurl{10.1007/s00030-008-6030-7}},
  \url{https://doi.org/10.1007/s00030-008-6030-7}.

\bibitem{MR2559741}
{\sc C.~Canc{\`e}s}, {\em Finite volume scheme for two-phase flows in
  heterogeneous porous media involving capillary pressure discontinuities},
  M2AN Math. Model. Numer. Anal., 43 (2009), pp.~973--1001,
  \href{http://dx.doi.org/10.1051/m2an/2009032}{doi:\nolinkurl{10.1051/m2an/2009032}},
  \url{https://doi.org/10.1051/m2an/2009032}.

\bibitem{MR1370250}
{\sc C.~Carlenzoli and A.~Quarteroni}, {\em Adaptive domain decomposition
  methods for advection-diffusion problems}, in Modeling, mesh generation, and
  adaptive numerical methods for partial differential equations ({M}inneapolis,
  {MN}, 1993), vol.~75 of IMA Vol. Math. Appl., Springer, New York, 1995,
  pp.~165--186,
  \href{http://dx.doi.org/10.1007/978-1-4612-4248-2_9}{doi:\nolinkurl{10.1007/978-1-4612-4248-2_9}},
  \url{https://doi.org/10.1007/978-1-4612-4248-2_9}.

\bibitem{chavent1976new}
{\sc G.~Chavent}, {\em A new formulation of diphasic incompressible flows in
  porous media},  (1976), pp.~258--270. Lecture Notes in Math., 503.

\bibitem{chavent1986mathematical}
{\sc G.~Chavent and J.~Jaffr{\'e}}, {\em Mathematical models and finite
  elements for reservoir simulation: single phase, multiphase and
  multicomponent flows through porous media}, vol.~17, North Holland, 1986.

\bibitem{chen2006computational}
{\sc Z.~Chen, G.~Huan, and Y.~Ma}, {\em Computational methods for multiphase
  flows in porous media}, vol.~2 of Computational Science \& Engineering,
  Society for Industrial and Applied Mathematics (SIAM), Philadelphia, PA,
  2006,
  \href{http://dx.doi.org/10.1137/1.9780898718942}{doi:\nolinkurl{10.1137/1.9780898718942}},
  \url{https://doi.org/10.1137/1.9780898718942}.

\bibitem{MR1297465}
{\sc L.~C. Cowsar, J.~Mandel, and M.~F. Wheeler}, {\em Balancing domain
  decomposition for mixed finite elements}, Math. Comp., 64 (1995),
  pp.~989--1015,
  \href{http://dx.doi.org/10.2307/2153480}{doi:\nolinkurl{10.2307/2153480}},
  \url{https://doi.org/10.2307/2153480}.

\bibitem{Douglas1983}
{\sc J.~Douglas, Jr.}, {\em Finite difference methods for two-phase
  incompressible flow in porous media}, SIAM J. Numer. Anal., 20 (1983),
  pp.~681--696,
  \href{http://dx.doi.org/10.1137/0720046}{doi:\nolinkurl{10.1137/0720046}},
  \url{https://doi.org/10.1137/0720046}.

\bibitem{MR2206441}
{\sc G.~Ench{\'e}ry, R.~Eymard, and A.~Michel}, {\em Numerical approximation of
  a two-phase flow problem in a porous medium with discontinuous capillary
  forces}, SIAM J. Numer. Anal., 43 (2006), pp.~2402--2422,
  \href{http://dx.doi.org/10.1137/040602936}{doi:\nolinkurl{10.1137/040602936}},
  \url{https://doi.org/10.1137/040602936}.

\bibitem{MR2989844}
{\sc A.~Ern, I.~Mozolevski, and L.~Schuh}, {\em Corrigendum to
  ``{D}iscontinuous {G}alerkin approximation of two-phase flows in
  heterogeneous porous media with discontinuous capillary pressures''
  [{C}omput. {M}ethods {A}ppl. {M}ech. {E}ngrg. 199 (2010) 1491--1501]
  [mr2630157]}, Comput. Methods Appl. Mech. Engrg., 245/246 (2012),
  pp.~348--349,
  \href{http://dx.doi.org/10.1016/j.cma.2012.05.011}{doi:\nolinkurl{10.1016/j.cma.2012.05.011}},
  \url{https://doi.org/10.1016/j.cma.2012.05.011}.

\bibitem{MR3392446}
{\sc R.~Eymard, T.~Gallou{\"e}t, and R.~Herbin}, {\em {$\mathcal{RT}_k$} mixed
  finite elements for some nonlinear problems}, Math. Comput. Simulation, 118
  (2015), pp.~186--197,
  \href{http://dx.doi.org/10.1016/j.matcom.2014.11.013}{doi:\nolinkurl{10.1016/j.matcom.2014.11.013}},
  \url{https://doi.org/10.1016/j.matcom.2014.11.013}.

\bibitem{FreyGeorge2008}
{\sc P.~J. Frey and P.-L. George}, {\em Mesh generation}, ISTE, London; John
  Wiley \& Sons, Inc., Hoboken, NJ, second~ed., 2008,
  \href{http://dx.doi.org/10.1002/9780470611166}{doi:\nolinkurl{10.1002/9780470611166}},
  \url{https://doi.org/10.1002/9780470611166}.
\newblock Application to finite elements.

\bibitem{MR2218966}
{\sc M.~J. Gander}, {\em Optimized {S}chwarz methods}, SIAM J. Numer. Anal., 44
  (2006), pp.~699--731,
  \href{http://dx.doi.org/10.1137/S0036142903425409}{doi:\nolinkurl{10.1137/S0036142903425409}},
  \url{https://doi.org/10.1137/S0036142903425409}.

\bibitem{MR2344706}
{\sc M.~J. Gander, L.~Halpern, and F.~Magoul{\`e}s}, {\em An optimized
  {S}chwarz method with two-sided {R}obin transmission conditions for the
  {H}elmholtz equation}, Internat. J. Numer. Methods Fluids, 55 (2007),
  pp.~163--175,
  \href{http://dx.doi.org/10.1002/fld.1433}{doi:\nolinkurl{10.1002/fld.1433}},
  \url{https://doi.org/10.1002/fld.1433}.

\bibitem{MR2235750}
{\sc M.~J. Gander, C.~Japhet, Y.~Maday, and F.~Nataf}, {\em A new cement to
  glue nonconforming grids with {R}obin interface conditions: the finite
  element case}, in Domain decomposition methods in science and engineering,
  vol.~40 of Lect. Notes Comput. Sci. Eng., Springer, Berlin, 2005,
  pp.~259--266,
  \href{http://dx.doi.org/10.1007/3-540-26825-1_24}{doi:\nolinkurl{10.1007/3-540-26825-1_24}},
  \url{https://doi.org/10.1007/3-540-26825-1_24}.

\bibitem{ganis2014global}
{\sc B.~Ganis, K.~Kumar, G.~Pencheva, M.~F. Wheeler, and I.~Yotov}, {\em A
  global {J}acobian method for mortar discretizations of a fully implicit
  two-phase flow model}, Multiscale Model. Simul., 12 (2014), pp.~1401--1423,
  \href{http://dx.doi.org/10.1137/140952922}{doi:\nolinkurl{10.1137/140952922}},
  \url{https://doi.org/10.1137/140952922}.

\bibitem{MR3771343}
{\sc Z.~Ge and M.~Ma}, {\em Multirate iterative scheme based on multiphysics
  discontinuous {G}alerkin method for a poroelasticity model}, Appl. Numer.
  Math., 128 (2018), pp.~125--138,
  \href{http://dx.doi.org/10.1016/j.apnum.2018.02.003}{doi:\nolinkurl{10.1016/j.apnum.2018.02.003}},
  \url{https://doi.org/10.1016/j.apnum.2018.02.003}.

\bibitem{MR3144798}
{\sc T.-T.-P. Hoang, J.~Jaffr{\'e}, C.~Japhet, M.~Kern, and J.~E. Roberts},
  {\em Space-time domain decomposition methods for diffusion problems in mixed
  formulations}, SIAM J. Numer. Anal., 51 (2013), pp.~3532--3559,
  \href{http://dx.doi.org/10.1137/130914401}{doi:\nolinkurl{10.1137/130914401}},
  \url{https://doi.org/10.1137/130914401}.

\bibitem{hoang2016space}
{\sc T.-T.-P. Hoang, C.~Japhet, M.~Kern, and J.~E. Roberts}, {\em Space-time
  domain decomposition for reduced fracture models in mixed formulation}, SIAM
  J. Numer. Anal., 54 (2016), pp.~288--316,
  \href{http://dx.doi.org/10.1137/15M1009651}{doi:\nolinkurl{10.1137/15M1009651}},
  \url{https://doi.org/10.1137/15M1009651}.

\bibitem{HOANG2017366}
{\sc T.-T.-P. Hoang, C.~Japhet, M.~Kern, and J.~E. Roberts}, {\em Space-time
  domain decomposition for advection-diffusion problems in mixed formulations},
  Math. Comput. Simulation, 137 (2017), pp.~366--389,
  \href{http://dx.doi.org/10.1016/j.matcom.2016.11.002}{doi:\nolinkurl{10.1016/j.matcom.2016.11.002}},
  \url{https://doi.org/10.1016/j.matcom.2016.11.002}.

\bibitem{hoteit2008efficient}
{\sc H.~Hoteit and A.~Firoozabadi}, {\em An efficient numerical model for
  incompressible two-phase flow in fractured media}, Advances in Water
  Resources, 31 (2008), pp.~891--905,
  \href{http://dx.doi.org/10.1016/j.advwatres.2008.02.004}{doi:\nolinkurl{10.1016/j.advwatres.2008.02.004}},
  \url{http://www.sciencedirect.com/science/article/pii/S0309170808000353}.

\bibitem{hoteit2008numerical}
{\sc H.~Hoteit and A.~Firoozabadi}, {\em Numerical modeling of two-phase flow
  in heterogeneous permeable media with different capillarity pressures},
  Advances in Water Resources, 31 (2008), pp.~56--73,
  \href{http://dx.doi.org/10.1016/j.advwatres.2007.06.006}{doi:\nolinkurl{10.1016/j.advwatres.2007.06.006}},
  \url{http://www.sciencedirect.com/science/article/pii/S030917080700108X}.

\bibitem{kelley1995}
{\sc C.~T. Kelley}, {\em Iterative methods for linear and nonlinear equations},
  vol.~16 of Frontiers in Applied Mathematics, Society for Industrial and
  Applied Mathematics (SIAM), Philadelphia, PA, 1995,
  \href{http://dx.doi.org/10.1137/1.9781611970944}{doi:\nolinkurl{10.1137/1.9781611970944}},
  \url{https://doi.org/10.1137/1.9781611970944}.
\newblock With separately available software.

\bibitem{MR2837398}
{\sc J.~Kou and S.~Sun}, {\em On iterative {IMPES} formulation for two-phase
  flow with capillarity in heterogeneous porous media}, Int. J. Numer. Anal.
  Model. Ser. B, 1 (2010), pp.~20--40.

\bibitem{2017arXiv170901644L}
{\sc S.~Lee and M.~F. Wheeler}, {\em Enriched {G}alerkin methods for two-phase
  flow in porous media with capillary pressure}, J. Comput. Phys., 367 (2018),
  pp.~65--86,
  \href{http://dx.doi.org/10.1016/j.jcp.2018.03.031}{doi:\nolinkurl{10.1016/j.jcp.2018.03.031}},
  \url{https://doi.org/10.1016/j.jcp.2018.03.031}.

\bibitem{MRST-Lie2014}
{\sc K.-A. Lie}, {\em An Introduction to Reservoir Simulation Using MATLAB:
  User guide for the Matlab Reservoir Simulation Toolbox (MRST)}, SINTEF ICT,
  Norway, 2014.

\bibitem{lie_2019}
{\sc K.-A. Lie}, {\em An Introduction to Reservoir Simulation Using MATLAB/GNU
  Octave: User Guide for the MATLAB Reservoir Simulation Toolbox (MRST)},
  Cambridge University Press, 2019,
  \href{http://dx.doi.org/10.1017/9781108591416}{doi:\nolinkurl{10.1017/9781108591416}}.

\bibitem{list2018upscaling}
{\sc F.~List, K.~Kumar, I.~S. Pop, and F.~A. Radu}, {\em Upscaling of
  unsaturated flow in fractured porous media}, arXiv preprint arXiv:1807.05993,
   (2018).

\bibitem{MR3489128}
{\sc F.~List and F.~A. Radu}, {\em A study on iterative methods for solving
  {R}ichards' equation}, Comput. Geosci., 20 (2016), pp.~341--353,
  \href{http://dx.doi.org/10.1007/s10596-016-9566-3}{doi:\nolinkurl{10.1007/s10596-016-9566-3}},
  \url{https://doi.org/10.1007/s10596-016-9566-3}.

\bibitem{MR1208381}
{\sc J.~Mandel}, {\em Balancing domain decomposition}, Comm. Numer. Methods
  Engrg., 9 (1993), pp.~233--241,
  \href{http://dx.doi.org/10.1002/cnm.1640090307}{doi:\nolinkurl{10.1002/cnm.1640090307}},
  \url{https://doi.org/10.1002/cnm.1640090307}.

\bibitem{Mishra2017479}
{\sc S.~Mishra}, {\em Numerical methods for conservation laws with
  discontinuous coefficients}, in Handbook of numerical methods for hyperbolic
  problems, vol.~18 of Handb. Numer. Anal., Elsevier/North-Holland, Amsterdam,
  2017, pp.~479--506.

\bibitem{MR3926838}
{\sc K.~Mitra and I.~S. Pop}, {\em A modified {L}-scheme to solve nonlinear
  diffusion problems}, Comput. Math. Appl., 77 (2019), pp.~1722--1738,
  \href{http://dx.doi.org/10.1016/j.camwa.2018.09.042}{doi:\nolinkurl{10.1016/j.camwa.2018.09.042}},
  \url{https://doi.org/10.1016/j.camwa.2018.09.042}.

\bibitem{MR2997427}
{\sc I.~Mozolevski and L.~Schuh}, {\em Numerical simulation of two-phase
  immiscible incompressible flows in heterogeneous porous media with capillary
  barriers}, J. Comput. Appl. Math., 242 (2013), pp.~12--27,
  \href{http://dx.doi.org/10.1016/j.cam.2012.09.045}{doi:\nolinkurl{10.1016/j.cam.2012.09.045}},
  \url{https://doi.org/10.1016/j.cam.2012.09.045}.

\bibitem{peszynska2002mortar}
{\sc M.~Peszy{\'n}ska, M.~F. Wheeler, and I.~Yotov}, {\em Mortar upscaling for
  multiphase flow in porous media}, Comput. Geosci., 6 (2002), pp.~73--100,
  \href{http://dx.doi.org/10.1023/A:1016529113809}{doi:\nolinkurl{10.1023/A:1016529113809}},
  \url{https://doi.org/10.1023/A:1016529113809}.

\bibitem{MR2079503}
{\sc I.~S. Pop, F.~Radu, and P.~Knabner}, {\em Mixed finite elements for the
  {R}ichards' equation: linearization procedure}, J. Comput. Appl. Math., 168
  (2004), pp.~365--373,
  \href{http://dx.doi.org/10.1016/j.cam.2003.04.008}{doi:\nolinkurl{10.1016/j.cam.2003.04.008}},
  \url{https://doi.org/10.1016/j.cam.2003.04.008}.

\bibitem{MR3800042}
{\sc F.~A. Radu, K.~Kumar, J.~M. Nordbotten, and I.~S. Pop}, {\em A robust,
  mass conservative scheme for two-phase flow in porous media including
  {H}{\"o}lder continuous nonlinearities}, IMA J. Numer. Anal., 38 (2018),
  pp.~884--920,
  \href{http://dx.doi.org/10.1093/imanum/drx032}{doi:\nolinkurl{10.1093/imanum/drx032}},
  \url{https://doi.org/10.1093/imanum/drx032}.

\bibitem{reichenberger2006mixed}
{\sc V.~Reichenberger, H.~Jakobs, P.~Bastian, and R.~Helmig}, {\em A
  mixed-dimensional finite volume method for two-phase flow in fractured porous
  media}, Advances in Water Resources, 29 (2006), pp.~1020--1036,
  \href{http://dx.doi.org/10.1016/j.advwatres.2005.09.001}{doi:\nolinkurl{10.1016/j.advwatres.2005.09.001}},
  \url{http://www.sciencedirect.com/science/article/pii/S0309170805002150}.

\bibitem{robertsJean}
{\sc J.~E. Roberts and J.-M. Thomas}, {\em Mixed and hybrid methods}, in
  Handbook of numerical analysis, {V}ol. {II}, Handb. Numer. Anal., II,
  North-Holland, Amsterdam, 1991, pp.~523--639.

\bibitem{MR3351781}
{\sc I.~Rybak, J.~Magiera, R.~Helmig, and C.~Rohde}, {\em Multirate time
  integration for coupled saturated/unsaturated porous medium and free flow
  systems}, Comput. Geosci., 19 (2015), pp.~299--309,
  \href{http://dx.doi.org/10.1007/s10596-015-9469-8}{doi:\nolinkurl{10.1007/s10596-015-9469-8}},
  \url{https://doi.org/10.1007/s10596-015-9469-8}.

\bibitem{MR3771899}
{\sc D.~Seus, K.~Mitra, I.~S. Pop, F.~A. Radu, and C.~Rohde}, {\em A linear
  domain decomposition method for partially saturated flow in porous media},
  Comput. Methods Appl. Mech. Engrg., 333 (2018), pp.~331--355,
  \href{http://dx.doi.org/10.1016/j.cma.2018.01.029}{doi:\nolinkurl{10.1016/j.cma.2018.01.029}},
  \url{https://doi.org/10.1016/j.cma.2018.01.029}.

\bibitem{seus2017lineartwophase}
{\sc D.~Seus, F.~A. Radu, and C.~Rohde}, {\em A linear domain decomposition
  method for two-phase flow in porous media}, in Numerical Mathematics and
  Advanced Applications ENUMATH 2017, F.~A. Radu, K.~Kumar, I.~Berre, J.~M.
  Nordbotten, and I.~S. Pop, eds., Cham, 2019, Springer International
  Publishing, pp.~603--614.

\bibitem{van1995effect}
{\sc C.~J. {Van Duijn}, J.~Molenaar, and M.~J. {De Neef}}, {\em The effect of
  capillary forces on immiscible two-phase flow in heterogeneous porous media},
  Transport in Porous Media, 21 (1995), pp.~71--93,
  \href{http://dx.doi.org/10.1007/BF00615335}{doi:\nolinkurl{10.1007/BF00615335}},
  \url{https://doi.org/10.1007/BF00615335}.

\bibitem{yotov2001interface}
{\sc I.~Yotov}, {\em Interface solvers and preconditioners of domain
  decomposition type for multiphase flow in multiblock porous media}, in
  Scientific computing and applications ({K}ananaskis, {AB}, 2000), vol.~7 of
  Adv. Comput. Theory Pract., Nova Sci. Publ., Huntington, NY, 2001,
  pp.~157--167.

\end{thebibliography}
\end{document}